\documentclass[opre,nonblindrev]{informs3}

\DoubleSpacedXI 

\usepackage{endnotes}
\let\footnote=\endnote
\let\enotesize=\normalsize
\def\notesname{Endnotes}%
\def\makeenmark{\hbox to1.275em{\theenmark.\enskip\hss}}
\def\enoteformat{\rightskip0pt\leftskip0pt\parindent=1.275em
  \leavevmode\llap{\makeenmark}}

\usepackage{natbib}
 \bibpunct[, ]{(}{)}{,}{a}{}{,}%
 \def\bibfont{\small}%

\usepackage[normalem]{ulem}

\TheoremsNumberedThrough     
\ECRepeatTheorems

\EquationsNumberedThrough    

\MANUSCRIPTNO{OPRE-2022-09-504} 

\usepackage[ruled]{algorithm2e}
\usepackage{verbatim}
\usepackage{enumitem}
\usepackage{bbm}
\usepackage{tikz}
\usetikzlibrary{positioning}
\usepackage{graphicx}
\usepackage{amsmath}
\usepackage{bm}
\usepackage{esvect}

\newcommand{\floor}[1]{\lfloor #1 \rfloor}
\newcommand{\ceil}[1]{\lceil #1 \rceil}
\usepackage{pgfplots}
\pgfplotsset{compat=1.18}
\usepackage{caption}

\usepackage{subcaption}
\usepackage{mathrsfs}  
\usepackage{float}
\usepackage{esvect}
\usepackage{amssymb}

\newcommand*{\QEDB}{\null\nobreak\hfill\ensuremath{\square}}

\newcommand\independent{\protect\mathpalette{\protect\independenT}{\perp}}
\def\independenT#1#2{\mathrel{\rlap{$#1#2$}\mkern2mu{#1#2}}}

\newcommand{\warehouse}[2]{
  \draw[fill=blue!40] (#1) rectangle ++(0.6, 0.6);
  \draw[fill=blue!40] (#1) ++ (0, 0.6) -- ++(0.3, 0.3) -- ++(0.3, -0.3);
}

\usepackage{xr}
\DeclareSymbolFont{extraup}{U}{zavm}{m}{n}
\DeclareMathSymbol{\varheart}{\mathalpha}{extraup}{86}
\DeclareMathSymbol{\vardiamond}{\mathalpha}{extraup}{87}

\newcommand{\vect}[1]{\boldsymbol{\mathbf{#1}}}

\newcommand{\new}[1]{{\color{black} #1}}

\newcommand{\simulation}[1]{{\color{black} #1}}

\newcommand{\yk}[1]{{\color{violet} [YK: #1]}}
\newcommand{\newob}[1]{{\color{magenta} [OB: #1]}}
\newcommand{\ob}[1]{{\color{magenta} #1}}
\newcommand{\ak}[1]{{\color{red} [AK: #1]}}

\usepackage{comment}
\usepackage{color-edits}

\addauthor{ob}{orange}

\begin{document}


\RUNAUTHOR{Besbes, Kanoria and Kumar}

\RUNTITLE{Dynamic Resource Allocation}

\TITLE{Dynamic Resource Allocation: Algorithmic Design Principles and Spectrum of Achievable Performances}

\ARTICLEAUTHORS{%
\AUTHOR{Omar Besbes}
\AFF{Columbia University, Graduate School of Business, New York, NY, 10027, \EMAIL{ob2105@columbia.edu}} 
\AUTHOR{Yash Kanoria}
\AFF{Columbia University, Graduate School of Business, New York, NY, 10027, \EMAIL{ykanoria@gmail.com}}
\AUTHOR{Akshit Kumar}
\AFF{Columbia University, Graduate School of Business, New York, NY, 10027, \EMAIL{ak4599@columbia.edu}}
}

\ABSTRACT{
    {Dynamic resource allocation problems are ubiquitous, arising in inventory management, order fulfillment, online advertising, and other applications. 
    We initially focus on one of the simplest models of online resource allocation: {the multisecretary problem}. In the multisecretary problem, a decision maker sequentially hires up to $B$ out of $T$ candidates, and candidate ability values are drawn i.i.d. from a distribution $F$ on $[0,1]$. First, we investigate fundamental limits on performance as a function of the value distribution under consideration. We quantify performance in terms of \emph{regret}, defined as the additive loss relative to the best performance achievable in hindsight. We present a novel fundamental regret lower bound scaling of $\Omega(T^{\frac{1}{2} - \frac{1}{2(1 + \beta)}})$ for distributions with gaps in their support, with $\beta$ quantifying the mass accumulation of types (values) around these gaps. This lower bound contrasts with the constant and logarithmic regret guarantees shown to be achievable in prior work, under specific assumptions on the value distribution. 
    Second, we introduce a novel algorithmic principle, Conservativeness with respect to Gaps ({\sf CwG}), which yields near-optimal performance with regret scaling of $\tilde{\mathcal{O}}(T^{\frac{1}{2} - \frac{1}{2(1 + \beta)}})$ for any distribution in a class parameterized by the mass accummulation parameter $\beta$. We then turn to operationalizing the {\sf CwG} principle across dynamic resource allocation problems. We study a general and practical algorithm, {\sf Repeatedly Act using Multiple Simulations (RAMS)}, which simulates possible futures to estimate a hindsight-based approximation of the value-to-go function. We establish that this algorithm inherits theoretical performance guarantees of algorithms tailored to the distribution of resource requests, including our {\sf CwG}-based algorithm, and find that it outperforms them in numerical experiments.}
}

\HISTORY{An earlier version of this paper appeared as an extended abstract in the Proceedings of the 23rd ACM Conference on Economics and Computation, EC'22 with the title ``The Multi-secretary Problem with Many Types".}

\maketitle

\section{Introduction}
\label{sec:intro}
\new{Online resource allocation provides a comprehensive framework for scenarios that involve allocating finite resources to requests arriving over time, with the objective of maximizing the overall reward. This model encompasses several well-studied problems, such as the multisecretary problem \citep{arlotto2019uniformly, bray2022logarithmic}, network revenue management \citep{talluri2006theory, vera2021bayesian, bumpensanti2020re}, and order fulfillment \citep{jasin2015lp}. 

Prior work mainly explores these problems under one of two distributional assumptions on the request types: (i) atomic distributions supported on a few points \citep{vera2021bayesian, bumpensanti2020re} and (ii) non-atomic distributions with contiguous support \citep{lueker1998average, bray2022logarithmic}. Under these cases, impressive constant and logarithmic regret guarantees have been established, where regret is defined as the expected difference between the total reward under the optimal hindsight policy and the total reward gathered under an online policy.}

\new{
However, for many important applications, neither of these two assumptions adequately capture reality. For instance, consider the order fulfillment problem encountered by e-commerce platforms like Amazon or Walmart. This is an online matching problem with spatially distributed demand (different zip codes or counties) with product inventory housed in various warehouses scattered across a geographic area. The fulfillment team aims to minimize cumulative shipping costs by dynamically matching each demand to a warehouse which has the item available. Warehouses have limited inventory, and decisions must be made in real-time. This problem can be framed within the online resource allocation problem paradigm. Yet, the aforementioned assumptions made in the prior literature do not capture key features of this setting: (i) the number of demand locations (types) is large (for instance, there are over 40,000 zip codes in the United States), and (ii) these demand locations are spatially clustered with gaps (regions with no demand), a natural characteristic of geographical landscapes such as rivers, mountains, deserts, etc. Hence, atomic distributions with a low number of types or  non-atomic distributions with contiguous support fail to capture the salient features of such a problem. Aside from modeling concerns, the near-optimal algorithms developed for each of the two classes of distributions mentioned above are tailored to that particular class of distributions. 
 
The above motivation leads us to the following research questions: \emph{(i) What (request type) distribution features drive achievable performance, and how does regret scale as a function of the underlying distribution? (ii) What algorithmic principles allow one to achieve optimal regret scaling? (iii) Is there a unifying near-optimal algorithm that is agnostic to the underlying distribution's features?}

To isolate and examine key performance drivers, we will initially focus on one of the simplest online resource allocation problems: the multisecretary problem, which is a special case of both the network revenue management problem as well as the online matching (order fulfillment) problem (we refer to Appendix \ref{onapp:fulfillment-problem-mapping} for a more extensive discussion on the latter connection).  In the multisecretary problem, a decision-maker (DM) with a budget to hire $B$ secretaries is presented with a series of $T$ independent values representing candidate abilities. The DM must make irrevocable ``accept'' (i.e., hire) or ``reject'' decisions on the fly, aiming to maximize the (expected) sum of the chosen candidates' abilities. 

  We make three main contributions.  The first two are in the context of the multisecretary problem: fundamental lower bounds on regret, and an algorithmic principle to achieve the optimal regret scaling. Our third contribution is a unifying and practical algorithm for achieving near optimal regret performance in general resource allocation problems. We now elaborate on these contributions.}
\begin{enumerate}[label = \emph{(\roman*)}]
    \item \textit{Drivers of regret}: In the context of the multisecretary problem, we identify a novel fundamental driver of regret which is characterized by a parameter $\beta$, which quantifies the mass accumulation of types around gaps (interval with zero probability mass). Using this parameter $\beta$ we characterize a broad class of distributions with gaps, which we refer to as $(\beta, \varepsilon_{0}, \delta)$-clustered distributions (cf. Definition \ref{def:beta-epsilon-delta-clustered}). The class of $(\beta, \varepsilon_{0}, \delta)$-clustered distributions is a superset of the class of discrete distributions \citep{arlotto2019uniformly}, and the class of non-atomic distributions with continuous support over $[0,1]$ and density uniformly bounded away from zero \citep{bray2022logarithmic}.
    We establish a universal lower bound (for any policy) on the growth rate of the regret as a function of the parameter $\beta$ which quantifies how mass accumulates around gaps. In particular, we establish that any policy must incur $\Omega(T^{\frac{1}{2} - \frac{1}{2(1 + \beta)}})$ regret in the worst-case (cf. Theorem \ref{thm:CwG-regret-lower-bound}) for a $(\beta, \varepsilon_{0}, \delta)$-clustered distribution. This is in stark contrast to prior results which prove regret scaling of $\Theta(1)$ \citep{arlotto2019uniformly} for the case of  distributions with a few discrete types and $\Theta(\log T)$ \citep{bray2022logarithmic} for a special class of non-atomic distributions. We also show that our lower bound on the regret scaling is achievable up to polylogarithmic factors. To the best of our knowledge, ours is the first result of its kind; notably the regret scaling we establish is polynomial in $T$ for $\beta > 0$ and an entire spectrum of regret scalings are possible. As $\beta$ increases, so does the exponent $\frac{1}{2} - \frac{1}{2(1 + \beta)}$ (from 0 to 1/2), characterizing the ``hardness'' of the problem instance.

    \item \textit{Algorithmic Principle}: It turns out the workhorse certainty equivalent ({\sf CE}) policy is insufficient to deal with general type distributions which have gaps in the support, already in  the case of the multisecretary problem. For such distributions, we introduce a new algorithmic principle we call \emph{Conservativeness with respect to gaps (CwG)}; which makes a crucial modification to the {\sf CE} policy. The idea is that if at any time the {\sf CE} threshold is close to the boundary of a gap, {\sf CwG} instead uses the gap as the acceptance threshold to avoid incurring large regret in the future. We establish that this enables the policy to mitigate the risk of incurring large regret (in the event that the threshold for the hindsight optimal falls on the opposite side of that gap). We use this principle to design a near-optimal algorithm, dubbed \textsf{CwG}, for the $(\beta, \varepsilon_{0}, \delta)$-clustered distributions. Its worst-case regret scales as $\tilde{\cal O}(T^{\frac{1}{2} -  \frac{1}{2(\beta + 1)}})$, matching the scaling of the lower bound in $T$ up to polylogarithmic terms (cf. Theorem \ref{thm:general-cwg-upper-bound}). For the case of a few discrete types, our algorithm recovers bounded regret, as in \cite{arlotto2019uniformly} (cf. Corollary \ref{cor:discrete-distribution-bound}). For the special class of non-atomic distributions with density bounded away from zero, {\sf CwG} is identical to {\sf CE} since there are no gaps and we recover the logarithmic regret scaling result of \cite{lueker1998average} and \cite{bray2022logarithmic} (cf. Corollary \ref{cor:cr-upper-bound-example}).  

    \item \simulation{\textit{Unifying Algorithm}: Returning to general resource allocation problems, we propose a versatile algorithm called \textsf{Repeatedly Act using Multiple Simulations (RAMS)}, which offers a practical and data-driven approach to resource allocation. At each $t$, {\sf RAMS} simulates multiple future demand scenarios.
    Each possible allocation decision at $t$ results in different cumulative rewards in hindsight, in each demand scenario. {\sf RAMS} greedily selects the allocation decision which maximizes the average over scenarios of the cumulative reward in hindsight.
    Unlike previous algorithms, {\sf RAMS} does not require to be tuned to specific distribution features, and 
    by its design can organically leverage the data-driven simulations of the future which are typically available in practical applications.
    In terms of performance, we establish a meta result (Theorem \ref{thm:meta-performance-rams}) that shows that {\sf RAMS} is guaranteed to inherit the regret performance guarantee of any algorithm satisfying certain conditions (specified in Theorem \ref{thm:meta-performance-rams}).
    This result, in conjunction with Theorem \ref{thm:general-cwg-upper-bound}, implies that {\sf RAMS} is near-optimal for the multisecretary problem and naturally incorporates the \emph{conservativeness with respect to gaps} principle. Furthermore, our meta theorem, together with existing results on other algorithms in the literature, tells us that {\sf RAMS} is near-optimal in a variety of settings for NRM and Order Fulfillment problems. 
    }
\end{enumerate}

\subsection{Related Literature}
The classical secretary problem was introduced by \cite{cayley1875mathematical} and \cite{moser1956problem}. 
The multisecretary variant of the above problem was initially studied by \citet{kleywegt1998dynamic} and \citet{kleinberg2005multiple}. Recently, \citet{arlotto2019uniformly} showed that, when the distribution of types is discrete, regret is bounded uniformly for all values of the number of candidates $T$ and the hiring budget $B$, where the constant may scale with the reciprocal of the minimum probability mass on any type.  In order to prove this result, they devise an adaptive policy called the Budget-Ratio (BR) policy where they compare the ratio of the remaining budget to the remaining number of candidates to interview and make the hire/reject decision by comparing the budget ratio to some fixed thresholds. This regret guarantee, in conjunction with a lower bound on regret from \citet{kleinberg2005multiple} yields a tight understanding of the class of distributions supported on a few discrete types. \new{Note that the classical secretary problem and its generalization considered in \cite{kleinberg2005multiple} do not assume the knowledge of the reward distribution. 
 However following the  work of \cite{arlotto2019uniformly}, the variant of multisecretary with distributional knowledge has also been referred to as the multisecretary problem and we will also employ this terminology.} 

At the other extreme, for a continuum of types, \cite{lueker1998average,bray2022logarithmic} show that instead of the regret being uniformly bounded, the best possible scaling for a certain class of non-atomic distributions with contiguous support is $\Theta(\log T)$ (\cite{bray2022logarithmic} shows that this is true for the more general network revenue management problem as well). In the context of the multisecretary problem, they devise a simple threshold policy based on the budget ratio to achieve this regret scaling. However, the class of non-atomic distributions considered in these papers requires the probability density function to be bounded away from zero. \new{In a parallel line of inquiry, the set of distributions examined by Blumrosen and Holenstein (2008) bears close resemblance to our own. Yet, there are marked differences in the settings and results. Specifically, Blumrosen and Holenstein (2008) concentrate on the auction setting involving a single item and restrict their study to continuous distributions.}

The multisecretary problem is a special case of a broader class of network revenue management (NRM) problems, or more broadly dynamic resource constrained reward collection  problems; see  \cite{balseiro2021survey} for a recent survey and unified modeling framework for this class of problems.  There is a wide variety of applications in auction theory \citep{kleinberg2005multiple}, online resource allocation \citep{kleywegt1998dynamic, talluri2006theory}, order fulfillment \citep{jasin2015lp}, among others. Note that this literature typically assumes a small number of types. 

\cite{vera2021bayesian,vera2021online} generalized the arguments in \cite{arlotto2019uniformly} to a broader class of online packing and online matching problems and proved a uniform regret guarantee across all values of capacity $B$ and time horizon $T$. They developed a technique called \emph{compensated coupling} and used it to prove a constant regret guarantee without requiring any non-degeneracy assumptions. \cite{bumpensanti2020re} also proved constant regret guarantees for a class of NRM problems, however their algorithm and proof techniques differ from those of \cite{vera2021bayesian, vera2021online}. While all these papers impressively establish constant regret bounds, all of them assume a few discrete types, and their regret bounds scale polynomially in the number of types. However in many practical systems, the number of types is, in fact, large. 

\simulation{Simulation-based algorithms have been studied in the network revenue management literature \citep{talluri1999randomized, kunnumkal2012randomized}, albeit without any regret guarantees. The idea in these papers is to solve multiple stochastic optimization problems with different realizations instead of a single fluid relaxation and average the shadow prices of the different optimization problems and implement a bid-price control. 
Recently, \citet{freund2019good} and \citet{sinclair2022hindsight} have used related ideas to develop algorithms for online bin packing with 
a few types.}

Another line of research connected to our work is on prophet inequalities, in particular $k$-unit prophet inequalities ($k$ corresponds to the budget $B$ described earlier). The $k$-unit prophet inequality problem, originally studied in \cite{hajiaghayi2007automated},  analyzes the competitive ratio which is defined as the ratio of the expected performance of an algorithm to the expected performance of the hindsight optimal in the worst case over the reward distributions, where the focus is on deriving tight guarantees in terms of $k$. The seminal work of \cite{alaei2014bayesian} proved a guarantee of $1 - 1/\sqrt{k + 3}$ on the competitive ratio and since then this result has been improved upon by  \cite{chawla2020static} and \cite{jiang2022tight}. One key distinction between this stream and our work is that we consider i.i.d values from a known distribution, which allows to prove stronger guarantees on the regret. The competitive ratio results above would imply a regret scaling of $\Theta(\sqrt{T})$, whereas we show that if the distribution is known and i.i.d, it is possible to do better even under the worst-case when the budget $B$ scales linearly in $T$ (cf. Theorem \ref{thm:general-cwg-upper-bound}).

\textbf{Organization of the paper.} Section \ref{sec:model} describes the model. In Section \ref{sec:fundamental-limits}, we describe a general family of distributions, dubbed $(\beta, \varepsilon_{0}, \delta)$-clustered distributions, and provide novel fundamental limits on regret scaling. In Section \ref{sec:cwg-principle} we state our key  \emph{conservativeness with respect to gaps (CwG)} algorithmic principle and provide near-optimal regret scaling for $(\beta, \varepsilon_{0}, \delta)$-clustered distributions in the context of the multisecretary problem. In Section \ref{sec:rams}, we discuss our unifying algorithm \textsf{RAMS}. We conclude in Section \ref{sec:conclusion}. Due to space constraints, all proofs have been relegated to the appendix.

\section{Model}
\label{sec:model}
\new{We consider a dynamic resource allocation problem with a {\it known} finite time horizon $T$. 
There are $d$ resources and the decision maker is endowed with an initial budget vector $B \in \mathbb{R}^d$ for the resources. At each time $t = 1, 2, \dots, T$, a request $\theta_t$ is drawn independently from a \emph{type} set $\Theta$ via some distribution $F$ which is \emph{known} to the decision maker. Upon observing a request $\theta_t$, the decision maker takes an action $a_t \in \mathcal{A}(B_t, \theta_t)$ where $\mathcal{A}(B_t, \theta_t)$ is the set of feasible actions at time $t$ which depends on the remaining budget $B_t$ and the request $\theta_t$. Let $\mathcal{A} \triangleq \cup_{B \geq \bm{0}} \cup_{\theta \in \Theta} \mathcal{A}(B, \theta)$ denote the set of all possible actions. Upon taking an action $a_t$, the decision maker collects a reward $r_t$ which depends on the request $\theta_t$ and the action $a_t$. We denote by $r:\Theta \times \mathcal{A} \to \mathbb{R}_{\geq 0}$  the reward function. Taking an action consumes resources and the amount of resource consumed depends on the request $\theta$, and is denoted by a consumption function $c : \Theta \times \mathcal{A} \to \mathbb{R}^d$ where $c_k(\theta, a)$ is the amount of $k$-th resource consumed when the request is $\theta$ and action is $a$. Given a request $\theta_{t}$ and action $a_t$, the remaining budget is updated as per $B_{t + 1} = B_{t} - c(\theta_t, a_t)$; the action $a_t$ is required to be such that each coordinate of $B_{t+1}$ is non-negative. We assume that there is a null action $a_0 \in \mathcal{A}$ which consumes no resources and generates no reward, i.e., $r(\theta, a_0) = 0$ for all $\theta \in \Theta$ and $c(\theta, a_0) = 0_{d \times 1}$ for all $\theta \in \Theta$. Further, we will assume that $|r(\theta, a)| \leq 1$ and $\|c(\theta, a)\|_\infty \leq 1$ for all $\theta \in \Theta$ and $a \in \mathcal{A}$.
}

A policy is said to be an \emph{online} (non-anticipating) policy if the decision on the $t$-th request is based only on the request $\theta_t$ at time $t$, the past requests,  $\{\theta_{j}\}_{j = 1}^{t - 1}$ and the history of the actions $\{a_j\}_{j = 1}^{t - 1}$ up to the time $t$. Let $U_{1}, U_{2}, \dots, U_{T}$ be a sequence of random variables that are independent and uniformly distributed over $[0,1]$ and independent of the requests $\theta_{1}, \theta_{2}, \dots, \theta_{T}$.  (The $U$s will allow us to accommodate randomized policies.) Define the filtration $\mathcal{F}_{t} = \sigma(\theta_{1}, U_{1}, \theta_{2}, U_{2}, \dots, \theta_{t}, U_{t} )$ for all $t \in [T]$. A feasible online policy ${\pi}$ is a sequence of $\{\mathcal{F}_{t}: t \in [T]\}$-measurable random variables $\{a^{\pi}_{1}, a^{\pi}_{2}, \dots, a^{\pi}_{T}\}$ such that $\sum_{t = 1}^{T} c(\theta_t, a_t^{\pi}) \leq B$ almost surely. We define the set of feasible online policies as $\Pi(B,T)$. For any feasible and online policy $\pi \in \Pi(B,T)$, define $R_{t}^{\pi} = \sum_{k = 1}^{t} r(\theta_k, a_k^{\pi}), \forall t \in [T]$ to be the accumulated reward up to  time $t$. The total expected reward under a policy $\pi \in \Pi(B,T)$ is given by $V^{\pi}_{1}(B, T) = \mathbb{E}\left[R_{T}^{\pi} \right] = \mathbb{E}\left[\sum_{t = 1}^{T} r(\theta_t, a_t^{\pi}) \right]$. 
Fix $T \in \mathbb{N}$ and $B \in \mathbb{R}_{\geq 0}^d$, the objective is to maximize the total expected reward given by $V_1^{\star}(B, T) = \sup_{\pi \in \Pi(B,T)} V_{1}^{\pi}(B,T)$.

Next we consider the hindsight ({\sf hs}), full-information version of the problem in which the requests $\vect{\theta}_{\geq 1} = \{\theta_{1}, \theta_{2}, \dots, \theta_{T}\}$ are known apriori. In the hindsight  setting, the problem essentially reduces to solving  $V_{1}^{\sf hs}\left(B, T;\vect{\theta}_{\geq 1}\right) = \max_{\bm{a}}\{\sum_{t = 1}^{T} r(\theta_t,a_t): (a_{1}, a_{2}, \dots, a_{T}) \in |\mathcal{A}|^{T} \text{ and } \sum_{t = 1}^{T} c(\theta_t, a_t) \leq B\}$ and the total expected value by the hindsight optimal problem is given as $V_1^{\sf hs}(B,T) = \mathbb{E}\left[V_1^{\sf hs}\left(B, T; \vect{\theta}_{\geq 1}\right) \right]$. It trivially follows that $V^{\sf hs}_1(B, T) \geq V_{1}^{\pi}(B, T), \forall \pi \in \Pi(B,T)$ and $\forall B \in [T]$. To measure the performance of a feasible online policy $\pi \in \Pi(B,T)$, we consider the hindsight problem as a benchmark and define the (expected) regret of the policy ${\pi}$ as the difference between the expected value of the hindsight problem and the expected value attained by the policy ${\pi}$ i.e.,
    $\text{Regret}(B,T; {\pi}) \triangleq V^{\textsf{hs}}_1(B,T) - V_{1}^{{\pi}}(B,T).$
We also define the (minimum achievable, expected) regret as the difference between the expected value of the hindsight problem and the expected value under the optimal online policy $\pi^\star \in \Pi(B,T)$.
\begin{align*}
    \text{Regret}(B,T) &= \inf_{{\pi} \in \Pi(B,T)} \text{Regret}(B,T; {\pi}) = V^{\textsf{hs}}_{1}(B,T) - V_1^{\star}(B,T).
\end{align*}
In what follows, we will focus on characterizing the growth rate of  $\text{Regret}(B,T)$ as a function of $T$ and the characteristics of the  distribution of types. 
\new{
Next we  discuss the three important classes of online resource allocation problems.
\paragraph{Network Revenue Management.} In this problem each request $\theta = (r_{\theta},\bm{c}_{\theta})$ is presented with a single reward $r_{\theta} \geq 0$ and a consumption vector $\bm{c}_{\theta} \in \mathbb{R}^d$. We have that $\mathcal{A} = \{a_0 = \text{reject}, a_1 = \text{accept} \}$. The reward and consumption functions are given as 
\begin{align*}
    &r(\theta, \text{reject}) = 0, \ \ c(\theta, \text{reject}) = \bm{0}_{d \times 1} \\
    &r(\theta, \text{accept}) = r_\theta, \ \ c(\theta, \text{accept}) = \bm{c}_\theta. 
\end{align*}

\paragraph{Online Matching (Order Fulfillment).} In this problem each request $\theta = r_\theta$ is presented with a vector of rewards $\bm{r}_{\theta} \in \mathbb{R}^{d}$.
Each request wants to consume at most one unit of any single resource. The action set is $\mathcal{A} = \{a_0, a_1, \dots, a_d\}$ where $a_k$ denotes that the request is matched to resource $k$ with $a_0$ being the null action denoting that the request is rejected. The reward and consumption functions are given as
\begin{align*}
    &r(\theta, a_0) = 0 , \ \ c(\theta, a_0) = \bm{0}_{d \times 1} \\
    &r(\theta, a_k) = \bm{r}_{\theta, k}, \ \ c(\theta, a_k) = \bm{e}_k, \qquad \forall k \in \{1,2,\dots, d\}
\end{align*}
where $\bm{r}_{\theta, k}$ denotes the $k$-th coordinate of $\bm{r}_{\theta}$ and $\bm{e}_k$ is a $d$-dimensional vector with the $k$-th coordinate being one and all other coordinates being zero.

\paragraph{Multisecretary Problem.} For the case of one resource ($d = 1$), network revenue management and online matching are equivalent problems and this special case is referred to as the multisecretary problem. We have that 
$c(\theta, \textup{accept})=1$ for all $\theta \in \Theta$. In the context of the multisecretary problem, the request type (equivalently, reward) will be referred to as the candidate ability. 
}

\section{Fundamental Limits on Achievable Performance}
\label{sec:fundamental-limits}
\new{To delve deeper into the intrinsic drivers of performance, we initially focus on the multisecretary problem -- a cornerstone model in online resource allocation. Clearly, any lower bound established for the multisecretary problem directly translates into a lower bound for a broader range of online resource allocation problems like NRM and online matching. We now define two classes of distributions under which the multisecretary problem has been previously studied.
\begin{assumption}[Small Number of Types]
    \label{ass:ms-finite-types}
    The type (reward) distribution $F$ is supported on a finite set 
    and the rewards are assumed to be in the interval $[0,1]$.
\end{assumption}
\begin{remark}
    Many prior works refer to this as the ``finite types setting'', and establish constant regret guarantees \citep[see, e.g.,][]{arlotto2019uniformly, bumpensanti2020re, vera2021bayesian}. However, these guarantees scale linearly with the number of types. Hence, they are most relevant when the size of discrete types set is {\it small}. To emphasize this aspect, we use the phrases ``small number of types'' or ``small discrete set'' or ``few types'' to describe this setting.
\end{remark}

\begin{assumption}[Infinitely Many Types with density bounded away from zero]
    \label{ass:ms-infinite-types}
    The type (reward) distribution $F$ is supported on an infinite set and $F$ admits a density $f$ which is bounded from below and above, i.e., there exist $0 < \underline{\nu} \leq \bar{\nu}< \infty$ such that $\underline{\nu} \leq f(\theta) \leq \bar{\nu}$ for all $\theta \in \Theta$. The rewards are assumed to be in the interval $[0,1]$.
\end{assumption}

To interpolate between these two class of distributions, we will introduce a general class of distributions which will capture the distributions with a few types and infinitely many types with bounded density as special cases.
}

\subsection{General Class of Distributions For the multisecretary Problem}
\label{sec:beta-epsilon-distribution}

We will anchor our analysis around a general family of distributions which allow for gaps in the type space and can capture as special cases discrete distributions as well as the non-atomic distributions with density uniformly bounded away from zero. We call this family $(\beta, \varepsilon_{0}, \delta)$-clustered distributions. 
For any $q \in [0,1]$, we define $F^{-1}(q) \triangleq \inf \{v: F(v) {\ge} q\}$. 



\begin{definition}[$(\beta, \varepsilon_{0}, \delta)$-clustered distributions]
\label{def:beta-epsilon-delta-clustered}
Fix $\beta \in [0, \infty)$, $\varepsilon_{0} \in (0,1]$ and $\delta \in [0,1]$. A distribution $F$ is said to be $(\beta, \varepsilon_{0},\delta)$-clustered if there exists $n \in \mathbb{N}\cup \{ 0 \}$ and \emph{gap quantiles}
$q^\star_0 = 0 < q^\star_1 < \dots < q^\star_n < q^\star_{n+1}=1$ such that we have 
\begin{enumerate}[label = (\alph*)]
    \item (Generalized cluster ``density'' requirement) $\forall i \in [n + 1], \forall q, \tilde{q} \in (q_{i-1}^\star, q_{i}^\star ]$, we have that $|F^{-1}(q) - F^{-1}(\tilde{q})| \leq C|q - \tilde{q}|^{\frac{1}{\beta + 1}} + \delta$ for some constant $C < \infty$.
    \item (Cluster size requirement) $q_{i}^\star - q_{i-1}^\star \geq \varepsilon_{0}, \forall i \in [n + 1]$.
\end{enumerate}
\end{definition}

Let $\mathscr{F}_{\beta, \varepsilon_{0}, \delta}$ denote the class of $(\beta, \varepsilon_{0}, \delta)$-clustered distributions. 
This class  includes a wide variety of distributions. An important sub-class  is the one with $\delta = 0$, which we denote by $\mathscr{F}_{\beta, \varepsilon_0}$.  We refer to  distributions in this subclass as $(\beta, \varepsilon_0)$-clustered. 


Define $H_i \triangleq [F^{-1}((q_{i-1}^\star)^+ ), F^{-1}(q_{i}^\star)]$, for all $i \in [n+1]$, where $F^{-1}(q^+) \triangleq \lim_{\epsilon \to 0^+} F^{-1}(q+\epsilon)$. 
We will refer to the $(H_i)$'s as {\it mass clusters} or just {\it clusters}.
We will use the term \emph{gaps} to refer to the complementary intervals $G_i \triangleq (F^{-1}(q^\star_i), F^{-1}((q^\star_i)^+))$ for $i \in [n]$, and the intervals at the extremes $G_0 = [0, F^{-1}(0^+)), G_{n+1}= (F^{-1}(1), 1]$, since they contain no probability mass. The requirement (a) 
can be thought of as a within-cluster ``density'' requirement, which becomes weaker as $\beta$ increases; we can think of $\beta$ as quantifying the within-cluster mass density (with a decreasing relationship). When $\delta=0$, this requirement corresponds to $F^{-1}$ being $(1/(\beta+1))$-H\"{o}lder continuous on the mass clusters. Requirement (b) is a cluster size requirement, $\varepsilon_0$ being the minimum cluster size; this requirement becomes more stringent as $\varepsilon_0$ increases. 
\new{The parameter $\delta$ provides us with additional flexibility in modelling our distributions. One such practically relevant class of distributions is the one with a large number of discrete types, which can be modelled using the parameter $\delta$ (cf. Example \ref{ex:many-small-discrete-types}).
In general, there is some flexibility on how the distributions are modelled, more specifically how the types are aggregated into clusters, and this is associated with a tradeoff between $\delta$ and $\varepsilon_0$ (and potentially $\beta$). Please refer to Appendix \ref{onapp:details-beta-varepsilon-delta-clustered} for more details.}

Next we present some examples of $(\beta, \varepsilon_{0}, \delta)$-clustered distributions including discrete distributions, as well  the uniform distribution, along with the appropriate choices of gap quantiles.
\begin{example}[Discrete Distributions]
\label{ex:discret-distributions}
 Consider a discrete distribution  \citep[as studied in][]{arlotto2019uniformly}.
 Let the support be $\{\theta_{1}, \theta_{2}, \dots, \theta_{n}\}$ with probability masses $\{p_{1}, p_{2}, \dots, p_{n}\}$. 
Assume that $0\leq \theta_{1} < \theta_2 < \dots < \theta_{n} \leq 1$. We make use of the natural choice of gap quantiles $q_i^\star = \sum_{j=1}^i p_j$ for all $i \in [n-1]$, leading to gaps $G_0=[0, a_1), G_i = (a_i, a_{i+1}) \; \forall i \in [m-1], G_n = (a_n, 1]$ and clusters $H_i = \{a_i\} \; \forall i \in [n]$. Now for $q, \tilde{q} \in (q_{i-1}^\star, q_i^\star] = Q_i$, we have that $|F^{-1}(q) - F^{-1}(\tilde{q})| = 0 \leq |q - \tilde{q}|$, i.e., the cluster density requirement is satisfied for $\beta = 0$ and $\delta = 0$. Defining $\varepsilon_{0} \triangleq \min\{p_{1}, p_{2}, \dots, p_{n}\}$ the cluster size requirement is satisfied. Therefore the discrete distribution belongs to the class of $(0, \varepsilon_{0})$-clustered distributions where $\varepsilon_{0}$ is the minimum probability mass in the support. 
\end{example}

\begin{example}[Non-atomic Distributions with Contiguous Support]
    \label{ex:bray}
 Consider the non-atomic distributions with pdf $f$ considered in \cite{bray2022logarithmic} (Assumption \ref{ass:ms-infinite-types}). Assume that there exists $\alpha_{0} > 0$ such that $f(x) \geq \alpha_{0}, \forall x \in [0,1]$. 
 (The uniform distribution over $[0,1]$ is a special case of these distributions with $f(x) = 1$ for all $x \in [0,1]$.) Such distributions are $(\beta = 0, \varepsilon_0 = 1, \delta = 0)$-clustered distributions with $n = 0$ gaps, i.e., $F^{-1}$ is $1$-H\"{o}lder continuous over the interval $(0,1]$ with the constant $C = 1/\alpha_0$.  The gap quantiles are only the trivial ones $q_0^\star = 0$ and $q_{1}^\star = 1$. There is a single mass cluster $H_1= [0,1]$ with mass $1$, which clearly satisfies the cluster density requirement with $\beta = 0, \varepsilon_0 = 1$ and $\delta = 0$.
\end{example}

\begin{figure}[ht]
    \centering
        \pgfplotsset{width = 0.48\linewidth, ylabel near ticks}
        \begin{tikzpicture}
            \begin{axis}[
                axis lines = left,
                xlabel = \(x\),
                ylabel = {\(f_{\beta}(x)\)},
            ]
            \addplot [domain=0:0.25, color=red, thick]{2};
            \addplot[domain = 0.75:1, color = red, thick, forget plot]{2};
            \addplot[domain = 0.25:0.75, color = red, thick, forget plot]{0};
            \addplot[mark=*,red, fill = white, thick, forget plot] coordinates {(0.25,0)};
            \addplot[mark=*,red, fill = white, thick, forget plot] coordinates {(0.75,0)};
            \addplot[mark=*,red, thick, forget plot] coordinates {(0.25,2)};
            \addplot[mark=*,red, thick, forget plot] coordinates {(0.75,2)};
            \addlegendentry{\(\beta = 0\)}
            \addplot[domain = 0:0.25, color = blue, thick, dashed]{4 - 16*x};
            \addplot[domain = 0.75:1, color = blue, thick, dashed, forget plot]{16*x - 12};
            \addlegendentry{\(\beta = 1\)}
            \addplot[domain = 0:0.25, color = black, thick, dotted]{96 * (x - 0.25)^2};
            \addplot[domain = 0.75:1, color = black, thick, dotted, forget plot]{96 * (x - 0.75) ^ 2};
            \addlegendentry{\(\beta = 2\)}
            \end{axis}
        \end{tikzpicture}
    \hfill
         \begin{tikzpicture}
            \begin{axis}[
                axis lines = left,
                xlabel = \(x\),
                ylabel = {\(F_{\beta}(x)\)},
                legend style={at={(0.72,0.02)},anchor=south west}
            ]

        \addplot [domain=0:0.25, color=red, thick]{2*x};
        \addplot[domain = 0.75:1, color = red, thick, forget plot]{2*x - 1};
        \addplot[domain = 0.25:0.75, color = red, thick, forget plot]{0.5};
        \addlegendentry{\(\beta = 0\)}

        \addplot[domain = 0:0.25, color = blue, thick, dashed]{0.5 - 8 * (x - 0.25) ^ 2};
        \addplot[domain = 0.75:1, color = blue, thick, dashed, forget plot]{0.5 + 8 * (x - 0.75)^ 2};
        \addplot[domain = 0.25:0.75, color = blue, thick, dashed, forget plot]{0.5};
        \addlegendentry{\(\beta = 1\)}

        \addplot[domain = 0:0.25, color = black, thick, dotted]{0.5 + 32 * (x - 0.25)^3};
        \addplot[domain = 0.75:1, color = black, thick, dotted,  forget plot]{0.5 + 32 * (x - 0.75) ^ 3};
        \addplot[domain = 0.25:0.75, color = black, thick, dotted, forget plot]{0.5};
        \addlegendentry{\(\beta = 2\)}
    \end{axis}
    \end{tikzpicture}
    \caption{ (L) PDF $f_{\beta}$ for $\beta = 0,1,2$. Notice the gap from $1/4$ to $3/4$, (R) CDF $F_{\beta}$ for $\beta = 0,1,2$.}
    \label{fig:example-beta-clustered-distribution}
\end{figure}
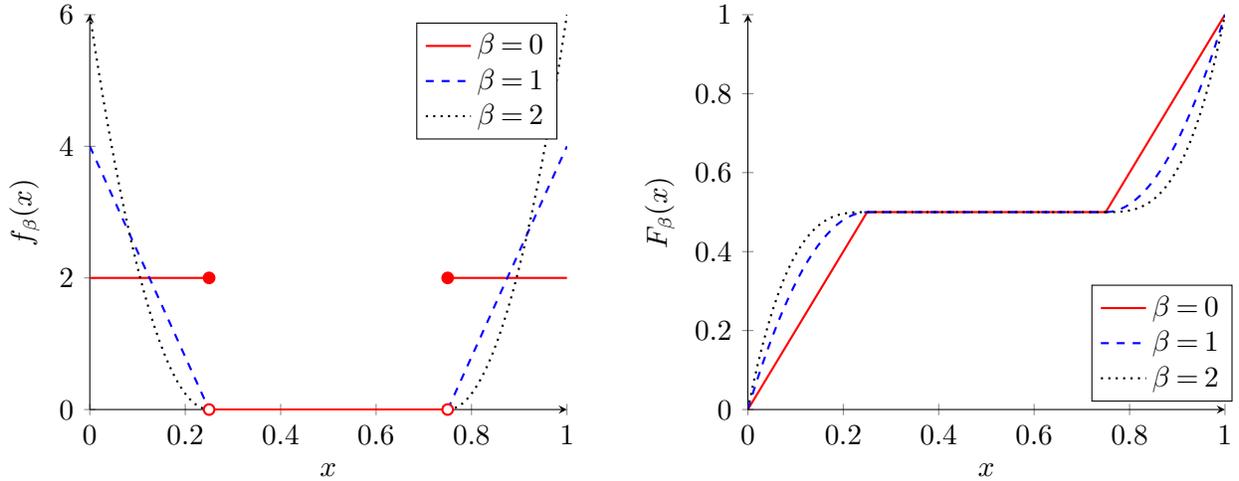

\begin{example}[A class of bimodal distributions]
\label{ex:bimodal-distribution}
 An example of a $(\beta, \varepsilon_{0})$-clustered distribution with $n = 1$ gap (with gap quantile $q_1^\star = 1/2$), for general $\beta \geq 0$ and $\varepsilon_{0} = 1/2$, which we will make use of to prove our lower bound results is presented below:
\begin{align}
    \label{eq:example-beta-epsilon-clustered-distribution}
    F_{\beta}(x) = \begin{cases}-2\cdot 4^\beta \cdot \left(\frac{1}{4} - x \right)^{\beta + 1} + \frac{1}{2}  & \quad 0 \leq x \leq \frac{1}{4} \\ \frac{1}{2} &\quad \frac{1}{4} \leq x \leq \frac{3}{4} \\ 2\cdot 4^\beta \cdot \left( x - \frac{3}{4}\right)^{\beta + 1} + \frac{1}{2} & \quad \frac{3}{4} \leq x \leq 1\end{cases}, \ \ \ \
    F_{\beta}^{-1}\left(q\right) = \begin{cases} \frac{1 - (1 - 2q)^{\frac{1}{\beta + 1}}}{4}, &\quad 0 \leq q \leq \frac{1}{2} \\ 
    \left[\frac{1}{4}, \frac{3}{4}\right] & \quad q = \frac{1}{2} \\
    \frac{(2q - 1)^{\frac{1}{\beta + 1}} + 3}{4}, & \quad \frac{1}{2} < q \leq 1\end{cases}
\end{align}
It is easy to see that $F^{-1}_{\beta}$ in \eqref{eq:example-beta-epsilon-clustered-distribution} is a $(\beta, 1/2)$-clustered distribution, with one gap $G_{1} = (1/4,3/4)$ and clusters $H_1 = [0, 1/4]$ and $H_2 = [3/4,1]$. Refer to Figure~\ref{fig:example-beta-clustered-distribution} for a plot of the density $f_{\beta}$ and the CDF $F_{\beta}$ of the $(\beta, 1/2)$-clustered distribution defined in \eqref{eq:example-beta-epsilon-clustered-distribution}.
\end{example}

Observe that $(\beta, \varepsilon_{0})$-clustered distributions already allow us to capture not only the previously studied distributions such as distributions with few types and continuous distributions (with density bounded below), but also a mixture of atomic and non-atomic distributions with gaps. As mentioned previously, the parameter $\delta$ provides us with additional flexibility to model distributions with a large number of discrete types, which may be of practical relevance. One such example is that of {\it many small discrete types} which we present below.

\begin{example}[Many Small Discrete Types]
    \label{ex:many-small-discrete-types}
    Fix a small $\delta > 0$ and consider a discrete distribution with \emph{many small} types supported on the points $\mathcal{S} = \{0,\delta, 2\delta, \dots, 1/4\} \cup \{3/4, 3/4 + \delta, 3/4 + 2\delta, \dots, 1\}$ with probability mass $2\delta$ on each of the points in $\mathcal{S}$. This constitutes a setting with many small discrete types since there are a large number of atomic types (separated by small empty intervals) and the probability mass of each type is small, i.e., it is proportional to $\delta$. This instance of many small discrete types captures the salient feature of the order fulfillment problem that there are a large number of demand types (e.g., zipcodes) with each demand type having small probability mass and these demand types are spatially clustered with possibly large gaps between different clusters of demand types. 
    As $\delta \to 0$, we recover the bimodal uniform distribution $F_0$ in the limit. One can similarly consider similar many-small-discrete-type analogs for other $(\beta, \varepsilon_0)$-clustered distributions. Note that the many small types need not be uniformly spaced. We require that the maximum distance between the discretized grid points be at most $\delta$.  In such discretizations, we have some flexibility in choosing which empty intervals to classify as ``gaps''. In the case of {many small} types, if the size of the empty intervals (due to discretization) is at most $\delta$ then we may consider the entire clump of these {many small} types as belonging to one cluster (say, $H_{i}$) and hence one quantile interval ($Q_{i}$). 
\end{example}

\subsection{Fundamental Lower bound on Performance}
\label{sec:universal-lower-bound}
\new{In this section we present a novel driver of regret scaling: the shape of the candidate ability (or value) distribution around gaps which is characterized by the parameter $\beta \in [0, \infty)$ and show that for $\beta > 0$, polynomial regret scaling is unavoidable. To focus on the scaling with parameter $\beta$, we fix $\delta = 0$ and $\varepsilon_0 = 1/2$. }
\begin{theorem}[universal lower bound]
\label{thm:CwG-regret-lower-bound}
Fix $\delta = 0, \varepsilon_0 = 1/2$ and consider any $\beta \in [0, \infty)$. Then there exists a candidate-ability distribution $F \in \mathscr{F}_{\beta, \varepsilon_{0}}$, a time horizon $T_{0} < \infty$, a universal constant $c > 0 $ such that, for all $T \geq T_{0}$ and for any online policy $\pi \in \Pi(B,T)$, we have that
\begin{align*}
    \sup_{B \in [T]}  \text{Regret}(B,T; \pi) \geq \textcolor{black}{(c / (1 + \beta))}  T^{\frac{1}{2} - \frac{1}{2(\beta + 1)}}  \mathbbm{1}\{\beta > 0\} + c  \log T \mathbbm{1}\{\beta = 0\}.
\end{align*}
\end{theorem}

This theorem provides an impossibility result: it says that for any fixed $\beta \in [0, \infty)$, 
there exists a distribution 
for which no online policy can achieve a better regret scaling than the one presented in Theorem~\ref{thm:CwG-regret-lower-bound}. This lower bound also highlights that the fundamental limits of the regret scaling are governed by the parameter $\beta$ which characterizes the curvature of the distribution around the gap boundaries. We observe that as $\beta \to \infty$, the scaling of regret approaches $\sqrt{T}$; i.e., no matter the online policy, it will suffer regret nearly as large as that of a simple non adaptive policy. Hence $\beta$ can be seen as characterizing the ``hardness" of an instance. The parameter $\beta$ has a physical interpretation as well. It captures how mass accumulates in the type space. For some intuition, consider the $F_{\beta}$ distribution described in \eqref{eq:example-beta-epsilon-clustered-distribution} and consider the gap boundary at $3/4$. As we move from the boundary point $3/4$ to a distance $\delta$ into the adjacent cluster, i.e., to $3/4+\delta$, the probability mass accrued grows as $C\delta^{\beta + 1}$ for some universal constant $C > 0$. Alternately, to accrue a probability mass of $\varepsilon$, we need to move a distance $C\varepsilon^{\frac{1}{\beta + 1}}$ from the boundary $3/4$ into the adjacent cluster. Therefore as $\beta$ increases, the distance one needs to travel to collect a probability mass of $\varepsilon$ also increases and this property is what makes the instances harder as $\beta$ increases. 

For $\beta = 0$, our lower bound follows from \cite{bray2022logarithmic}. To establish our bound for $\beta > 0$, we consider the distributions $F_{\beta}$ defined in \eqref{eq:example-beta-epsilon-clustered-distribution}. At a high level, we consider two events of $\Omega(1)$ probability -- one is a perturbation of the other -- under one event (denoted as $\mathcal{H}$), there are more than the expected number of arrivals with values at least $3/4$ ("high" types) and hence the hindsight threshold is (slightly) more than $3/4$, and on the other event (denoted as $\mathcal{L}$), there are fewer than expected number of arrivals with value at least $3/4$ and hence the hindsight threshold is (slightly) less than $1/4$. While the hindsight optimal policy does well on both the events, the optimal online policy can only do well on one or the other but not in both. We show that any online algorithm must make at least $\Omega(\sqrt{T})$ mistakes on at least one of the two events, and leveraging how the mass accumulates over space (characterized by Definition \ref{def:beta-epsilon-delta-clustered}), on may show that the cost of each of these mistakes $\Omega(T^{-\frac{1}{2(\beta + 1)}})$.
Combining the two gives us that the cumulative regret scales as $\Omega(T^{\frac{1}{2} - \frac{1}{2(\beta + 1)}})$.   We elaborate on this in the formal proof in Appendix \ref{onapp:universal-lower-bound}. 

\section{Algorithmic Design Principles for Near Optimal Performance}
\label{sec:cwg-principle}
Having established a spectrum of fundamental performance boundaries, it is natural to inquire if it is possible to achieve these limits, and if so, what algorithms are capable of attaining these fundamental limits. A prevalent algorithmic principle in the network revenue management literature is the Certainty Equivalent (CE) heuristic. This approach solves a deterministic approximation of a stochastic optimization problem by substituting random variables with their expected values. Given its widespread use, the CE heuristic emerges as a natural initial candidate for analysis and characterization of achievable performance. In this section, we will focus on the {\sf CE} heuristic for non-atomic distributions to avoid any tie-breaking issues which are present for atomic distributions. For the multisecretary problem, the {\sf CE} heuristic is defined as follows: at each time $t$ (before the arrival of request $\theta_t$), given a remaining budget $B_t$ and remaining number of time steps $T - t + 1$, we compute the budget ratio $B_t / (T - t + 1)$ and accept the request $\theta_t$ if and only if $r(\theta_t, \text{accept}) \geq F^{-1}(1 - B_{t} / (T - t + 1))$. Note that the {\sf CE} heuristic employs an adaptive threshold at each time $t$.


\subsection{Failure of the {\sf CE} policy under many types with gaps}
\label{sec:ce-policy-failure}
Indeed, in the case of non-atomic distributions with density uniformly bounded away from zero, \citet{lueker1998average} and \citet{bray2022logarithmic} showed that {\sf CE} achieves $\mathcal{O}(\log T)$ regret, and that this is the best scaling achievable. However, it turns out that as soon as one introduces a gap in these non-atomic distributions (as in Example \ref{ex:bimodal-distribution}), the performance of {\sf CE} degrades significantly. This phenomenon is documented in the proposition below.

\begin{proposition}[Failure of {\sf CE}]
\label{thm:ce-regret-lower-bound} 
Fix any {$\eta \in (0,1)$ and $\varepsilon \in (0,1/2]$}. Suppose the candidate-ability distribution $F$ is \emph{any} non-atomic distribution 
that has a gap of length at least $\eta$, i.e., $\exists c \in (0,1 - \eta)$ such that $F(c) = F((c + \eta)^{-} )$, and such that there is mass at least $\epsilon$ on each side of the gap, i.e., $\min\{F(c), 1 - F(c)\} \geq \varepsilon$. Then for the \textup{\sf CE} policy, there exists $T_{0} \equiv T_{0}(\varepsilon) < \infty$, a constant $c \equiv c(\eta, \varepsilon) >0$ and $B \in [T]$ such that $\text{Regret}(B, T; \textup{ \sf CE}) \geq c \sqrt{T}$ for all $T \geq T_{0}$.
\end{proposition}

The regret of the {\sf CE} policy increases dramatically if there is a gap in the types, even when one maintains the uniform distribution of types (or any other distribution) outside of the gap. As a matter of fact, the regret scaling is as large as that of a non-adaptive policy. \new{The result in Proposition \ref{thm:ce-regret-lower-bound} is analogous to the results for few types in the literature. 
 The main driver of $\Omega(\sqrt{T})$ regret scaling for both the many types with gaps and finite types settings is \emph{degeneracy}, i.e., situations where the dual variables corresponding to the initial fluid model LP are not unique. This issue is well documented in the setting with finitely many types \citep{bumpensanti2020re, vera2021bayesian}, but also manifests in the case of non-atomic distributions with gaps. As such the proof of Proposition \ref{thm:ce-regret-lower-bound} follows from the proof of the analogous result for finitely many types in  \citep[Proposition 2]{bumpensanti2020re}.}



\subsection{Conservativeness with respect to gaps}

We observed that the {\sf CE} policy breaks down for distributions with many types and ``gaps'' (intervals) of absent types; it suffers $\Omega(\sqrt{T})$ regret, as large as that of a non-adaptive algorithm. We identified that the main driver for the $\Omega(\sqrt{T})$ regret of the {\sf CE} policy is the presence of gaps. To solve this issue, we introduce a new algorithmic principle which we call ``conservativeness with respect to gaps'' ({\sf CwG}), and use it to provably achieve near optimal regret scalings for the $(\beta, \varepsilon_{0}, \delta)$-clustered distributions which allow for gaps. The idea of {\sf CwG} is that if there is a risk that the acceptance threshold based on {\sf CE} will move across a given gap in the future, then {\sf CwG} uses that gap as the acceptance threshold instead of using the {\sf CE}-based threshold. Based on the {\sf CwG} principle, we devise a new policy with the same name, which we present in Algorithm \ref{alg:CwG}.

\begin{algorithm}
	\SetAlgoNoLine
	\KwIn{Time Horizon $T$, Hiring Budget $B$, $(\beta, \varepsilon_{0}, \delta)$-clustered dist. $F$ with gaps $G_{i} = (a_{i}, b_{i})$.}
	\SetKwInOut{Initialize}{Initialize}
	\Initialize{$B_{1} = B, q_{i}^{\star} = F(a_{i}) = F(b_{i}^{-}), \forall i \in [n]$, $\Tilde{T} = \max\{0, T - \floor{ 64 \log(1 / \varepsilon_{0}) / \varepsilon_{0}^{2}}\}$}
	\SetKwFunction{SOMETHING(dfsd,fgd)}{nbb}
	    \For{$t = 1$ \KwTo $\Tilde{T}$}{
	    $p_{t}^{\textsf{CE}} = 1 - \frac{B_{t}}{T - t + 1}$\\
	    $\mathcal{S}_{t} = \Big\{i: p_{t}^{\textsf{CE}} \in \mathcal{B}\left(q_{i}^{\star}, \sqrt{\frac{2\log (T - t + 1)}{T - t + 1}} \right) \Big\}$ \\
	    \eIf{$\mathcal{S}_{t} = \emptyset$}{
	        $p_{t}^{\textsf{CwG}} = p_{t}^{\textsf{CE}}$
	    }{
	        $j_{t}^{\star} = \argmin_{i \in \mathcal{S}_{t}} |p_{t}^{\textsf{CE}} - q_{i}^{\star}|$ \\
	        $p_{t}^{\textsf{CwG}} = q_{j_{t}^{\star}}^{\star}$
	    }
	    Observe a candidate of ability $\theta_{t}$ and form the set $\mathcal{I}_{t} = \{q \in [0,1]: F^{-1}(q) = \theta_{t}\}$ \\
	    Let $X_{t}$ be a uniform sample from the set $\mathcal{I}_{t}$ \\
	    \eIf{$X_{t} \geq p_{t}^{\textup{\textsf{CwG}}}$ and $B_{t} > 0$}{
	        Hire the candidate and $B_{t + 1} \gets B_{t} - 1$}{
        	Reject the candidate and $B_{t + 1} \gets B_{t}$
	    }
	}
	Define $p_{\Tilde{T} + 1}^{\textsf{CE}} = 1 - \frac{B_{\Tilde{T} + 1}}{T - \Tilde{T}}$ \\
	\For{$t = \Tilde{T} + 1$ \KwTo $T$}{
	    Observe a candidate of ability $\theta_{t}$ and form the set $\mathcal{I}_{t} = \{q \in [0,1]: F^{-1}(q) =  \theta_{t}\}$ \\
	    Let $X_{t}$ be a uniformly random sample from the set $\mathcal{I}_{t}$ \\
	    \eIf{$X_{t} \geq p_{\Tilde{T} + 1}^{\textup{\textsf{CE}}}$ and $B_{t} > 0$}{
	        Hire the candidate and $B_{t + 1} \gets B_{t} - 1$}{
        	Reject the candidate and $B_{t + 1} \gets B_{t}$
	    }
	}
	\caption{Conservativeness with respect to Gaps (CwG)}
	\label{alg:CwG}
\end{algorithm}

The algorithm operates in two phases. For simplicity, assume that $T \geq \floor{64\log(1 / \varepsilon_{0}) / \varepsilon_{0}^{2}}$. We begin by describing the first phase. For the first $\Tilde{T} \triangleq T - \floor{64\log(1 / \varepsilon_{0}) / \varepsilon_{0}^{2}}$ steps the algorithm uses the {\sf CwG} principle, where  
if the re-solving threshold $p_{t}^{\textsf{CE}}$ is close to a gap, we modify it by instead using the quantile corresponding to the boundary of the gap as our acceptance threshold $p_{t}^{\textsf{CwG}}$. 
It remains to clarify how the quantile threshold $p_{t}^{\textsf{CwG}}$ translates to an accept/reject decision for the arrival at $t$. After observing the type $\theta_{t}$, we form the set of corresponding quantiles $\mathcal{I}_{t}$. If $\mathcal{I}_{t}$ is a singleton (this is the case if $\theta_t$ does not lie at an atom of $F$) then we have that its unique element $X_{t} = F(\theta_{t})$. 
If $\theta_t$ lies at an atom of $F$, the set $\mathcal{I}_{t}$ is a corresponding interval (recall Example~\ref{ex:discret-distributions}). If $p_{t}^{\textsf{CwG}} \notin \mathcal{I}_t$ then the hire/reject decision is unambiguous. The only case of ambiguity is $p_{t}^{\textsf{CwG}} \in \mathcal{I}_t$. To handle this case, we make use of randomization to break ties by drawing $X_{t}$ uniformly from the interval $\mathcal{I}_t$, and hiring the candidate only if the $X_{t}$ is weakly greater than $p_{t}^{\textsf{CwG}}$.

We now describe the second phase of the algorithm. In the final $\ceil{64\log(1 / \varepsilon_{0}) / \varepsilon_{0}^2 }$ time steps, the radius $\sqrt{2\log \tau / \tau}$ (where $\tau = T - t + 1$ is the number of remaining time steps) by which we measure the closeness of {\sf CE} threshold $p_{t}^{\textsf{CE}}$ and the gap quantiles $\{q_{i}^{\star}\}_{i = 1}^{n}$ becomes too large, i.e. $\sqrt{2 \log \tau / \tau} > \varepsilon_0 / 2$. This results in more than one gap quantiles being in the $\sqrt{2 \log \tau / \tau}$ neighborhood of $p_t^{\sf CE}$ which in turn makes the choice of $p_t^{\sf CwG}$ ambiguous and further complicates the regret analysis. In order to avoid this ambiguity and simplify the analysis, 
we employ a static allocation policy in the second phase: we solve for the certainty equivalent threshold $p_{\tilde{T}+1}^{\textsf{CE}}$ at time $\tilde{T}+1$, and use that threshold for the remaining $\ceil{64\log(1 / \varepsilon_{0}) / \varepsilon_{0}^{2}}$ time steps.

\subsubsection{Performance Analysis}
\begin{theorem}
\label{thm:general-cwg-upper-bound}
For any $\beta \in [0, \infty), \varepsilon_{0} \in (0,1]$ and $\delta \in (0,1]$, suppose the candidate-ability distribution $F$ with associated gaps is $(\beta, \varepsilon_{0}, \delta)$-clustered. 
Then for all $T \in \mathbb{N}$ and for all $B \in [T]$, there exists a universal constant $C< \infty$ such that the regret of the {\sf CwG} policy is upper bounded as 
    \begin{align}
    \label{eq:general-cwg-upper-bound}
    \text{Regret}(B,T;\textup{\sf CwG}) &\leq \underbrace{C(1 + 1 / \beta) (\log T)^{\frac{1}{2} + \frac{1}{2(\beta + 1)}} T^{\frac{1}{2} - \frac{1}{2(\beta + 1)}} \cdot \mathbbm{1}\{\beta > 0\} + C(\log T)^{2} \mathbbm{1}\{\beta = 0\}}_{(\spadesuit)} \nonumber\\
    &\quad  + \underbrace{C \delta \sqrt{T \log T}}_{(\vardiamond)} + \underbrace{C\sqrt{\log(1 / \varepsilon_{0})} / \varepsilon_{0}}_{(\varheart)}.
\end{align}
\end{theorem}


\new{
{\it Discussion of Theorem \ref{thm:general-cwg-upper-bound}}. The regret upper bound can be decomposed as shown in \eqref{eq:general-cwg-upper-bound}, where each of the terms has a different driver. The terms in $(\spadesuit)$ are driven by the shape of the reward distribution around gaps and is characterized by the parameter $\beta \in [0, \infty)$. Comparing the term $(\spadesuit)$ to the lower bound in Theorem \ref{thm:CwG-regret-lower-bound}, we note that the scaling of the upper bound matches the scaling of the lower bound in $T$ up to a polylogarithmic factor and hence the proposed {\sf CwG} policy is near-optimal. In the case of the {\sf CE} policy, we had identified that the main driver of its worst case regret of $\Theta(\sqrt{T})$ was the presence of gaps in the distribution of candidate abilities. Theorem~\ref{thm:general-cwg-upper-bound} tells us that one can overcome the difficulty introduced by gaps in the distribution by using the {\sf CwG} principle that we devised.
The term in $(\vardiamond)$ is driven by the parameter $\delta$ which allows us to model distributions with many {\it small} discrete types (cf. Example \ref{ex:many-small-discrete-types}). We will typically assume that $\delta$ is small and may scale as $o(1/\sqrt{T})$. Note that for the extreme cases of a few types (cf. Example \ref{ex:discret-distributions}) or continuous distributions (cf. Example \ref{ex:bray}), we have that $\delta = 0$ and hence the term in $(\vardiamond)$ disappears.
The term in $(\varheart)$ is driven by the minimum probability mass $\varepsilon_0$ and is typically assumed to be a constant in $(0,1]$. The contribution of $(\varheart)$ is attributable to the regret accrued due to the static allocation rule employed in Algorithm \ref{alg:CwG} in the last $\ceil{64 \log(1 / \varepsilon_0) / \varepsilon_0^2}$. In terms of scaling of $(\varheart)$, it matches up to polylogarithmic factors the lower bound on regret scaling of $\Omega(1/\varepsilon_0)$ presented in Lemma 1 of \cite{arlotto2019uniformly}.
}

\begin{corollary}
\label{cor:cr-upper-bound-example}
Suppose the candidate-ability distribution is $F_{0}$ where $F_{0}$ is as defined in \eqref{eq:example-beta-epsilon-clustered-distribution} with $\beta = 0$. Then we have that for all $T \in \mathbb{N}$ and for all $B\in [T]$ the regret of our {\sf CwG} policy is upper bounded as 
$\text{Regret}(B,T; \textup{\sf CwG}) \leq C(\log T)^{2}$ for the universal constant $C < \infty$ in Theorem~\ref{thm:general-cwg-upper-bound}.
\end{corollary}
\paragraph{Discussion of Corollary \ref{cor:cr-upper-bound-example}.} This corollary follows immediately from Theorem \ref{thm:general-cwg-upper-bound} by setting $\beta = 0$ and $\delta = 0$. The distribution $F_{0} = \text{Uniform}([0,1/4] \cup [3/4,1])$ is a natural variant of the uniform distribution with a gap. Corollary \ref{cor:cr-upper-bound-example} shows that regret of CwG scales as $\mathcal{O}((\log T)^{2})$ for the distribution $F_{0}$. 
This is a significant improvement on the $\Omega(\sqrt{T})$ regret scaling of the {\sf CE} policy for the same distribution $F_{0}$, and the regret of the {\sf CwG} policy for the $F_{0}$ distribution is only a $\log T$ factor larger than the regret for the uniform distribution. The key takeaway from Corollary \ref{cor:cr-upper-bound-example} in conjunction with Proposition \ref{thm:ce-regret-lower-bound} is that the presence of gaps is not a fundamental driver of the achievable regret performance, and one can overcome the difficulty posed by gaps by using the {\sf CwG} principle.

\begin{corollary}[Constant Regret for discrete distributions]
\label{cor:discrete-distribution-bound}
Suppose the candidate-ability distribution is $F$ where $F$ is a discrete distribution as desribed in Example \ref{ex:discret-distributions}. 
Then, for all $T \in \mathbb{N}$ and for all $B \in [T]$, we have $\text{Regret}(B,T; \textup{\sf CwG}) \leq C \sqrt{\log (1 / \varepsilon_{0})} / \varepsilon_{0}$ for a universal constant $C < \infty$.
\end{corollary}

\begin{remark}
\label{rmk:discrete-distribution-constant-regret}
The discrete distribution considered in Example \ref{ex:discret-distributions} belongs the class of $(0, \varepsilon_{0})$-clustered distributions and hence from Corollary \ref{cor:cr-upper-bound-example}, it follows that $\text{Regret}(B,T) = \mathcal{O}\left((\log T)^{2} \right)$. However, recall from Example \ref{ex:discret-distributions} that for discrete distributions we have that $\forall i \in [n + 1], \forall q, \tilde{q} \in Q_{i}, |F^{-1}(q) - F^{-1}(\tilde{q}) | = 0$, and this distinguishes discrete distributions from general $(0, \varepsilon_{0})$-clustered distributions. This distinction allows us to obtain stronger regret guarantees than the one implied by Corollary \ref{cor:cr-upper-bound-example} and recover the result of \citet{arlotto2019uniformly}. The proof of Corollary \ref{cor:discrete-distribution-bound} follows by modifying the analysis leading to Theorem~\ref{thm:general-cwg-upper-bound}. The modifications enable us to eliminate the $C(\log T)^2$ term in the regret bound in Theorem 2. We defer the details to Appendix \ref{sec:proof-cwg-discrete-type}. 
\end{remark}

\begin{corollary}[Regret for non-atomic distribution with contiguous support]
\label{cor:beta-one-clustered-regret}
For any $\beta \in [0, \infty)$, $\varepsilon_{0} = 1$, and $\delta = 0$, suppose the candidate-ability distribution $F$  is $(\beta, \varepsilon_0 = 1, \delta = 0)$-clustered ($F$ has no non-trivial gaps). Then for all $T \in \mathbb{N}$ and for all $B \in [T]$, there exists a universal constant $C  < \infty$ such that the regret of our CwG policy is 
\begin{align*}
\text{Regret}(B,T; \textup{CwG}) \leq C\left(1 + \frac{1}{\beta} \right) T^{\frac{1}{2} - \frac{1}{2(1 + \beta)}} \mathbbm{1}\{\beta > 0\} + C \log T \mathbbm{1}\{\beta = 0\}
\end{align*}
\end{corollary}

\paragraph{Discussion of Corollary \ref{cor:beta-one-clustered-regret}:} 
This corollary follows immediately from Theorem \ref{thm:general-cwg-upper-bound} by setting $\varepsilon_0 = 1$ and $\delta = 0$, except for some polylogarithmic factors.
The class of $(\beta,\varepsilon_0 = 1, \delta = 0)$-clustered distributions allows for the pdf $f$ to be zero at some points. 
An example of such a distribution is given by $\tilde{F}_{\beta}(x) = \left(0.5 - 2^{\beta}\left(0.5 - x \right)^{\beta + 1} \right) \mathbbm{1}\{x \leq 0.5\} + \left(0.5 + 2^{\beta} \left(x - 0.5 \right)^{\beta + 1}\right) \mathbbm{1}\{x > 0.5\}$ where the pdf $f$ is zero at $x = 0.5$. Since there are no non-trivial gaps for the distribution $\tilde{F}_{\beta}$, we choose to treat the whole interval $[0,1]$ as a single cluster and hence have $\varepsilon_0 = 1$. It can be easily verified that $\tilde{F}_{\beta}$ satisfies the ``cluster density requirement'' in Definition \ref{def:beta-epsilon-delta-clustered} with $\delta = 0$. Note that the distribution $\tilde{F}_{\beta}$ is not admissible under the assumptions of \cite{bray2022logarithmic} for $B = T/2$ and $\beta > 0$. Since there are no gaps of positive length in $(\beta,1)$-clustered distributions, the {\sf CwG} policy boils down to the {\sf CE} policy. If the probability density function $f$ is bounded below by a constant, we have $\beta = 0$ and we recover the $\mathcal{O}(\log T)$ scaling in \citet{bray2022logarithmic}. If $f$ is zero at some points, then the regret scaling is determined by $\beta$ which quantifies how the mass accumulates around types where $f$ is zero. This result, in conjunction with Theorem \ref{thm:CwG-regret-lower-bound}, proves that the {\sf CE} policy is near-optimal in the absence of non-trivial gaps.

\new{
\subsection{Achieving Conservativeness with respect to Gaps via a Simulation-based Policy}
\label{sec:cwg-via-simulation}
In Algorithm~\ref{alg:CwG}, if the re-solving threshold $p_{t}^{\textsf{CE}}$ at time $t= T-\tau+1$ was within $\sqrt{2\log \tau/ \tau}$ of a gap, we modified it as by instead using the quantile corresponding to the boundary of the gap  as our acceptance threshold.  An alternative to this method is a simulation-based approach, which we'll outline next, followed by a full treatment in the next section.}

\begin{figure}[h!]
    \centering
    \begin{subfigure}{0.40\linewidth}
        \centering
        \begin{tikzpicture}
            \draw[->, thick] (0,0) -- (6,0);
            \draw[->, thick] (0,0) -- (0,3);
            \draw[-, thick] (0,1.5) -- (2,1.5);
            \draw[-,thick] (3,1.5) -- (5,1.5);
            \draw[-,dotted,thick] (5,1.5) -- (5,0);
            \draw[-,dotted, thick] (2,1.5) -- (2,0);
            \draw[-,dotted, thick] (3,1.5) -- (3,0);
            \node at (0,-0.3) {$0$};
            \node at (5,-0.3) {$1$};
            \node at (5.8,-0.3) {\scriptsize $\theta$};
            \node at (0.4, 2.7) {\scriptsize $f(\theta)$};
            \draw[-, ultra thick, red] (1.8,2) -- (1.8,-0.3);
            \node at (1.8,2.4) { \textcolor{red}{$F^{-1}(p_t^{\sf CE})$}};
            \draw[->, ultra thick] (1.8, 2.1) -- (2,2.1);
            \draw[-, ultra thick, blue] (2,2) -- (2,-0.3);
            \node at (2.1,-0.6) { \textcolor{blue}{$F^{-1}(p_t^{\sf CwG})$}};
        \end{tikzpicture}
        \caption{{\footnotesize Implementation of CwG via Algorithm~\ref{alg:CwG}}}
        \label{fig:explicit-cwg}
    \end{subfigure}
    \begin{subfigure}{0.58\textwidth}
    \centering
        \begin{tikzpicture}
            \draw[->, thick] (0,0) -- (6,0);
            \draw[->, thick] (0,0) -- (0,3);
            \draw[-, thick] (0,1.5) -- (2,1.5);
            \draw[-,thick] (3,1.5) -- (5,1.5);
            \draw[-,dotted,thick] (5,1.5) -- (5,0);
            \draw[-,dotted, thick] (2,1.5) -- (2,0);
            \draw[-,dotted, thick] (3,1.5) -- (3,0);
            \node at (0,-0.3) {$0$};
            \node at (5,-0.3) {$1$};
            \node at (5.8,-0.3) {\scriptsize $\theta$};
            \node at (0.4, 2.7) {\scriptsize $f(\theta)$};
            \draw[-, ultra thick, red] (1.8,2) -- (1.8,-0.3);
            \node at (1.8,2.4) {\textcolor{red}{$F^{-1}(p_t^{\sf CE})$}};
            \draw[->, ultra thick] (1.8, 2.1) -- (2.25,2.1);
            \draw[-, ultra thick, dashed, teal] (1.9, 2) -- (1.9, -0.3);
            \node at (2.2,-0.5) {\textcolor{teal}{$\bm{\theta}_{[B_t]}^{(1)}$}};
            \draw[-, ultra thick, dashed, teal] (1.65, 2) -- (1.65, -0.3);
            \node at (1.6,-0.5) {\textcolor{teal}{$\bm{\theta}_{[B_t]}^{(2)}$}};
            \draw[-, ultra thick, dashed, teal] (3.2,2) -- (3.2,-0.3);
            \node at (3.2,-0.5) { \textcolor{teal}{$\bm{\theta}_{[B_t]}^{(3)}$}};
            \draw[-, ultra thick, teal] (2.25,2) -- (2.25,-0.3);
        \end{tikzpicture}
        \caption{{\footnotesize CwG arising organically from a simulation-based approach}}
        \label{fig:implicit-cwg}
    \end{subfigure}
    \caption{Implemention of the {\sf CwG} principle using two algorithmic approaches.}
    \label{fig:ms-sim-example}
\end{figure}
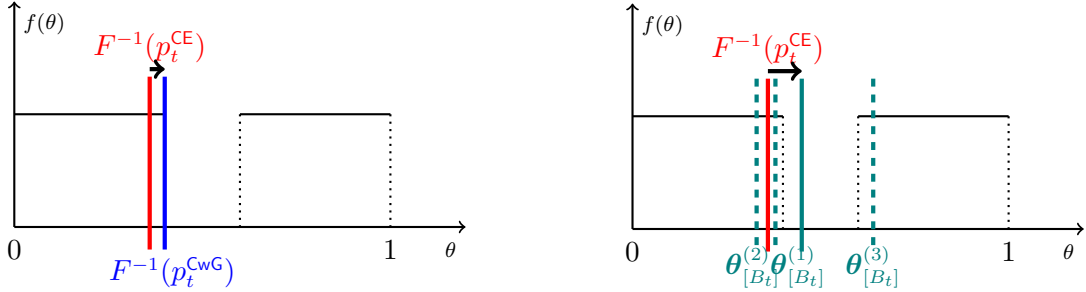

\new{Consider the bimodal uniform distribution, described by \eqref{eq:example-beta-epsilon-clustered-distribution} with $\beta = 0$. Assume the {\sf CE} threshold at time $t$, denoted as $p_{t}^{\textsf{CE}}$, is $1/2 - \epsilon$, where $\epsilon$ is sufficiently small ($\epsilon < \sqrt{2 \log \tau / \tau}$, where $\tau = T - t + 1$). 
Under Algorithm \ref{alg:CwG}, the {\sf CwG} quantile threshold is set to $p_t^{\sf CwG} = 1/2$. Consequently, only abilities with values of at least $3/4$ will be accepted at time $t$. This is illustrated in Figure \ref{fig:explicit-cwg}, where the threshold 
shifts from $F^{-1}(p_{t}^{\sf CE}) = 1/4 - 2\epsilon$ (in red) to $F^{-1}(p_t^{\sf CwG}) = 1/4$ (in blue). 
}

\new{On the other hand, consider the following simulation-based approach: simulate multiple future demand scenarios. For the $i$-th simulated scenario, let $\bm{\theta}^{(i)}_{[B_t]}$ denote the value of the $B_t$-th largest candidate ability on the simulated sample path $\tilde{\bm{\theta}}_{\geq t + 1}^{(i)}$, where $B_t$ is the remaining budget. The candidate with ability $\theta_t$ is accepted if $\theta_t \geq K_t^{-1} \sum_{k = 1}^{K_t} \bm{\theta}_{[B_t]}^{(i)}$ where $K_t$ is the number of scenarios. Figure \ref{fig:implicit-cwg} illustrates this simulation-based approach using three simulated demand scenarios, with the $B_t$-th largest value in each of the demand scenarios (denoted as $\{\bm{\theta}_{[B_t]}^{(i)}\}_{i = 1}^3$) being depicted as the dashed green lines. 
The average of these values (depicted as a solid green line in Figure \ref{fig:implicit-cwg}) falls within the gap interval $(1/4,3/4)$, resulting in only abilities of at least $3/4$ being accepted at time $t$. The simulation-based approach yields the same action as the carefully crafted {\sf CwG} policy (Algorithm \ref{alg:CwG}). Interestingly, as we will later explore in Section \ref{sec:rams}, this simulation-based approach inherits the regret guarantee of the {\sf CwG} policy (cf. Corollary \ref{cor:cwg-rams-guarantee}), and outperforms the {\sf CwG} policy in numerical experiments (cf. Figure \ref{fig:cwg_ce_comparison}). It is worth noting that the $\{\bm{\theta}_{[B_t]}^{(i)}\}_{i = 1}^3$ values represent shadow prices for the single resource (the hiring budget) under the three different demand scenarios. In the simulation-based approach, the candidate ability $\theta_t$ is accepted if its reward $\theta_t$ exceeds the approximated average shadow price, obtained by averaging the shadow prices over multiple demand scenarios, i.e., $\sum_{i = 1}^3 \bm{\theta}_{[B_t]}^{(i)} / 3$. 
Importantly, as we present next, this simulation-based idea is not limited to the multisecretary problem but can be applied more broadly to dynamic resource allocations, such as network revenue management and online matching, and notably inherits performance guarantees which hold for \emph{any} algorithm satisfying certain conditions in these settings.
}

\section{Unifying Algorithm:~Repeatedly Act using Multiple Simulations}
\label{sec:rams}

\simulation{
In this section, we will lift the idea of using simulations to drive decisions from the multisecretary setting to the broader class of NRM and online matching problems. We dub the resulting natural and versatile simulation-based algorithm {\sf Repeatedly Act using Multiple Simulations (RAMS)}. Prior to formally presenting {\sf RAMS}, we will establish some notations.
Let $V_{t}^{\textsf{hs}}(B_t; \vect{\theta}_{\geq t})$ denote the hindsight optimal value for a given tail sequence of requests $\bm{\theta}_{\geq t}\triangleq \{\theta_t, \dots, \theta_{T}\}$ and remaining budget $B_t$, 
\begin{align}
    \label{eq:hindsight-value-function}
    V_{t}^{\textsf{hs}}(B_t; \vect{\theta}_{\geq t}) \equiv \max_{\mathbf{a} \in |\mathcal{A}|^{T - t + 1}} \sum_{k= t}^T r(\theta_k, a_k) \ \ \text{s.t.} \sum_{k= t}^T c(\theta_k, a_k) \leq B_t.
\end{align}

Furthermore, it is natural to define $V_{T + 1}^{\textsf{hs}}(B_{T + 1}, \emptyset) \equiv 0, \forall B_{T+1}$. 
We will assume access to a simulator $\mathcal{S}$ which takes as input
a history $\mathcal{H}$ 
of request arrivals and random seed $U$ and produces a simulated demand scenario. 
Here a demand scenario is a tail sequence of requests $\vect{\theta}_{\geq t+1}$; we remark that the order of requests in a tail sequence will not matter to {\sf RAMS}, since it will perform a hindsight-based calculation. 
Note that the assumption of access to a simulator is a weaker and more practical assumption than knowledge of the distribution $F$. 
This permits {\sf RAMS} to be a data-driven algorithm where distributional knowledge $F$ is replaced by a high fidelity simulator based on historical data. 
Additionally, while most of our previous discussion was focused on a stationary setting with i.i.d requests, {\sf RAMS} could be 
applicable in non-stationary settings where the request types may have some form of temporal correlations, corresponding to the reality of many applications. This is due to the fact that {\sf RAMS} is completely agnostic to the underlying type distribution.
}

\subsection{Algorithmic Description}
\simulation{
The basic idea behind \textsf{RAMS} is as follows: given the remaining budget $B_t$ at time $t$, upon observing a request $\theta_{t}$, we \emph{simulate} $K_{t}$ sample paths of the future 
denoted as $\{\tilde{\vect{\theta}}_{\geq t + 1}^{(i)}\}_{i = 1}^{K_t}$. On each of the these {simulated} sample paths $\vect{\tilde{\theta}}^{(i)}_{\geq t + 1}$, we compute the maximum achievable cumulative reward in hindsight under each possible action $a \in \mathcal{A}(B_t, \theta_t)$ at time $t$, denoted by $Q_t^{\sf hs}(B_t, a; \tilde{\vect{\theta}}_{\geq t}^{(i)})$ where $\tilde{\vect{\theta}}_{\geq t}^{(i)} \triangleq \{\theta_{t}\} \cup \tilde{\vect{\theta}}_{\geq t + 1}^{(i)}$. For each action $a \in \mathcal{A}(B_t, \theta_t)$ we average over the $K_t$ simulated sample paths, and choose the action which maximizes the average cumulative reward, i.e., $\argmax_{a \in \mathcal{A}(B_t, \theta_t)} K_t^{-1} \sum_{i = 1}^{K_t} Q_t^{\sf hs}(B_t, a;  \tilde{\vect{\theta}}_{\geq t + 1}^{(i)})$. We formally describe \textsf{RAMS} in Algorithm \ref{alg:rams}.
}

\begin{algorithm}[h]
   	\SetAlgoNoLine
	\KwIn{Time Horizon $T$, Budget $B \in \mathbb{R}_+^d$, simulator $\mathcal{S}$, Sequence of number of simulated sample paths $\{K_{t}\}_{t = 1}^{T}$}
	\SetKwInOut{Initialize}{Initialize}
	\Initialize{$B_{1} = B, \mathcal{H} = \emptyset$}
	\SetKwFunction{SOMETHING(dfsd,fgd)}{nbb}
	    \For{$t = 1$ \KwTo $T$}{
	    Observe the request $\theta_{t}$  \\
        $\mathcal{H} \gets \mathcal{H} \cup \{\theta_t\}$ \\
        Make $K_t$ conditionally independent calls to the simulator $\mathcal{S}$ with 
        history $\mathcal{H}$ and random seed $U \sim \text{Unif}([0,1])$ (denote the $K_t$ simulated sample paths of requests as $\{\tilde{\vect{\theta}}_{\geq t + 1}^{(i)}\}_{i = 1}^{K_t}$.)
         \\
        \For{$i = 1$ \KwTo $K_{t}$}{
            \For{$a \in \mathcal{A}(B_t, \theta_t)$}{
            \vspace{-1em}
            \begin{align}   
                \label{eq:rams-opt-Q}
                \hspace{-3em} Q_{t}^{\sf hs}\left(B_{t}, a;  \tilde{\vect{\theta}}_{\geq t }^{(i)}\right) &= r(\theta_t, a)  + \bigg\{\max_{(a_k)_{k > t}} \sum_{k > t}^{} r(\tilde{\theta}_k^{(i)}, a_k) 
                \text{ s.t.} \sum_{k > t}^{} c(\Tilde{\theta}_k^{(i)}, a_k) \leq B_{t} - c(\theta_t, a) \bigg\} 
            \end{align}
            }
        }
        Take the action $a_{t} = \argmax_{a \in \mathcal{A}{(B_t, \theta_t)}} K_t^{-1} \sum_{i = 1}^{K_t} {Q}_{t}^{\sf hs}\left(B_{t}, a; \tilde{\vect{\theta}}_{\geq t}^{(i)}\right)$ \\
        $B_{t + 1} \gets B_{t} - c(\theta_t, a_t)$
        }
	\caption{\textsf{Repeatedly Act using Multiple Simulations (RAMS)}}
	\label{alg:rams}
\end{algorithm}

\simulation{
For a feasible online policy $\pi$, given a state $B_t$ {and an action $a$ which is feasible in that state $a \in \mathcal{A}(B_t, \theta_t) \subseteq \mathcal{A}$}, define the following $Q$-function
\begin{align*}
    Q^{\pi}_t(B_t, a; \theta_t) &= r(\theta_t, a) + \mathbb{E}\left[ \sum_{k = t + 1}^T r(\theta_k, a_k^{\pi}) \right] , \ \
    Q^{\star}_t(B_t, a; \theta_t) = \max_{{\pi \in \Pi(B_t-c(\theta_t,a), T-t+1)}} \ \  Q^{\pi}_t(B_t, a; \theta_t)\, .
\end{align*}
The action under the optimal online policy is $\argmax_{a \in \mathcal{A}} Q_{t}^{\star}(B_t, a; \theta_t)$, 
however computing this dynamic programming solution may be infeasible in general. Instead {\sf RAMS} utilizes the ``hindsight-based'' approximation to the $Q$-function, estimated from simulated futures, $K_t^{-1}\sum_{i = 1}^{K_{t}} Q_{t}^{\sf hs}(B_t, a; \tilde{\vect{\theta}}_{\geq t}^{(i)}) \approx \mathbb{E}_{\vect{\theta}_{\geq t + 1}}\left[Q_t^{\sf hs}(B_t, a; \vect{\theta}_{\geq t}) \right]$ as a proxy to make allocation decisions. Note that $\mathbb{E}_{\vect{\theta}_{\geq t + 1}}\left[Q_t^{\sf hs}(B_t, a; \vect{\theta}_{\geq t}) \right] \geq Q_t^{\star}(B_t, a; \theta_t)$ and from \eqref{eq:rams-opt-Q}, we have that $Q_t^{\sf hs}(B_t, a; \tilde{\vect{\theta}}_{\geq t}^{(i)}) = V_{t + 1}^{\sf hs}(B_t - c(\theta_t, a); \tilde{\vect{\theta}}_{\geq t + 1}^{(i)}) + r(\theta_t, a)$.
Next we define {\it marginal compensation}  for a given action $a$ at time $t$ \citep{vera2021bayesian}.  Intuitively speaking, {\it marginal compensation} is the minimum payment one must make to an agent who knows the future to persuade that agent to take action $a$ at time $t$ on a realized sample path. 
\begin{definition}[Marginal Compensation]
    \label{def:marginal-compensation}
    Given budget $B_t \geq \bm{0}$ and tail sequence of requests $\vect{\theta}_{\geq t}$ for some $t \in [T]$, for any action $a \in \mathcal{A}(B_t, \theta_t)$, we define
    \begin{align}
        \label{eq:marginal-compensation-sample-path}
        \partial \mathcal{R}_{t}( B_{t}, a;  \vect{\theta}_{\geq t}) &\triangleq V_{t}^{\textsf{hs}}(B_{t}; \vect{\theta}_{\geq t}) - \left[V_{t + 1}^{\textsf{hs}}\left( B_{t} - c(\theta_t, a); \vect{\theta}_{\geq t + 1}\right) + r\left(\theta_{t}, a \right) \right] \\
        \label{eq:marginal-compensation-sample-path-averaged}
        \partial \mathcal{R}_{t}(B_{t}, a) &\triangleq \mathbb{E}_{\vect{\theta}_{\geq t}}\left[\partial \mathcal{R}_{t}(B_{t}, a; \vect{\theta}_{\geq t}) | B_{t} \right].
    \end{align}
\end{definition}
We refer to $\partial \mathcal{R}_{t}( B_{t}, a;  \vect{\theta}_{\geq t})$ as marginal compensation and $\partial \mathcal{R}_{t}(B_{t}, a)$ as the expected marginal compensation.
A key fact from \citep[][Lemma 1]{vera2021bayesian} is that the expected regret of a policy can be decomposed as the sum of the expected marginal compensations for the actions taken by the policy, as formalized below 
\begin{lemma}
    \label{lem:cwg-regret-decomposition}
    For all $T \in [N]$ and budget $B \in [T]$, consider any online policy $\pi \in \Pi(B,T)$ and let $B_{t}^{\pi}$ denote the remaining budget at time $t$ under policy $\pi$. Then we have that 
    \begin{align}
        \text{Regret}(B,T; \pi) &= \sum_{t = 1}^{T} \mathbb{E}_{B_{t}^{\pi}}\left[\partial \mathcal{R}_{t}(B_{t}^{\pi}, a_{t}^{\pi}) \right]. 
    \end{align}
\end{lemma}

\begin{lemma}[{\sf RAMS} is equivalent to minimizing expected marginal compensation]
    \label{lem:rams-equivalent-marginal-compensation-minimization}
    Given a budget $B_t$, request $\theta_t$ and a collection of simulated sample paths $\{\tilde{\bm{\theta}}_{\geq t + 1}^{(i)}\}_{i = 1}^{K_t}$, {\sf RAMS} takes an action $a_t \in \mathcal{A}(B_t, \theta_t)$ at time $t$ which minimizes the simulation-based estimate of expected marginal compensation, i.e. $a_t = \argmin_{a \in \mathcal{A}(B_t, \theta_t)} K_t^{-1} \sum_{i = 1}^{K_t} \partial \mathcal{R}_{t}(B_t, a; \tilde{\bm{\theta}}^{(i)}_{\geq t} )$ where $\tilde{\bm{\theta}}_{\geq t}^{(i)} = \{\theta_t\} \cup \tilde{\bm{\theta}}_{\geq t+ 1}^{(i)}$.
\end{lemma}
Lemma \ref{lem:rams-equivalent-marginal-compensation-minimization} follows immediately from \eqref{eq:rams-opt-Q} and \eqref{eq:marginal-compensation-sample-path} and provides an alternate description of {\sf RAMS}.
}

\subsection{Performance Analysis: Meta Theorem for \textsf{RAMS}}
\simulation{
Since the expected regret of the policy is the sum of the expected marginal compensations (Lemma \ref{lem:cwg-regret-decomposition}), and RAMS performs a simulation-based minimization of the expected marginal compensation (Lemma \ref{lem:rams-equivalent-marginal-compensation-minimization}), it follows that RAMS provides the ``best achievable'' regret performance (in a certain sense). This reasoning is formalized in the following meta theorem. 
\begin{theorem}[Meta Performance of \textsf{RAMS}]
    \label{thm:meta-performance-rams}
    Consider an online resource allocation problem with horizon $T$, number of resources $d$, initial budget $B \in \mathbb{R}^d$, a finite action set $\mathcal{A}$ and request distribution $F$ as defined in Section \ref{sec:model}.
    Assume the following
    \begin{itemize}
        \item[(i)] There exists an algorithm {\sf ALG} for the online resource allocation problem such that the expected marginal compensation is uniformly bounded at each $1 \leq t \leq T $ as per $\sup_{B_t \geq \bm{0}} \partial \mathcal{R}_{t}(B_t, a_t^{\sf ALG}) \leq \Delta_t({\sf ALG})$ 
        where $B_t$ is the remaining budget at time $t$ and $a_t^{\sf ALG}$ is the action under {\sf ALG}.
        \item[(ii)] There exists a constant $\mathcal{C} \equiv \mathcal{C}(F) < \infty$ such that {the marginal compensation in a time step is uniformly bounded by $\mathcal{C}$, i.e.,} $\sup_{B_t, a, \bm{\theta}_{\geq t}}\partial \mathcal{R}_{t}\left(B_t, a; \bm{\theta}_{\geq t} \right) \leq \mathcal{C}$ for all $t \geq 1$.
    \end{itemize}
    Let $K_t$ denote the number of simulated sample paths drawn at time $t$. Then for any $\eta > 2$, there exists a constant $C \equiv C(\eta, |\mathcal{A}|, \mathcal{C}(F)) < \infty$, such that 
    \begin{align*}
        \text{Regret}(B,T; \textup{\sf RAMS}) \leq  \sum_{t = 1}^T \Delta_{t}(\textup{\sf ALG}) + C \sum_{t = 1}^{T} K_{t}^{-\frac{1}{\eta}} 
    \end{align*}
\end{theorem}
}

\paragraph{Discussion of Theorem \ref{thm:meta-performance-rams}.}
\simulation{
Note that while the theorem has been stated for the {\it i.i.d} setting, Theorem \ref{thm:meta-performance-rams} can also apply to non-stationary settings with some form of temporal correlations. Theorem \ref{thm:meta-performance-rams} states that the regret of \textsf{RAMS} can be broken down into two components: $\Delta_{t}(\textsf{ALG})$ and $K_{t}^{-\frac{1}{\eta}}$. The former term $\Delta_{t}(\textsf{ALG})$ follows from the assumed uniform (over the states) upper bound on the expected compensation $\partial \mathcal{R}_t(B_t, a_t^{\sf ALG})$ under algorithm \textsf{ALG}, while the latter term $K_{t}^{-\frac{1}{\eta}}$ is due to the finite number of simulated sample paths. Theorem~\ref{thm:meta-performance-rams} states that \textsf{RAMS} inherits -- up to sampling error -- the best (uniform) regret guarantee which holds for any algorithm. 
Our numerical observations show that \textsf{RAMS} outperforms regret-optimal algorithms tailored for specific distributions or problem contexts, without the need for tuning (see Section \ref{sec:numerical}). Notably, neither \textsf{RAMS} nor the meta theorem (Theorem \ref{thm:meta-performance-rams}) require prior knowledge of these optimized algorithms. As long as there exist algorithms that satisfy assumption $(i)$ and that $(ii)$ holds, \textsf{RAMS} achieves the same regret scaling. 

 We highlight that there exist algorithms developed in this  and prior work for different problem settings which satisfy assumption \emph{(i)} (cf. Corollaries \ref{cor:cwg-rams-guarantee}-\ref{cor:online-matching-constant}). Coming to assumption \emph{(ii)}, in the context of network revenue management problem, this assumption holds under mild conditions, as captured in the following claim.} 
\simulation{
\begin{claim}
    \label{claim:sufficient-condition-compensation-uniform-bound}
    In the context of the NRM problem, for any request type $\theta = (r_{\theta}, \bm{c}_\theta) \in \Theta$, assume that the consumption vector $\bm{c}_{\theta}$ is bounded i.e., $\underline{\nu} \leq \|\bm{c}_{\theta}\|_{\infty} \leq \bar{\nu}$ for $0 < \underline{\nu} \leq \bar{\nu} < \infty$. Then we have that $\sup_{B_t, a, \bm{\theta}_{\geq t}} \partial \mathcal{R}_t\left(B_t, a; \bm{\theta}_{\geq t} \right) \leq d r_{\max} \bar{\nu}  / \underline{\nu} \triangleq {\cal C}(F)$ where $d$ is the number of resources and $r_{\max} \equiv \max_{\theta \in \Theta} r_{\theta} \leq 1$ (by assumption).
\end{claim}
Note that the sufficient condition in Claim \ref{claim:sufficient-condition-compensation-uniform-bound} permits  many (or infinitely many) consumption types, in contrast to the typical assumption in the prior literature of a small number of consumption types  \citep[with some notable exceptions][]{lueker1998average, arlotto2020logarithmic, li2022online, bray2022logarithmic}.
}

\simulation{Combining Theorem \ref{thm:meta-performance-rams} with analyses of specific algorithms, we can show that {\sf RAMS} achieves the same regret scaling as that of the {\sf CwG} algorithm (Algorithm \ref{alg:CwG}) for the class of $(\beta, \varepsilon_0, \delta)-$clustered distributions (Corollary \ref{cor:cwg-rams-guarantee}). Zooming out from the multisecretary problem, we consider the more general network revenue management and online matching problems. We show that under the assumption of a small number of discrete types, {\sf RAMS} achieves bounded regret scaling for both the network revenue management (Corollary \ref{cor:nrm-regret-profiles}(a)) and online matching (Corollary \ref{cor:online-matching-constant}). Under infinitely many types and some structural assumptions, {\sf RAMS} achieves logarithmic (Corollary \ref{cor:nrm-regret-profiles}(b)) and log-squared regret (Corollary \ref{cor:nrm-regret-profiles}(c)) scaling for the general NRM problem in line with state of the art algorithms presented in \cite{bray2022logarithmic} and \cite{jiang2022degeneracy} respectively. 
Detailed assumptions and corollaries are presented in Appendix \ref{onapp:recovering-existing-guarantees-rams} due to space constraints.
}


\subsection{Connection of {\sf RAMS} to prior work}
\new{ Due to the equivalence of {\sf RAMS} to minimizing the expected compensation at each time period (cf. Lemma \ref{lem:rams-equivalent-marginal-compensation-minimization}), {\sf RAMS} follows the ``Bayes Selector" principle developed in \cite{vera2021bayesian}. However, the focus of \cite{vera2021bayesian} is on settings with a few types and hence their algorithm has been tailored for such settings, whereas {\sf RAMS} is a very general algorithm which does not require any knowledge of the underlying assumptions on the type space.

In the context of network revenue management, {\sf RAMS} is a refined version of the dual averaging policy proposed in \cite{talluri1999randomized}, where dual prices are computed for multiple demand scenarios and the allocation decisions are made by averaging these dual prices over the different scenarios. Under {\sf RAMS}, given a remaining budget $B_t$, a request $\theta_t$ is accepted if $K_t^{-1} \sum_{i = 1}^{K_t} Q_t^{\sf hs}(B_t, \text{accept}; \tilde{\bm{\theta}}_{\geq t}^{(i)}) \geq K_t^{-1} \sum_{i = 1}^{K_t} Q_t^{\sf hs}(B_t, \text{reject}; \tilde{\bm{\theta}}_{\geq t}^{(i)})$. 
Assume that there exists a dual vector $\mu(B_t; \tilde{\bm{\theta}}_{\geq t + 1}^{(i)})$ for \eqref{eq:hindsight-value-function} with tail sequence $\tilde{\bm{\theta}}_{\geq t + 1}^{(i)}$ such that first order approximation of $V_{t + 1}^{\sf hs}(B_t; \tilde{\bm{\theta}}_{\geq t + 1}^{(i)})$ is good, i.e., $V_{t + 1}^{\sf hs}(B_t; \tilde{\bm{\theta}}^{(i)}_{\geq t + 1}) - V_{t + 1}^{\sf hs}(B_t - c(\theta_t, \text{accept}); \tilde{\bm{\theta}}^{(i)}_{\geq t + 1}) \approx \mu(B_t; \tilde{\bm{\theta}}_{\geq t + 1}^{(i)})^{\top} c(\theta_t, \text{accept})$. Then, using \eqref{eq:hindsight-value-function}, \eqref{eq:rams-opt-Q} and the fact that $r(\theta_t, \text{reject}) = 0$ and $c(\theta_t, \text{reject}) = \bm{0}$,  under {\sf RAMS}, the request $\theta_t$ is accepted if 
\begin{align*}
    r(\theta_t, \text{accept}) &\geq
    \frac{1}{K_t} \sum_{i = 1}^{K_t} \left( V_{t + 1}^{\sf hs}(B_t; \tilde{\bm{\theta}}^{(i)}_{\geq t + 1}) - V_{t + 1}^{\sf hs}(B_t - c(\theta_t, \text{accept}); \tilde{\bm{\theta}}^{(i)}_{\geq t + 1})    \right) \\
    &\approx \frac{1}{K_t} \sum_{i = 1}^{K_t} \mu(B_t; \tilde{\bm{\theta}}_{\geq t + 1}^{(i)})^{\top} c(\theta_t, \text{accept}) = \underbrace{\left(\frac{1}{K_t} \sum_{i = 1}^{K_t} \mu(B_t; \tilde{\bm{\theta}}_{\geq t + 1}^{(i)}) \right)^{\top}}_{\text{average dual price for bid price control}}  c(\theta_t, \text{accept}).
\end{align*}
Therefore, assuming that the first order approximation is good, {\sf RAMS} will accept the request $\theta_t$ if the reward exceeds the sum of the average dual prices for the resources it consumes, and this resembles the bid price control policy \citep{talluri1998analysis}.
Thus we see that dual averaging is, in fact, an 
approximate version of {\sf RAMS} for settings in which individual actions have a ``small'' impact, and our theoretical backing for {\sf RAMS} (Theorem~\ref{thm:meta-performance-rams}) provides new justification for why dual averaging should work well in such settings.
Dual averaging is very practical and requires only a small adaptation of  dual-based dynamic resource allocation systems based on model predictive control, which are typical in the industry, e.g., in supply chain optimization. Specifically, 
it only requires the construction of multiple demand scenarios. The hindsight problem for each scenario can be solved in parallel (using the existing MPC solver as is) and then a simple dual averaging layer can be inserted before the decision making layer. 
}

\simulation{{\sf RAMS} can be viewed as the manifestation in our setting of the so-called Multi Forecast--Model Predictive Control ({\sf MF-MPC}) policy which appears in the control literature, e.g., see \cite{shen2021incremental} and citations therein. In {\sf MF-MPC}, one constructs multiple plausible forecasts of the future, termed \emph{scenarios}, and constructs a different plan for each of the possible scenarios, while imposing the constraint that the plans must agree on the present action to be chosen. This process is repeated each time an action is to be chosen. The connection with {\sf MF-MPC} further reveals an illuminating interpretation for {\sf RAMS}: 
Suppose all uncertainty about the future will be resolved right after the current action is chosen. What current action is optimal in this proxy problem? This is the action chosen by {\sf RAMS} at each time; after all, by definition, {\sf RAMS} solves the Bellman equation for this proxy problem. This interpretation throws light on the approximation underlying {\sf RAMS}, and may help us --in future work-- to understand how well {\sf RAMS} (or, more generally, any compensation-based approach) can approximate the optimal MDP solution in a given setting.}

\subsection{Numerical Simulations}
\label{sec:numerical}
\simulation{
We perform numerical experiments under different assumptions and for different problem classes. For the multisecretary problem, we study the performance of the {\sf CwG} algorithm for different distributions (Figure \ref{fig:cwg_diff_dists}),  compare the performance of {\sf CE}, {\sf CwG} and {\sf RAMS} for the bimodal uniform distribution $F_0$ (Figure \ref{fig:cwg_ce_comparison}),  and  study the impact of $\beta > 0$ (Figure \ref{fig:beta-comparison}). In addition,  we consider the general network revenue management problem with a few types and two resources and compare the performance of previous algorithms with that of {\sf RAMS} (Figure \ref{fig:nrm_comparison}). In each of the settings that we consider, we vary the time horizon $T$, and consider a budget of $B = T/2 \times \mathbbm{1}_{d \times 1}$ where $d$ is the number of resources. We note that this starting budget leads to the worst-case regret scaling for the instances with gaps which we consider. Overall, our simulation results confirm our theoretical predictions, including the importance of the conservativeness with respect to gaps principle, and demonstrate superior numerical performance of the {\sf RAMS} algorithm.
}

\simulation{
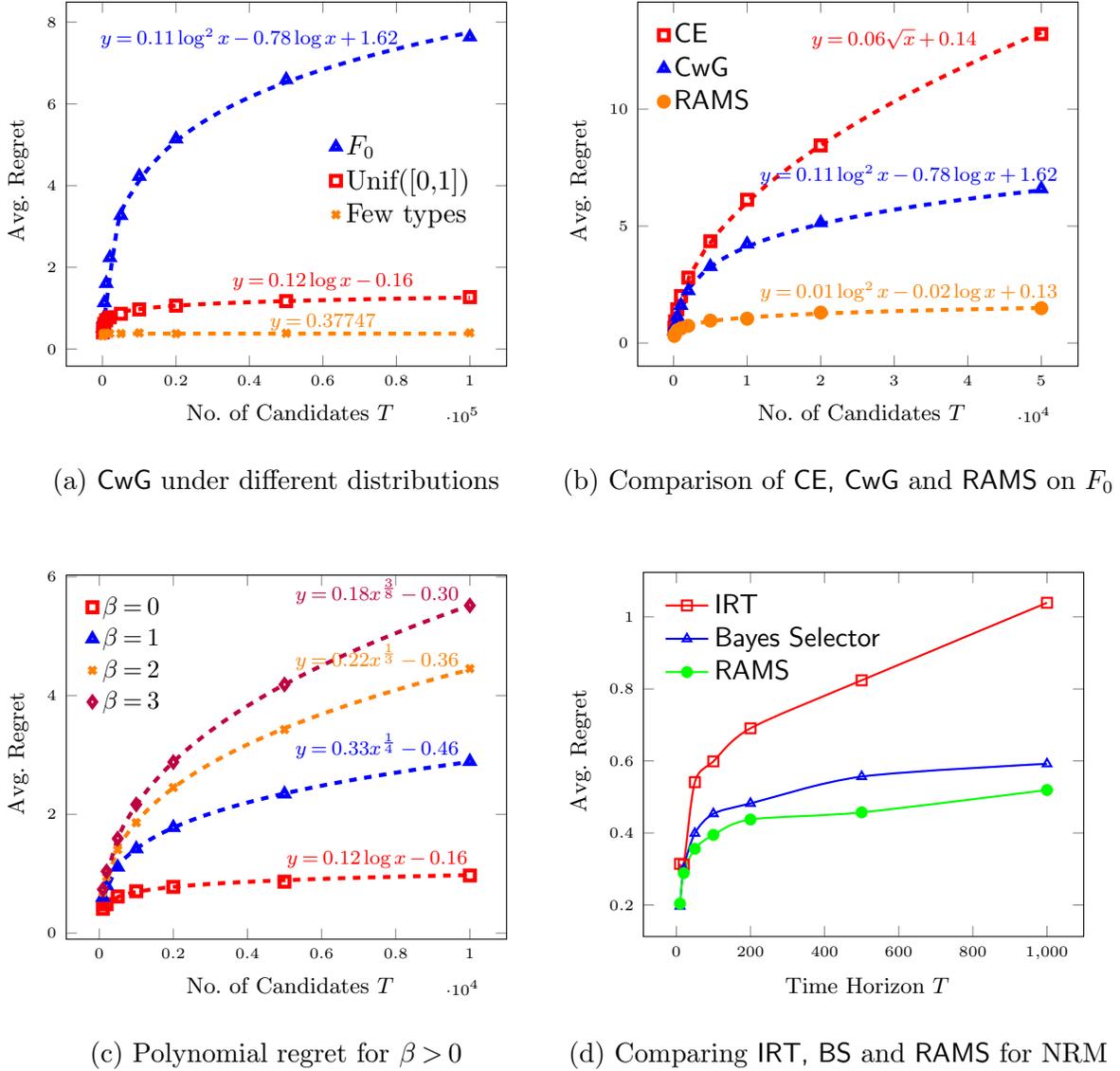
\begin{figure}[ht]
    \centering
    \begin{subfigure}{0.47\linewidth}
        \pgfplotsset{ width=1\linewidth,
        ylabel near ticks,
                 legend style={
                 draw = none,
                 fill=none,
                 cells={anchor=west},
                 legend pos=north west,
                 font = },
                 label style={font=\footnotesize},
                 tick label style={font=\tiny},
                }
    \begin{tikzpicture}
        \begin{axis}[legend cell align = {left}, legend style={at={(axis cs:60000,2.75)},anchor=south west}, xlabel = No. of Candidates $T$,ylabel = Avg. Regret, domain=100:100000]
        
        \addplot[scatter, only marks, ultra thick, blue, mark = triangle, scatter/classes = {f0 = {}}] file[skip first] {revision_plot_data_files/ccr_beta_0.dat};
        
        \addplot[ scatter, only marks, ultra thick, red, mark = square, scatter/classes = {unif = {}}] file[skip first] {revision_plot_data_files/uniform.dat};
        
        \addplot[scatter, only marks, ultra thick, orange, mark = x, scatter/classes = {discrete83 = {}}] file[skip first] {revision_plot_data_files/discrete_epsilon_83.dat};
        
        \addplot[smooth, dashed, blue, ultra thick] {0.1133824 * ln(x) * ln(x) - 0.77297728 * ln(x) + 1.62 };
        \node[rotate = 0] at (axis cs: 40000, 7.6) {\scriptsize \textcolor{blue}{$y = 0.11 \log^2 x - 0.78 \log x + 1.62$}};
        
        \addplot[smooth, dashed, red, ultra thick] {0.12351855 * ln(x) - 0.16194357};
        \node[] at (axis cs: 60000, 1.65) {\scriptsize \textcolor{red}{$y = 0.12 \log x - 0.16$}};

        \addplot[smooth, dashed, orange, ultra thick] {0.37747};
        \node[] at (axis cs: 60000, 0.65) {\scriptsize \textcolor{orange}{$y = 0.37747$}};
        
        \legend{$F_0$, Unif(\text{[0,1]}), $\text{Few types}$}
    \end{axis}
    \end{tikzpicture}
    \caption{{\sf CwG} under different distributions}
    \label{fig:cwg_diff_dists}
    \end{subfigure}
    \begin{subfigure}{0.47\linewidth}
        \pgfplotsset{ width=1\linewidth,
        ylabel near ticks,
                 legend style={
                 draw = none,
                 fill=none,
                 cells={anchor=west},
                 legend pos=north west,
                 font = },
                 label style={font=\footnotesize},
                 tick label style={font=\tiny},
                }
    \begin{tikzpicture}
        \begin{axis}[ legend cell align = {left}, legend pos = north west, xlabel = No. of Candidates $T$,ylabel = Avg. Regret, domain = 100:50000]
        \addplot[scatter, only marks ,ultra thick, red, mark = square, scatter/classes = {ce = {}}] file[skip first] {revision_plot_data_files/ce_beta_gap_0.dat};
        
        \addplot[scatter, only marks, ultra thick,blue, mark = triangle, scatter/classes = {cwg = {}}] file[skip first] {revision_plot_data_files/cwg_beta_0.dat};
        
        \addplot[scatter, only marks, ultra thick, orange, mark = *, scatter/classes = {rams = {}}]file[skip first] {revision_plot_data_files/ms_gap_beta_0_sim.dat};

        \addplot[red, dashed, ultra thick]{0.0587725 * sqrt(x) +  0.14013};
        \node[] at (axis cs: 30000, 13) {\scriptsize \textcolor{red}{$y = 0.06\sqrt{x} + 0.14$}};

        \addplot[blue, dashed, ultra thick]{0.1133824 * ln(x) * ln(x) - 0.77297728 * ln(x) + 1.62 };
        \node[] at (axis cs: 32000, 7.2) {\scriptsize \textcolor{blue}{$y = 0.11\log^2 x - 0.78 \log x + 1.62$}};

        \addplot[orange, dashed, ultra thick]{0.01388464 * ln(x) * ln(x) -0.02373132 * ln(x) +  0.12961871 };
        \node[] at (axis cs: 32000, 2.2) {\scriptsize \textcolor{orange}{$y = 0.01\log^2 x - 0.02 \log x + 0.13$}};
        
        \legend{{\sf CE} , {\sf CwG} , \textsf{RAMS}}
    \end{axis}
    \end{tikzpicture}
    \caption{Comparison of {\sf CE, CwG} and {\sf RAMS} on $F_0$}
    \label{fig:cwg_ce_comparison}
    \end{subfigure}
    
    \vspace{2em}
    
    \begin{subfigure}{0.47\linewidth}
    \pgfplotsset{ width=1\linewidth,
        ylabel near ticks,
                 legend style={
                 draw = none,
                 fill=none,
                 cells={anchor=west},
                 legend pos=north west,
                 font = \small},
                 label style={font=\footnotesize},
                 tick label style={font=\tiny},
                }
         \begin{tikzpicture}
        \begin{axis}[xlabel = No. of Candidates $T$ ,ylabel = Avg. Regret, domain = 100:10000
        ]

        \addplot[scatter, only marks, ultra thick, red, mark = square, scatter/classes = {0 = {}}] file[skip first] {revision_plot_data_files/ce_beta_0.dat};

        \addplot[scatter, only marks, ultra thick, mark = triangle, blue , scatter/classes = {1 = {}}] file[skip first] {revision_plot_data_files/ce_beta_1.dat};
        
        \addplot[scatter, only marks, ultra thick, orange, mark = x, scatter/classes = {2 = {}}] file[skip first] {revision_plot_data_files/ce_beta_2.dat};
    
        \addplot[scatter, only marks, ultra thick, purple, mark = diamond, scatter/classes = {3 = {}}] file[skip first] {revision_plot_data_files/ce_beta_3.dat};

        \addplot[smooth, dashed, red, ultra thick] {0.12351855 * ln(x) - 0.16194357};
        \node[] at (axis cs: 7500, 1.25) {\scriptsize \textcolor{red}{$y = 0.12 \log x - 0.16$}};

        \addplot[smooth, dashed, blue, ultra thick] {0.33414154 * x^(1/4) - 0.45982557};
        \node[] at (axis cs: 7500, 3.15) {\scriptsize \textcolor{blue}{$y = 0.33 x^{\frac{1}{4}} - 0.46$}};
        
        \addplot[smooth, dashed, orange, ultra thick] { 0.22252733 * x^(1/3) -0.3580183};
        \node[] at (axis cs: 7500, 4.65) {\scriptsize \textcolor{orange}{$y = 0.22 x^{\frac{1}{3}} - 0.36$}};

        \addplot[smooth, dashed, purple, ultra thick] { 0.18389254 * x^(3/8) -0.29846566};
        \node[] at (axis cs: 7500, 5.75) {\scriptsize \textcolor{purple}{$y = 0.18 x^{\frac{3}{8}} - 0.30$}};
        
        \legend{$\beta = 0$, $\beta = 1$,  $\beta = 2$ , $\beta = 3$};
    \end{axis}
    \end{tikzpicture}
    \caption{Polynomial regret for $\beta > 0$}
    \label{fig:beta-comparison}
    \end{subfigure}
    \begin{subfigure}{0.47\linewidth}
        \pgfplotsset{ width=1\linewidth,
        ylabel near ticks,
                 legend style={
                 draw = none,
                 fill=none,
                 cells={anchor=west},
                 legend pos=north west,
                 font = },
                 label style={font=\footnotesize},
                 tick label style={font=\tiny},
                }
    \begin{tikzpicture}
        \begin{axis}[ legend cell align = {left}, legend pos = north west, xlabel = Time Horizon $T$,ylabel = Avg. Regret]
        \addplot[smooth, thick, red, mark = square] file[skip first] {revision_plot_data_files/nrm_finite_irt.dat};
        \addplot[smooth, thick,blue, mark = triangle] file[skip first] {revision_plot_data_files/nrm_finite_bayes_selector.dat};
        \addplot[smooth,thick, green, mark = *]file[skip first] {revision_plot_data_files/nrm_finite_sim_constant_50.dat};
        \legend{{\sf IRT} , {\sf Bayes Selector} , \textsf{RAMS}}
    \end{axis}
    \end{tikzpicture}
    \caption{Comparing {\sf IRT, BS} and {\sf RAMS} for NRM}
    \label{fig:nrm_comparison}
    \end{subfigure}
    
    \caption{(a) Illustrates the performance of {\sf CwG} for different distributions, (b) compares the performance of the {\sf CE}, {\sf CwG} and {\sf RAMS} policies on $F_{0} = \text{Unif}([0,1/4] \cup [3/4,1])$, (c) highlights the polynomial regret scaling (with the exponent dependent on $\beta$) for the gapless variant of the $F_\beta$ distribution given in \eqref{eq:example-beta-epsilon-clustered-distribution}, (d) compares the performance of {\sf IRT}, {\sf Bayes Selector} and {\sf RAMS} for NRM with a few types.}
    \label{fig:numerical_comparison}
\end{figure}
}
\simulation{
\textbf{Figure \ref{fig:cwg_diff_dists}.~}We numerically study the regret scaling of the {\sf CwG} policy as a function of the time horizon $T$ for different distributions. The distributions we consider are: (i) bimodal uniform distribution $F_0 = \text{Uniform}([0,1/4] \cup [3/4,1])$, (ii) the uniform distribution over $[0,1]$ and (iii) a discrete distribution over a few types $\{0.25, 0.5, 0.75\}$ and the probability mass being $1/3$ for each of the points. We numerically evaluate the average regret for different number of candidates $T$ (with the budget varying as $B = T/2$) and fit a curve (as shown in the dashed lines) to observe the regret scaling. For each of the three distributions considered, we empirically observe that the regret scaling is consistent with our theoretical guarantees as implied by Corollary \ref{cor:cr-upper-bound-example} (log squared regret) for the bimodal distribution, Corollary \ref{cor:beta-one-clustered-regret} (logarithmic regret) for the uniform distribuion and Corollary \ref{cor:discrete-distribution-bound} (bounded regret) for the discrete distribution with few types.

}

\simulation{
\textbf{Figure \ref{fig:cwg_ce_comparison}.}~
We numerically study the average regret scaling of the {\sf CE}, {\sf CwG} and {\sf RAMS} policy for the bimodal uniform distribution $F_0 = \text{Unif}([0,1/4] \cup [3/4,1])$ with gap in the interval $[1/4,3/4]$. We fit a curve (as shown in dashed lines) to observe the regret scaling. For each of the three policies considered, we empirically observe that the regret scaling is consistent with our theoretical guarantees as implied by Proposition \ref{thm:ce-regret-lower-bound} for the {\sf CE} policy, Corollary \ref{cor:cr-upper-bound-example} for the {\sf CwG} policy and Corollary \ref{cor:cwg-rams-guarantee} for the {\sf RAMS} policy. While both {\sf CwG} and {\sf RAMS} have the same regret scaling, we observe that {\sf RAMS} has superior numerical performance over {\sf CwG} since {\sf RAMS} is designed to minimize the compensation and hence the regret, whereas {\sf CwG} is designed to optimize only the scaling of the compensation (and hence the regret scaling).

}

\simulation{
\textbf{Figure \ref{fig:beta-comparison}.}~
To assess the influence of the parameter $\beta$, we examine the performance of the {\sf CE} (equivalently {\sf CwG}) algorithm on the \emph{gapless} version of the $F_{\beta}$ distribution, as described in \eqref{eq:example-beta-epsilon-clustered-distribution} for $\beta \in \{0,1,2,3\}$. From Theorem \ref{thm:CwG-regret-lower-bound} and Corollary \ref{cor:beta-one-clustered-regret}, we know that {\sf CE} has the optimal regret scaling. We fit a curve (shown in dashed lines) to the empirical average regret for different values of time horizon $T$ and observe that the regret for $\beta \in \{1,2,3\}$ scales polynomially in the time horizon with the exponent given by $\frac{1}{2} - \frac{1}{2(1 + \beta)}$ and this is consistent with our guarantees in Corollary \ref{cor:beta-one-clustered-regret}.

}

\simulation{
\textbf{Figure \ref{fig:nrm_comparison}.~} We consider an NRM problem with two resources and six types. The types $\theta = (r_\theta, \bm{c}_\theta)$ are given as
    $\xi_{1} = (1.0, [1,0])$,  $\xi_{2} = (0.6, [1,0])$, $\xi_{3} = (1, [0,1])$,  $\xi_{4} = (0.5, [0,1])$, $\xi_{5} = (0.9, [1,1])$, $\xi_{6} = (0.8, [1,1])$.
The requests arrive i.i.d with $\mathbb{P}\left( \theta_{t} = \xi_{j}\right) = 0.2, \forall j \in \{1,2,3,4\}$ and $\mathbb{P}\left( \xi_{t} = \xi_{j} \right) = 0.1, \forall j \in \{5,6\}$. 
We compare the performance of {\sf RAMS} against two near optimal algorithms - {\sf Infrequent Resolving with Thresholding (IRT)} \citep{bumpensanti2020re} and {\sf Bayes Selector (BS)} \citep{vera2021bayesian}. We observe that for all the three algorithms that we consider, the regret increases initially but converges to a constant for sufficiently large $T$. We observe that amongst all the three algorithms considered, {\sf RAMS} either matches or improves upon the algorithms.
}

\section{Conclusion}
\label{sec:conclusion}
\new{
In this work, we considered dynamic resource allocation problems and investigated the impact of distributional assumptions on algorithmic performance. By focusing on the multisecretary problem, we gained valuable insights into the fundamental drivers and limits of algorithmic regret performance. We identified a novel driver of regret, characterized by the parameter $\beta$, which measures the concentration of types around gaps. We introduced the Conservativeness with respect to Gaps (CwG) principle, and used it to develop an innovative algorithmic approach that mitigates the limitations of the widely used certainty-equivalent (CE) policy. The CwG principle, along with its associated CwG algorithm, achieves near-optimal regret scaling of $\tilde{O}(T^{\frac{1}{2} - \frac{1}{2(1 + \beta)}})$ for a broad class of distributions with gaps parameterized by $\beta$. Furthermore, we analyzed the natural {\sf Repeatedly Act using Multiple Simulations (RAMS)} algorithm, which offers a general-purpose solution for online resource allocation problems (not just the multisecretary problem), which is applicable to any distribution of requests. {\sf RAMS} is practical and data-driven, relying on simulated future demand scenarios to drive decision making. Heuristically speaking, {\sf RAMS} is equivalent to a bid price control policy where the bid prices are computed by averaging the shadow prices of the hindsight optimal problem for multiple scenarios. This requires a minor adaptation of existing dual-based systems which is an industry default. 
}

\simulation{
Recently, there has been a growing interest in studying online resource allocation problems in the presence of horizon uncertainty \citep{besbes2014dynamic, balseiro2022online, bai2023fluid, aouad2022nonparametric}. 
Specifically, \cite{bai2023fluid} demonstrate that by leveraging an alternative fluid benchmark, it is possible to achieve a sublinear regret scaling of $\mathcal{O}(\sqrt{T})$, through the use of a static policy. Nevertheless, a na\"ive  implementation of the {\sf RAMS} approach yields regret (relative to the alternative fluid benchmark of \cite{bai2023fluid}) that scales linearly. Whether {\sf RAMS} can be adapted to attain sublinear regret remains unknown. We leave the exploration of this intriguing question, as well as other related queries surrounding the development of near-optimal algorithms under horizon uncertainty, for future endeavors.
}


\textbf{Acknowledgements.} 
YK and AK gratefully acknowledge the support of the National Science Foundation via grant  CMMI-1653477. We thank the reviewers, the associate editor and the area chair for their thought provoking comments which greatly improved the paper.

\bibliographystyle{informs2014}
\bibliography{references}

\ECSwitch


\ECHead{Appendix}

The Appendix is organized as follows.
Appendix \ref{onapp:universal-lower-bound} provides the proof of the universal lower bound in Theorem \ref{thm:CwG-regret-lower-bound}.
Appendix \ref{onapp:analysis-cwg-policy} provides the proof of Theorem \ref{thm:general-cwg-upper-bound}.
Appendix \ref{sec:proof-cwg-discrete-type} and Appendix \ref{app:proof-corollary-beta-one-clustered} provides the proof of Corollaries \ref{cor:discrete-distribution-bound} and \ref{cor:beta-one-clustered-regret} respectively. 
Appendix \ref{onapp:recovering-existing-guarantees-rams} discusses how {\sf RAMS} is able to recover both the new and prior regret guarantees in the context of the multisecretary problem and the more general network revenue management and online matching problems.
Appendix \ref{onapp:details-rams} provides the proof of Theorem \ref{thm:meta-performance-rams} and the corollaries \ref{cor:cwg-rams-guarantee}, \ref{cor:nrm-regret-profiles} and \ref{cor:online-matching-constant}.
Appendix \ref{onapp:fulfillment-problem-mapping} provides a discussion of the connections between the order fulfillment problem and the multisecretary problem.
Appendix \ref{onapp:details-beta-varepsilon-delta-clustered} provides some details on the different possible clusterings for $(\beta, \varepsilon_0, \delta)$-clustered distributions.

\section{Proof of Theorem \ref{thm:CwG-regret-lower-bound}}
\label{onapp:universal-lower-bound}
\begin{figure}[h]
    \centering
    \begin{tikzpicture}
        \draw[-, thick] (0,0) -- (16,0);
        \draw[-, thick] (8,0.15) -- (8,-0.15);
        \node at (8,-0.5) {$\frac{1}{2}$};
        \draw[-, thick] (7.5, 0.15) -- (7.5, -0.15);
        \draw[-, thick] (8.5, 0.15) -- (8.5, -0.15);
        \draw[<->, thick] (7.5,0.25) -- (8.5,0.25);
        \node at (7.5, -0.5) {$\ell$};
        \node at (8.5,-0.5) {$u$};
        \node at (8, 0.5) {$g$};
        \draw[-,thick] (6,0.15) -- (6,-0.15);
        \draw[-,thick] (0,0.15) -- (0,-0.15);
        \draw[-, thick] (10,0.15) -- (10,-0.15);
        \draw[-, thick] (16,0.15) -- (16,-0.15);
        \node at (6.2,-0.5) {$\ell_{1}$};
        \node at (10, -0.5) {$u_{1}$};
        \draw[<->, thick] (6, 0.25) -- (7.5, 0.25);
        \draw[<->, thick] (8.5, 0.25) -- (10, 0.25);
        \node at (6.75, 0.5) {$\Delta_{\beta}$};
        \node at (9.25, 0.5) {$\Delta_{\beta}$};
        \draw[-, thick] (5,0.15) -- (5, -0.15);
        \draw[-, thick] (11,0.15) -- (11, -0.15);
        \draw[<->, thick] (5, 0.25) -- (6, 0.25);
        \draw[<->, thick] (10,0.25) -- (11, 0.25);
        \node at (5.5, 0.5) { $\tilde{\Delta}_{\beta}$};
        \node at (10.5, 0.5) { $\tilde{\Delta}_{\beta}$};
        \node at (5, -0.5) {$\ell_{2}$};
        \node at (11, -0.5) {$u_{2}$};
        \node at (0,-0.5) {$0$};
        \node at (16,-0.5) {$1$};
        \draw [thick, decoration={brace, mirror, raise=0.3cm}, decorate] (0,-0.5) -- (5,-0.5); 
        \node at (2.5, -1.2) {$\mathcal{I}_{L}$};
        \draw [thick, decoration={brace, mirror, raise=0.3cm}, decorate] (5,-0.5) -- (6,-0.5);
        \node at (5.5, -1.2) {$\mathcal{I}_{M_{1}}$};
        \draw [thick, decoration={brace, mirror, raise=0.3cm}, decorate] (6,-0.5) -- (7.5,-0.5);
        \node at (6.75, -1.2) {$\mathcal{I}_{M_{2}}$};
        \draw [thick, decoration={brace, mirror, raise=0.3cm}, decorate] (8.5,-0.5) -- (10,-0.5); 
        \node at (9.25, -1.2) {$\mathcal{I}_{M_{3}}$};
        \draw [thick, decoration={brace, mirror, raise=0.3cm}, decorate] (10,-0.5) -- (11,-0.5);
        \node at (10.5, -1.2) {$\mathcal{I}_{M_{4}}$};
        \draw [thick, decoration={brace, mirror, raise=0.3cm}, decorate] (11,-0.5) -- (16,-0.5);
        \node at (13.5, -1.2) {$\mathcal{I}_{H}$};
    \end{tikzpicture}
    \caption{Partition of the set $\mathcal{S} = [0,\ell] \cup [u,1]$ into disjoint set $\mathcal{I}_{L} = [0, \ell_{2}), \mathcal{I}_{M_{1}} = [\ell_{2}, \ell_{1}), \mathcal{I}_{M_{2}} = [\ell_{1}, \ell], \mathcal{I}_{M_{3}} = [u, u_{1}), \mathcal{I}_{M_{4}} = [u_{1}, u_{2}), \mathcal{I}_{H} = [u_{2}, 1]$, where $\ell_{1} \triangleq \ell - \Delta_{\beta}, \ell_{2} \triangleq \ell - c_{0} \Delta_{\beta}$, $u_{1} \triangleq u + \Delta_{\beta}, u_{2} \triangleq u + c_{0} \Delta_{\beta}$ and $\Tilde{\Delta}_{\beta} \triangleq (c_{0} - 1) \Delta_{\beta}$.}
    \label{fig:interval-partition}
\end{figure}
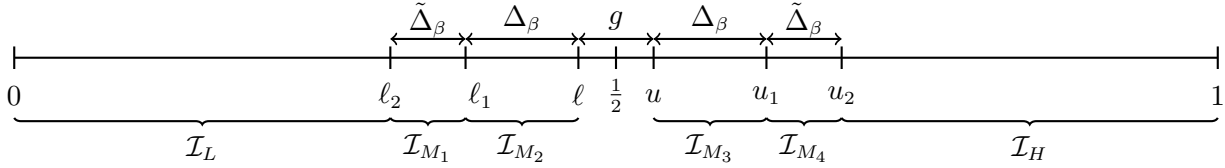

First we will consider the case of $\beta = 0$. For the uniform distribution over $[0,1]$, we have that $\beta = 0$ and from Proposition 4 of \cite{bray2022logarithmic}, Theorem \ref{thm:CwG-regret-lower-bound} follows for $\beta = 0$. Therefore our focus will on the case of $\beta > 0$. Fix $\beta >0$ and fix a number $g \geq 0$. In the context of Example \ref{ex:bimodal-distribution}, we have that $g = \frac{1}{2}$. Consider a distribution supported on the set $\mathcal{S} \triangleq \left[0, \ell \right] \cup \left[ u, 1\right]$ where $\ell \triangleq \frac{1}{2} - \frac{g}{2}$ and $u \triangleq \frac{1}{2} + \frac{g}{2}$. For $g = 1/2$, we have that $\ell = \frac{1}{4}$ and $u = \frac{3}{4}$. For a fixed $\beta > 0$ and $g, \ell, u$ as defined above, consider the following candidate ability distribution $F_{\beta, \ell, u}$,
    \begin{align}
        \label{eq:cdf-ability-distribution}
        F_{\beta, \ell, u}(x) = \begin{cases} -\frac{(\ell - x)^{1 + \beta }}{2\ell^{1 + \beta}} + \frac{1}{2} , & \quad 0 \leq x \leq \ell \\ 
        \frac{1}{2}, & \quad \ell \leq x \leq u \\
        \frac{(x - u)^{1 + \beta}}{2(1 - u)^{1 + \beta}} + \frac{1}{2}, & \quad u \leq x \leq 1
        \end{cases}
    \end{align}
    For $g > 0$, we can easily verify that $F_{\beta, \ell, u}$ is a $\left(\beta, \varepsilon_{0} = \frac{1}{2}\right)$-clustered distribution and for $g = 0$, $F_{\beta, \ell, u}$ is a $\left(\beta, \varepsilon_{0} = 1 \right)$-clustered distribution.
    Next, we will fix the time horizon $T > 0$ and set the budget $B \triangleq \floor{\frac{1}{2}T}$. Define $\Delta_{\beta} \triangleq T^{-\frac{1}{2(1 + \beta)}}(1 - u) = T^{-\frac{1}{2(1 + \beta)}} \ell$. Define $c_{0} \triangleq \left(\frac{129}{128} \right)^{\frac{1}{1 + \beta}} > 1$. Define the following quantities:
    \begin{align}
        \label{eq:notation}
        \alpha_{0} \triangleq c_{0} - 1 > 0, \ \ \Tilde{\Delta}_{\beta} \triangleq \alpha_{0} \Delta_{\beta}, \ \ \ell_{1} \triangleq \ell - \Delta_{\beta}, \ \ \ell_{2} \triangleq \ell - c_{0} \Delta_{\beta}, \ \ u_{1} \triangleq u + \Delta_{\beta}, \ \ u_{2} \triangleq u + c_{0} \Delta_{\beta}.
    \end{align}
    Now we will partition the set $\mathcal{S} \triangleq [0, \ell] \cup [u, 1]$ into the following sets (refer to Figure \ref{fig:interval-partition}):
    \begin{align*}
        \mathcal{I}_{L} = [0, \ell_{2}), \ \ \mathcal{I}_{M_{1}} = [\ell_{2}, \ell_{1}), \ \ \mathcal{I}_{M_{2}} = [\ell_{1}, \ell], \ \ \mathcal{I}_{M_{3}} = [u, u_{1}), \ \ \mathcal{I}_{M_{4}} = [u_{1}, u_{2}), \ \ \mathcal{I}_{H} = [u_{2}, 1]
    \end{align*}
    Further define the sets $\mathcal{I}_{M_{c}} \triangleq \mathcal{I}_{M_{2}} \cup \mathcal{I}_{M_{3}}$ and $\mathcal{I}_{M_{p}} \triangleq \mathcal{I}_{M_{1}} \cup \mathcal{I}_{M_{4}}$.

 Let $\vect{\theta}_{\geq 1}$ denote the sequence of candidate abilities and define $N(\mathcal{A}, t_{1}, t_{2})$ denote the number of candidate abilities in the set $\mathcal{A}$ that arrive in the time interval $[t_{1}, t_{2}]$. Formally, the random variable $N(\mathcal{A}, t_{1}, t_{2})$ is defined as 
 \begin{align}
     \label{eq:num-arrivals}
     N(\mathcal{A}, t_{1}, t_{2}) \triangleq \sum_{k = t_{1}}^{t_{2}} \mathbbm{1}\{\theta_{t} \in \mathcal{A}\}, \quad \forall \mathcal{A} \subseteq \mathcal{S}, t_{1}, t_{2} \in \{1,2, \dots, T\}
 \end{align}
 Let $\mu_{t_{1}}^{t_{2}}(A) \triangleq \mathbb{E}\left[N(\mathcal{I}_{H}, t_{1}, t_{2}) \right]$ denote the mean of the random variable $N(\mathcal{I}_{H}, t_{1}, t_{2})$.
    Next we define the following set of events:
    \begin{align}
        \label{event:high-type-first-half-low}
        \mathcal{H}_{1} &\triangleq \bigg\{\frac{T}{4} - \frac{\sqrt{T}}{2} \leq N\left(\mathcal{I}_{H}, 1, B\right) \leq \frac{T}{4}\bigg\} 
        \\
        \label{event:high-type-second-half-low}
        \mathcal{H}_{2} &\triangleq \bigg\{\frac{T}{4} - 4\sqrt{T} \leq N\left(\mathcal{I}_{H}, B + 1, T \right) \leq \frac{T}{4} - 3 \sqrt{T} \bigg\}, \\
        \label{event:high-type-second-half-high}
        \tilde{\mathcal{H}}_{2} &\triangleq \bigg\{\frac{T}{4} + \frac{\sqrt{T}}{2} \leq N\left(\mathcal{I}_{H}, B + 1, T \right) \leq \frac{T}{4} + \frac{3\sqrt{T}}{2} \bigg\}, \\
        \label{event:middle-center-type-first-half}
        \mathcal{C}_{1} &\triangleq \bigg\{\frac{\sqrt{T}}{4} \leq N(\mathcal{I}_{M_{c}}, 1, B) \leq \sqrt{T} \bigg\}, \\
        \label{event:middle-center-type-second-half}
        \mathcal{C}_{2} &\triangleq \bigg\{ \frac{\sqrt{T}}{4} \leq N(\mathcal{I}_{M_{c}}, B + 1, T)  \leq \sqrt{T} \bigg\}, \\ 
        \label{event:middle-periphery-type-first-half}
        \mathcal{P}_{1} &\triangleq \bigg\{\frac{\sqrt{T}}{256} \leq N(\mathcal{I}_{M_{p}}, 1, B) \leq \frac{\sqrt{T}}{64} \bigg\}, \\
        \label{event:middle-periphery-type-second-half}
        \mathcal{P}_{2} &\triangleq \bigg\{ \frac{\sqrt{T}}{256} \leq N(\mathcal{I}_{M_{p}}, B + 1, T)  \leq \frac{\sqrt{T}}{64} \bigg\}.
    \end{align}
    Further we define the events $\mathcal{H} \triangleq \mathcal{H}_{1} \cap \mathcal{H}_{2}, \Tilde{\mathcal{H}} \triangleq \mathcal{H}_{1} \cap \tilde{\mathcal{H}}_{2},  \mathcal{C} \triangleq \mathcal{C}_{1} \cap \mathcal{C}_{2}$ and $\mathcal{P} \triangleq \mathcal{P}_{1} \cap \mathcal{P}_{2}$.

    \paragraph{Discussion of the Hindsight Optimal.} Conditional on the event $\Tilde{\mathcal{H}} \cap \mathcal{C} \cap \mathcal{P}$, we have that the total number of arrivals in the set $\mathcal{I}_{H}$ is more than the budget $B$ i.e. $N(\mathcal{I}_{H}, 1,T) \geq \frac{1}{2}T \geq B$ and hence the hindsight optimal must reject \emph{all} the arrivals in the set $\mathcal{I}_{L} \cup \mathcal{I}_{M_{p}} \cup \mathcal{I}_{M_{c}}$ and \emph{possibly some} arrivals in the set $\mathcal{I}_{H}$. However, conditional on the event $\mathcal{H} \cap \mathcal{C} \cap \mathcal{P}$, we have that total number of arrivals in the set $\mathcal{I}_{H} \cup \mathcal{I}_{M_{c}} \cup \mathcal{I}_{M_{p}}$ is less than the budget $B$ i.e, $N(\mathcal{I}_{H}, 1, T) + N(\mathcal{I}_{M_{c}}, 1, T) + N(\mathcal{I}_{M_{p}},1,T) \leq \frac{1}{2}T - \frac{31}{32}\sqrt{T} < \floor{\frac{1}{2}T} = B$ for sufficiently large $T$ and hence the hindsight optimal must accept \emph{all} the arrivals in the set $\mathcal{I}_{H} \cup \mathcal{I}_{M_{p}} \cup \mathcal{I}_{M_{c}}$ and \emph{possibly some} arrivals in the set $\mathcal{I}_{L}$.

    Let $N^{\textsf{DP}}(\mathcal{A}, t_{1}, t_{2})$ denote the number of accepted candidates by the {\sf DP} (optimal dynamic programming policy) with ability in the set $\mathcal{A}$ and they arrive in the time interval $[t_{1}, t_{2}]$ which we formally define as:
    \begin{align}
        N^{\textsf{DP}}(\mathcal{A}, t_{1}, t_{2}) \triangleq \sum_{k = t_{1}}^{t_{2}} \mathbbm{1}\{\theta_{t} \in \mathcal{A}, \pi_{t}^{\textsf{DP}} = \text{accept}\}, \quad \forall \mathcal{A} \subseteq \mathcal{S}, t_{1}, t_{2} \in \{1,2,\dots, T\}
     \end{align}
     Define the event $\mathcal{E}$ which says that under the optimal online policy, the number of accepted candidates up till time $B (= \floor{\frac{1}{2}T})$ is at least one eighth of the number of arrivals in set $\mathcal{I}_{M_{c}}$ up till time $B$, i.e,
     \begin{align}
         \label{event:dp-mistakes}
         \mathcal{E} \triangleq \bigg\{ N^{\textsf{DP}}(\mathcal{I}_{M_{c}}, 1, B) \geq \frac{N(I_{M_{c}},1,B)}{8} \bigg\}
     \end{align}
     \paragraph{Proof Strategy.} Our proof will proceed by considering the following events: (a) $\mathcal{E} \cap \tilde{\mathcal{H}} \cap \mathcal{C} \cap \mathcal{P}$ and (b) $\mathcal{E}^{c} \cap \mathcal{H} \cap \mathcal{C} \cap \mathcal{P}$. In case (a), from the discussion about the hindsight optimal policy, the hindsight optimal policy will reject all the arrivals in the set $\mathcal{I}_{M_{c}}$ but DP will accept at least $\frac{1}{32}\sqrt{T}$ arrivals in set $\mathcal{I}_{M_{c}}$ in the time interval $[1,B]$. This will result in the DP \emph{incorrectly} rejecting at least $\frac{1}{32}\sqrt{T}$ arrivals in interval $\mathcal{I}_{H}$ and the cost of each of these mistakes is at least $\alpha_{0} \Delta_{\beta}$. In case (b), from the discussion about the hindsight optimal policy, the hindsight optimal policy will accept all the arrivals in the interval $\mathcal{I}_{M_{c}}$ but the DP accepts at most $\frac{1}{8}\sqrt{T}$ arrivals in the interval $\mathcal{I}_{M_{c}}$. This implies that at least $\frac{1}{8}\sqrt{T}$ arrivals in the interval $\mathcal{I}_{M_{c}}$ are \emph{incorrectly} rejected. This will result in the DP \emph{incorrectly} accepting at least $\frac{1}{8}\sqrt{T}$ arrivals in the interval $\mathcal{I}_{L}$ and the cost of each of these mistakes is again at least $\alpha_{0} \Delta_{\beta}$. Informally speaking, we can lower bound the expected regret as 
     \begin{align*}
         \text{Regret}(B, T; \textsf{DP}) \geq c\left( \mathbb{P}(\mathcal{E} \cap \tilde{\mathcal{H}} \cap \mathcal{C} \cap \mathcal{P} ) + \mathbb{P}\left(\mathcal{E}^{c} \cap \mathcal{H} \cap \mathcal{C} \cap \mathcal{P} \right) \right)\left[ \left(\text{\# of mistakes} \right) \times \left( \text{cost/mistake}\right)\right]
     \end{align*}
     Assuming we can show that $\mathbb{P}\left(\mathcal{E} \cap \tilde{\mathcal{H}} \cap \mathcal{C} \cap \mathcal{P} \right) + \mathbb{P}\left(\mathcal{E}^{c} \cap \mathcal{H} \cap \mathcal{C} \cap \mathcal{P} \right) \geq \gamma > 0$ for some $\gamma \in (0,1)$, we have that $\text{\# of mistakes are } \Omega(\sqrt{T})$ and cost of each mistake is $\Omega\left( T^{-\frac{1}{2(1 + \beta)}}\right)$. Combining all this will provide the lower bound guarantee as desired for $\beta > 0$. 

     Consider the random variable $\Lambda(B,T; \textsf{DP})$
     \begin{align}
         \Lambda(B,T; \textsf{DP}) &= \sum_{t = 1}^{T} \theta_{t} a_{t}^{\textsf{hs}} - \sum_{t = 1}^{T} \theta_{t} a_{t}^{\textsf{DP}}
     \end{align}

     Next we will formalize our proof strategy using the following two lemmas.

     \begin{lemma}
         \label{lem:conditional-regret-lb-1}
         Consider the event $\mathcal{E} \cap \Tilde{\mathcal{H}} \cap \mathcal{C} \cap \mathcal{P}$, then we have that $$\mathbb{E}\left[\Lambda(B,T; \textsf{DP}) | \mathcal{E} \cap \Tilde{\mathcal{H}} \cap \mathcal{C} \cap \mathcal{P} \right] \geq \frac{\alpha_{0} \ell}{32} T^{\frac{1}{2} - \frac{1}{2(1 + \beta)}},$$
         where $\alpha_{0} = \left(\frac{129}{128} \right)^{\frac{1}{1 + \beta}} - 1$ defined in \eqref{eq:notation} and $\ell = \frac{1}{2} - \frac{g}{2}$.
     \end{lemma}

     \begin{lemma}
         \label{lem:conditional-regret-lb-2}
         Consider the event $\mathcal{E}^{c} \cap \mathcal{H} \cap \mathcal{C} \cap \mathcal{P}$, then we have that
         $$\mathbb{E}\left[\Lambda(B,T; \textsf{DP}) | \mathcal{E}^{c} \cap \mathcal{H} \cap \mathcal{C} \cap \mathcal{P} \right] \geq \frac{\alpha_{0}}{8} T^{\frac{1}{2} - \frac{1}{2(1 + \beta)}},$$
         where $\alpha_{0} = \left(\frac{129}{128} \right)^{\frac{1}{1 + \beta}} - 1$ defined in \eqref{eq:notation}.
     \end{lemma}
     We defer the proofs of Lemmas \ref{lem:conditional-regret-lb-1} and \ref{lem:conditional-regret-lb-2} to Sections \ref{subsec:proof-lemma-conditional-regret-lb-1} and \ref{subsec:proof-lemma-conditional-regret-lb-2} respectively.
     Finally, we have that 
     \begin{align}
        \nonumber
         \text{Regret}(B, T; \textsf{DP}) &\stackrel{(a)}= \mathbb{E}\left[\Lambda(B,T; \textsf{DP}) \right], \\
         \nonumber
         &\stackrel{(b)}\geq \mathbb{E}\left[\Lambda(B,T; \textsf{DP}) | \mathcal{E} \cap \Tilde{\mathcal{H}} \cap \mathcal{C} \cap \mathcal{P} \right] \mathbb{P}\left(\mathcal{E} \cap \Tilde{\mathcal{H}} \cap \mathcal{C} \cap \mathcal{P} \right) \\ \nonumber & \quad + \mathbb{E}\left[\Lambda(B,T; \textsf{DP}) | \mathcal{E}^{c} \cap {\mathcal{H}} \cap \mathcal{C} \cap \mathcal{P} \right] \mathbb{P}\left(\mathcal{E}^{c} \cap {\mathcal{H}} \cap \mathcal{C} \cap \mathcal{P} \right), \\
         \label{eq:regret-lower-bound}
         &\stackrel{(c)}\geq \frac{\alpha_{0}}{32} T^{\frac{1}{2} - \frac{1}{2(1 + \beta)}} \left(\mathbb{P}\left( \mathcal{E} \cap \Tilde{\mathcal{H}} \cap \mathcal{C} \cap \mathcal{P}\right) + \mathbb{P}\left( \mathcal{E}^{c} \cap {\mathcal{H}} \cap \mathcal{C} \cap \mathcal{P} \right) \right),
     \end{align}
     where (a) follows from the definition of (expected) regret, (b) follows from total law of expectations, (c) follows from Lemmas \ref{lem:conditional-regret-lb-1} and \ref{lem:conditional-regret-lb-2}.
     
     Observe that $\Tilde{\mathcal{H}} = \mathcal{H}_{1} \cap \Tilde{\mathcal{H}}_{2}, \mathcal{H} = \mathcal{H}_{1} \cap \mathcal{H}_{2}$, $\mathcal{C} = \mathcal{C}_{1} \cap \mathcal{C}_{2}$ and $\mathcal{P} = \mathcal{P}_{1} \cap \mathcal{P}_{2}$ and moreover $ \mathcal{C}_{1} \independent \mathcal{C}_{2}, \mathcal{P}_{1} \independent \mathcal{P}_{2}$ and $\mathcal{H}_{1} \independent \mathcal{H}_{2}, \Tilde{\mathcal{H}}_{2}$ since the events $\mathcal{H}_{1}, \mathcal{C}_{1}, \mathcal{P}_{1}$ only depend on the arrivals in the time interval $[1,B]$ i.e. $\{\theta_{k}\}_{k = 1}^{B}$ whereas the events $\mathcal{H}_{2}, \Tilde{\mathcal{H}}_{2}, \mathcal{C}_{2}$ and $\mathcal{P}_{2}$ only depend on the arrivals in the time interval $[B+1, T]$ i.e. $\{\theta_{k}\}_{k = B + 1}^{T}$ and the arrivals by assumption are \emph{i.i.d}. Additionally, the events $\mathcal{E}, \mathcal{E}^{c}$ also only depend on the arrivals in the interval $[1,B]$ and hence are independent of $\mathcal{H}_{2}, \Tilde{\mathcal{H}}_{2}, \mathcal{C}_{2}$ and $\mathcal{P}_{2}$. Therefore, we have that 
     \begin{align}
        \nonumber
         \mathbb{P}&\left( \mathcal{E} \cap \Tilde{\mathcal{H}} \cap \mathcal{C} \cap \mathcal{P} \right) + \mathbb{P}\left( \mathcal{E}^{c} \cap \mathcal{H} \cap \mathcal{C} \cap \mathcal{P} \right) \\ 
         \nonumber
         &\stackrel{(a)}= \mathbb{P}\left( \mathcal{E} \cap \mathcal{H}_{1} \cap \tilde{\mathcal{H}}_{2} \cap \mathcal{C}_{1} \cap \mathcal{C}_{2} \cap \mathcal{P}_{1} \cap \mathcal{P}_{2} \right) + \mathbb{P}\left( \mathcal{E}^{c} \cap \mathcal{H}_{1} \cap {\mathcal{H}}_{2} \cap \mathcal{C}_{1} \cap \mathcal{C}_{2} \cap \mathcal{P}_{1} \cap \mathcal{P}_{2} \right) \\
         \nonumber
         &\stackrel{(b)}= \mathbb{P}\left(\mathcal{E} \cap \mathcal{H}_{1} \cap \mathcal{C}_{1} \cap \mathcal{P}_{1} \right) \mathbb{P}\left( \tilde{\mathcal{H}}_{2} \cap \mathcal{C}_{2} \cap \mathcal{P}_{2}\right) + \mathbb{P}\left( \mathcal{E}^{c} \cap \mathcal{H}_{1} \cap \mathcal{C}_{1} \cap \mathcal{P}_{1}\right) \mathbb{P}\left(\mathcal{H}_{2} \cap \mathcal{C}_{2} \cap \mathcal{P}_{2} \right) \\
         \nonumber
         &\stackrel{(c)}\geq \min\bigg\{\mathbb{P}\left(\Tilde{\mathcal{H}}_{2} \cap \mathcal{C}_{2} \cap \mathcal{P}_{2} \right), \mathbb{P}\left(\mathcal{H}_{2} \cap \mathcal{C}_{2} \cap \mathcal{P}_{2} \right)\bigg\} \cdot \left(\mathbb{P}\left(\mathcal{E} \cap \mathcal{H}_{1} \cap \mathcal{C}_{1} \cap \mathcal{P}_{1} \right) + \mathbb{P}\left(\mathcal{E}^{c} \cap \mathcal{H}_{1} \cap \mathcal{C}_{1} \cap \mathcal{P}_{1} \right) \right), \\
         \label{eq:prob-lower-bound}
         &\stackrel{(d)}= \min\bigg\{\mathbb{P}\left(\Tilde{\mathcal{H}}_{2} \cap \mathcal{C}_{2} \cap \mathcal{P}_{2} \right), \mathbb{P}\left(\mathcal{H}_{2} \cap \mathcal{C}_{2} \cap \mathcal{P}_{2} \right)\bigg\} \mathbb{P}\left( \mathcal{H}_{1} \cap \mathcal{C}_{1} \cap \mathcal{P}_{1}\right),
     \end{align}
     where (a) follows from the definition of $\mathcal{H}, \Tilde{\mathcal{H}}, \mathcal{C}$ and $\mathcal{P}$, (b) follows from the fact that $\mathcal{E} \cap \mathcal{H}_{1} \cap \mathcal{C}_{1} \cap \mathcal{P}_{1} \independent \Tilde{\mathcal{H}}_{2} \cap \mathcal{C}_{2} \cap \mathcal{P}_{2}$ and $\mathcal{E}^{c} \cap \mathcal{H}_{1} \cap \mathcal{C}_{1} \cap \mathcal{P}_{1} \independent \mathcal{H}_{2} \cap \mathcal{C}_{2} \cap \mathcal{P}_{2}$ using the arguments presented previously, (c) follows trivially, (d) follows from the law of total probability. 
     
     Now it suffices to show to that there exists a constant $\alpha > 0$ independent of $T$ such that for all $T$ sufficiently large, we have that  $\mathbb{P}(\Tilde{\mathcal{H}}_{2} \cap \mathcal{C}_{2} \cap \mathcal{P}_{2}), \mathbb{P}\left(\mathcal{H}_{2} \cap \mathcal{C}_{2} \cap \mathcal{P}_{2} \right), \mathbb{P}\left(\mathcal{H}_{1} \cap \mathcal{C}_{1} \cap \mathcal{P}_{1} \right) \geq \alpha$. Using a CLT argument, one can easily see that $\mathbb{P}\left(\mathcal{H}_{1}\right), \mathbb{P}(\tilde{\mathcal{H}}_{2}), \mathbb{P}\left(\mathcal{H}_{2} \right) \geq \alpha^{\prime} > 0$ and $\mathbb{P}\left(\mathcal{C}_{1}\right), \mathbb{P}\left(\mathcal{C}_{2}\right), \mathbb{P}\left( \mathcal{P}_{1}\right), \mathbb{P}\left(\mathcal{P}_{2} \right) \geq \alpha^{\prime}$ for $\alpha^{\prime \prime} \neq \alpha^{\prime}$. However the events $\mathcal{H}_{1}, \mathcal{C}_{1}, \mathcal{P}_{1}$ (similarly $\mathcal{H}_{2}, \mathcal{C}_{2}, \mathcal{P}_{2}$ and $\Tilde{\mathcal{H}}_{2}, \mathcal{C}_{2}, \mathcal{P}_{2}$) are correlated and hence proving $\mathbb{P}(\Tilde{\mathcal{H}}_{2} \cap \mathcal{C}_{2} \cap \mathcal{P}_{2}), \mathbb{P}\left(\mathcal{H}_{2} \cap \mathcal{C}_{2} \cap \mathcal{P}_{2} \right), \mathbb{P}\left(\mathcal{H}_{1} \cap \mathcal{C}_{1} \cap \mathcal{P}_{1} \right) \geq \alpha$ requires a conditioning argument which we will illustrate now. We will argue this the event $\mathcal{H}_{1} \cap \mathcal{C}_{1} \cap \mathcal{P}_{1} $ and the exact same argument works for the events $\Tilde{\mathcal{H}}_{2} \cap \mathcal{C}_{2} \cap \mathcal{P}_{2}$ and ${\mathcal{H}}_{2} \cap \mathcal{C}_{2} \cap \mathcal{P}_{2}$. We have that 
     \begin{align*}
         \mathbb{P}\left(\mathcal{H}_{1} \cap \mathcal{C}_{1} \cap \mathcal{P}_{1} \right) &\stackrel{(a)}= \mathbb{P}\left(\mathcal{C}_{1} \cap \mathcal{P}_{1}  \right) \mathbb{P}\left( \mathcal{H}_{1} | \mathcal{C}_{1} \cap \mathcal{P}_{1}\right), \\
         &\stackrel{(b)}= \mathbb{P}\left( \mathcal{H}_{1} \right) - \mathbb{P}\left(\mathcal{H}_{1} | \left( \mathcal{C}_{1} \cap \mathcal{P}_{1} \right)^{c} \right) \mathbb{P}\left( \left(\mathcal{C}_{1} \cap \mathcal{P}_{1} \right)^{c}\right), \\
         &\stackrel{(c)}\geq \mathbb{P}\left( \mathcal{H}_{1} \right) - \left(\mathbb{P}\left( \mathcal{C}_{1}^{c} \right) + \mathbb{P}\left( \mathcal{P}_{1}^{c}\right) \right),
     \end{align*}
     where (a) follows from the definition of conditional probability, (b) follows from the law of total probability i.e. $\mathbb{P}\left(\mathcal{H}_{1}\right) = \mathbb{P}\left(\mathcal{H}_{1} | \mathcal{C}_{1} \cap \mathcal{P}_{1} \right) \mathbb{P}\left(\mathcal{C}_{1} \cap \mathcal{P}_{1} \right) + \mathbb{P}\left( \mathcal{H}_{1} | \left(\mathcal{C}_{1} \cap \mathcal{P}_{1} \right)^{c} \right) \mathbb{P}\left( \left( \mathcal{C}_{1} \cap \mathcal{P}_{1} \right)^{c} \right)$ and (c) follows from the fact that $\mathbb{P}\left( \mathcal{H}_{1} | \left(\mathcal{C}_{1} \cap \mathcal{P}_{1} \right)^{c} \right)\mathbb{P}\left( \left( \mathcal{C}_{1} \cap \mathcal{P}_{1}\right)^{c}\right) \leq \mathbb{P}\left( \left( \mathcal{C}_{1} \cap \mathcal{P}_{1}\right)^{c}\right) \leq \mathbb{P}\left( \mathcal{C}_{1}^{c} \right) + \mathbb{P}\left( \mathcal{P}_{1}^{c} \right)$ where the first inequality follows from the fact that $\mathbb{P}\left( \mathcal{H}_{1} | \left(\mathcal{C}_{1} \cap \mathcal{P}_{1} \right)^{c} \right) \leq 1$ and the second inequality follows from the union bound. Using the exact same arguments we have that 
     \begin{align*}
         \mathbb{P}\left( \mathcal{H}_{2} \cap \mathcal{C}_{2} \cap \mathcal{P}_{2} \right) &\geq \mathbb{P}\left( \mathcal{H}_{2} \right) - \left(\mathbb{P}\left( \mathcal{C}_{2}^{c}\right) + \mathbb{P}\left( \mathcal{P}_{2}^{c} \right)\right), \\
         \mathbb{P}( \Tilde{\mathcal{H}}_{2} \cap \mathcal{C}_{2} \cap \mathcal{P}_{2}) &\geq \mathbb{P}(\Tilde{\mathcal{H}}_{2}) - \left(\mathbb{P}\left( \mathcal{C}_{2}^{c}\right) + \mathbb{P}\left( \mathcal{P}_{2}^{c} \right)\right).
     \end{align*}
     Next we present a few lemmas which would imply that $\mathbb{P}(\Tilde{\mathcal{H}}_{2} \cap \mathcal{C}_{2} \cap \mathcal{P}_{2} ), \mathbb{P}\left( \mathcal{H}_{2} \cap \mathcal{C}_{2} \cap \mathcal{P}_{2} \right), \mathbb{P}( \tilde{\mathcal{H}}_{2} \cap \mathcal{C}_{2} \cap \mathcal{P}_{2}) \geq 0.001$.
     \begin{lemma}
        \label{lem:high-type-concentration}
         There exists $T_{0} < \infty$ such that for all $T \geq T_{0}$, we have that $\mathbb{P}\left( \mathcal{H}_{1}\right), \mathbb{P}\left( \mathcal{H}_{2}\right), \mathbb{P}( \tilde{\mathcal{H}}_{2} ) \geq 0.003$ 
     \end{lemma}
     \begin{lemma}
        \label{lem:mid-type-concentration}
         There exists $T_{0} < \infty$ such that for all $T \geq T_{0}$, we have that $\mathbb{P}( \mathcal{C}_{1}^{c}), \mathbb{P}( \mathcal{C}_{2}^{c}), \mathbb{P}( {\mathcal{P}}_{1}^{c} ), \mathbb{P}( \mathcal{P}_{2}^{c} ) \leq 0.001$ 
     \end{lemma}
     We defer the proofs of Lemma \ref{lem:high-type-concentration} and \ref{lem:mid-type-concentration} to Appendix \ref{subsec:high-type-concentration} and \ref{subsec:niddle-type-concentration} respectively. Using Lemmas \ref{lem:high-type-concentration} and \ref{lem:mid-type-concentration} and \eqref{eq:prob-lower-bound}, we have that $\mathbb{P}( \mathcal{E} \cap \Tilde{\mathcal{H}} \cap \mathcal{C} \cap \mathcal{P} ) + \mathbb{P}\left( \mathcal{E}^{c} \cap \mathcal{H} \cap \mathcal{C} \cap \mathcal{P} \right) \geq 10^{-6}$, combined with \eqref{eq:regret-lower-bound} concludes the proof. \QEDB

\subsection{Proof of Lemma \ref{lem:conditional-regret-lb-1}}
\label{subsec:proof-lemma-conditional-regret-lb-1}
Recall the definition of the random variable $\Lambda(B,T;\textsf{DP}) = \sum_{k = 1}^{T} \theta_{k} \pi_{k}^{\sf hs} - \sum_{k = 1}^{T} \theta_{k}$. For a sequence of candidate ability arrivals $\bm{\theta}_{\geq 1}$, we can define the following random set of indices $\mathcal{J}^{\textsf{hs}}$ and $\mathcal{J}^{\textsf{DP}}$ as
     \begin{align}
        \label{def:random-indices-set}
         \mathcal{J}^{\textsf{hs}}(\mathcal{A}) &\triangleq \{k : \theta_{k} \in \mathcal{A} \text{ and } \pi_{k}^{\textsf{hs}} = 1\}, \ \ \mathcal{J}^{\textsf{DP}}(\mathcal{A}) \triangleq \{k : \theta_{k} \in \mathcal{A} \text{ and } \pi_{k}^{\textsf{DP}} = 1\}, \quad \forall \mathcal{A} \subseteq \mathcal{S}
     \end{align}
     Notice that we can equivalently write the sum of values chosen under the hindsight optimal and the DP policy as 
     \begin{align}
        \label{eq:offline-cost-decomposition}
         \sum_{t = 1}^{T} \theta_{t} \pi_{t}^{\textsf{hs}} &= \sum_{k \in \mathcal{J}^{\textsf{hs}}(\mathcal{I}_{L})} \theta_{k} + \sum_{k \in \mathcal{J}^{\textsf{hs}}(\mathcal{I}_{M_{p}})} \theta_{k} + \sum_{k \in \mathcal{J}^{\textsf{hs}}(\mathcal{I}_{M_{c}})} \theta_{k} + \sum_{k \in \mathcal{J}^{\textsf{hs}}(\mathcal{I}_{H})} \theta_{k} \\
         \label{eq:dp-cost-decomposition}
         \sum_{t = 1}^{T} \theta_{t} \pi_{t}^{\textsf{DP}} &= \sum_{k \in \mathcal{J}^{\textsf{DP}}(\mathcal{I}_{L})} \theta_{k} + \sum_{k \in \mathcal{J}^{\textsf{DP}}(\mathcal{I}_{M_{p}})} \theta_{k} + \sum_{k \in \mathcal{J}^{\textsf{DP}}(\mathcal{I}_{M_{c}})} \theta_{k} + \sum_{k \in \mathcal{J}^{\textsf{DP}}(\mathcal{I}_{H})} \theta_{k}
     \end{align}
     Now conditional on the event $\mathcal{E} \cap \tilde{\mathcal{H}} \cap \mathcal{C} \cap \mathcal{P}$, we have that $\mathcal{J}^{\textsf{hs}}(\mathcal{I}_{L}) = \mathcal{J}^{\textsf{hs}}(\mathcal{I}_{M_{p}}) = \mathcal{J}^{\textsf{hs}}(\mathcal{I}_{M_{c}}) = \emptyset$ and $|\mathcal{J}^{\textsf{hs}}(\mathcal{I}_{H})| = B$ and we have that $|\mathcal{J}^{\textsf{DP}}(\mathcal{I}_{M_{c}})| \geq \frac{1}{32}\sqrt{T}$. This follows from the fact under the event $\mathcal{E}$, the \textsf{DP} accepts at least $\frac{1}{8}N(\mathcal{I}_{M_{c}},1,B)$ and from the event $\mathcal{C}_{1}$, it follows that $N(\mathcal{I}_{M_{c}},1,B) \geq \frac{1}{4}\sqrt{T}$. Using this we have that 
     \begin{align}
     \label{eq:off-cost-simplified-1}
         \sum_{t = 1}^{T} \theta_{t} \pi_{t}^{\textsf{hs}} &= \sum_{k \in \mathcal{J}^{\textsf{hs}}(\mathcal{I}_{H})} \theta_{k} = \sum_{k \in \mathcal{J}^{\textsf{hs}}(\mathcal{I}_{H}) \backslash \mathcal{J}^{\textsf{DP}}(\mathcal{I}_{H})} \theta_{k} + \sum_{k \in \mathcal{J}^{\textsf{DP}}(\mathcal{I}_{H})} \theta_{k}
     \end{align}
     Since any online policy can select at most $B$ candidates and the offline policy will select the top $B$ candidates, we have that conditional on the event $\mathcal{E} \cap \tilde{\mathcal{H}} \cap \mathcal{C} \cap \mathcal{P}$,
     \begin{align}
        \label{eq:set-size-inequality}
        |\mathcal{J}^{\textsf{DP}}(\mathcal{I}_{H})| + |\mathcal{J}^{\textsf{hs}}(\mathcal{I}_{H}) \backslash \mathcal{J}^{\textsf{DP}}(\mathcal{I}_{H})| \geq |\mathcal{J}^{\textsf{DP}}(\mathcal{I}_{H})| + |\mathcal{J}^{\textsf{DP}}(\mathcal{I}_{M_{c}})| + |\mathcal{J}^{\textsf{DP}}(\mathcal{I}_{M_{p}})| + |\mathcal{J}^{\textsf{DP}}(\mathcal{I}_{L})|
     \end{align}
     Conditional on the event $\mathcal{E} \cap \Tilde{\mathcal{H}} \cap \mathcal{C} \cap \mathcal{P}$, we have that 
     \begin{align*}
         &\mathbb{E}\left[\Lambda(B,T; \textsf{DP}) | \mathcal{E} \cap \Tilde{\mathcal{H}} \cap \mathcal{C} \cap \mathcal{P}\right] \\ 
         &\stackrel{(a)}= \sum_{k \in \mathcal{J}^{\textsf{hs}}(\mathcal{I}_{H}) \backslash \mathcal{J}^{\textsf{DP}}(\mathcal{I}_{H})} \theta_{k} - \left( \sum_{k \in \mathcal{J}^{\textsf{DP}}(\mathcal{I}_{L})} \theta_{k} + \sum_{k \in \mathcal{J}^{\textsf{DP}}(\mathcal{I}_{M_{p}})} \theta_{k} + \sum_{k \in \mathcal{J}^{\textsf{DP}}(\mathcal{I}_{M_{c}})} \theta_{k}\right), \\
         &\stackrel{(b)}\geq |\mathcal{J}^{\textsf{hs}}(\mathcal{I}_{H}) \backslash \mathcal{J}^{\textsf{DP}}(\mathcal{I}_{H})| \left(u + c_{0} \Delta_{\beta} \right) - \left( \sum_{k \in \mathcal{J}^{\textsf{DP}}(\mathcal{I}_{L})} \theta_{k} + \sum_{k \in \mathcal{J}^{\textsf{DP}}(\mathcal{I}_{M_{p}})} \theta_{k} + \sum_{k \in \mathcal{J}^{\textsf{DP}}(\mathcal{I}_{M_{c}})} \theta_{k}\right), \\
         &\stackrel{(c)}\geq |\mathcal{J}^{\textsf{hs}}(\mathcal{I}_{H}) \backslash \mathcal{J}^{\textsf{DP}}(\mathcal{I}_{H})| \left(u + c_{0} \Delta_{\beta} \right) - \left(|\mathcal{J}^{\textsf{DP}}(\mathcal{I}_{L})| \ell + |\mathcal{J}^{\textsf{DP}}(\mathcal{I}_{M_{p}})| (u + c_{0} \Delta_{\beta}) + |\mathcal{J}^{\textsf{DP}}(\mathcal{I}_{M_{c}})| (u + \Delta_{\beta}) \right) , \\
         &\stackrel{(d)}\geq |\mathcal{J}^{\textsf{DP}}(\mathcal{I}_{L})|\left[\left(u + c_{0} \Delta_{\beta} \right) - \ell \right] + |\mathcal{J}^{\textsf{DP}}(\mathcal{I}_{M_{p}})| \left[\left(u + c_{0} \Delta_{\beta} \right) - \left( u + c_{0} \Delta_{\beta} \right) \right] \\ & \quad \quad + |\mathcal{J}^{\textsf{DP}}(\mathcal{I}_{M_{c}})| \left[\left(u + c_{0} \Delta_{\beta} \right) - \left(u + \Delta_{\beta} \right) \right] \\
         &\stackrel{(e)}= |\mathcal{J}^{\textsf{DP}}(\mathcal{I}_{L})|\left[\left(u + c_{0} \Delta_{\beta} \right) - \ell \right] + |\mathcal{J}^{\textsf{DP}}(\mathcal{I}_{M_{c}})| \alpha_{0} \Delta_{\beta} ,\\
         &\stackrel{(f)}\geq \frac{\alpha_{0} \ell}{32} T^{\frac{1}{2} - \frac{1}{2(1 + \beta)}},
     \end{align*}
     where (a) follows from \eqref{eq:dp-cost-decomposition} and \eqref{eq:off-cost-simplified-1}, (b) follows from the fact that $\sum_{k \in \mathcal{S}} a_{k} \geq |\mathcal{S}| \min_{k \in \mathcal{S}}\{a_{k}\}$ and by construction, for all the arrivals in the set $\mathcal{I}_{H}, \theta_{k} \geq u + c_{0} \Delta_{\beta}$, (c) follows similar to (b), (d) follows from \eqref{eq:set-size-inequality}, (e) follows from the definition of $\alpha_{0}$ in \eqref{eq:notation}, (f) follows from the fact that $|\mathcal{J}^{\textsf{DP}}(\mathcal{I}_{M_{c}})| \geq \frac{1}{32}\sqrt{T}$ due to the event $\mathcal{E}$ and from the fact that $|\mathcal{J}^{\textsf{DP}}(\mathcal{I}_{L})| \geq 0$. \QEDB

\subsection{Proof of Lemma \ref{lem:conditional-regret-lb-2}}
\label{subsec:proof-lemma-conditional-regret-lb-2}

    Recall the definitions of set of indices $\mathcal{J}^{\textsf{hs}}$ and $\mathcal{J}^{\textsf{DP}}$ from \eqref{def:random-indices-set} and decomposition of the sum of values chosen under hindsight optimal and the DP policy as given in \eqref{eq:offline-cost-decomposition} and \eqref{eq:dp-cost-decomposition}. Recall from the discussion of the hindsight optimal that under the event $ \mathcal{E}^{c} \cap \mathcal{H} \cap \mathcal{C} \cap \mathcal{P}$, the hindsight optimal will accept all the candidates with abilities in the set $\mathcal{I}_{H}, \mathcal{I}_{M_{p}}, \mathcal{I}_{M_{c}}$ and possibly some candidates in the set $\mathcal{I}_{L}$.
   Conditional on the event $\mathcal{E}^{c} \cap \mathcal{H} \cap \mathcal{C} \cap \mathcal{P}$, we have that,
    \begin{align*}
        B &\stackrel{(a)}= |\mathcal{J}^{\textsf{hs}}(\mathcal{I}_{H})| + |\mathcal{J}^{\textsf{hs}}(\mathcal{I}_{M_{p}})| + |\mathcal{J}^{\textsf{hs}}(\mathcal{I}_{M_{c}})| + |\mathcal{J}^{\textsf{hs}}(\mathcal{I}_{L})| \\
        &\stackrel{(b)}= |\mathcal{J}^{\textsf{hs}}(\mathcal{I}_{H}) \backslash \mathcal{J}^{\textsf{DP}}(\mathcal{I}_{H})| + |\mathcal{J}^{\textsf{DP}}(\mathcal{I}_{H})| + |\mathcal{J}^{\textsf{hs}}(\mathcal{I}_{M_{p}}) \backslash \mathcal{J}^{\textsf{DP}}(\mathcal{I}_{M_{p}})| + |\mathcal{J}^{\textsf{DP}}(\mathcal{I}_{M_{p}})| \\ & \quad \quad + |\mathcal{J}^{\textsf{hs}}(\mathcal{I}_{M_{c}}) \backslash \mathcal{J}^{\textsf{DP}}(\mathcal{I}_{M_{c}})| + |\mathcal{J}^{\textsf{DP}}(\mathcal{I}_{M_{c}})| + |\mathcal{J}^{\textsf{hs}}(\mathcal{I}_{L})|, \\
        &\stackrel{(c)}\geq |\mathcal{J}^{\textsf{DP}}(\mathcal{I}_{H})| + |\mathcal{J}^{\textsf{DP}}(\mathcal{I}_{M_{p}})| + |\mathcal{J}^{\textsf{DP}}(\mathcal{I}_{M_{c}})| + |\mathcal{J}^{\textsf{DP}}(\mathcal{I}_{L})|, \\
        &\stackrel{(d)}= |\mathcal{J}^{\textsf{DP}}(\mathcal{I}_{H})| + |\mathcal{J}^{\textsf{DP}}(\mathcal{I}_{M_{p}})| + |\mathcal{J}^{\textsf{DP}}(\mathcal{I}_{M_{c}})| + |\mathcal{J}^{\textsf{DP}}(\mathcal{I}_{L}) \backslash \mathcal{J}^{\textsf{hs}}(\mathcal{I}_{L})| + |\mathcal{J}^{\textsf{hs}}(\mathcal{I}_{L})|,
    \end{align*}
    where (a) follows from the fact the hindsight optimal will accept exactly $B$ candidates, (b) follows from the fact that for countable set $A,B$ such that $B \subseteq A$, we have that $|A| = |A \backslash B| + |B|$, (c) follows from the fact that any online policy will accept at most $B$ candidates, (d) follows for the same reason as (b). This implies the following inequality,
    \begin{align}
        \label{eq:set-size-inequality-2}
        |\mathcal{J}^{\textsf{hs}}(\mathcal{I}_{H}) \backslash \mathcal{J}^{\textsf{DP}}(\mathcal{I}_{H})| + |\mathcal{J}^{\textsf{hs}}(\mathcal{I}_{M_{p}}) \backslash \mathcal{J}^{\textsf{DP}}(\mathcal{I}_{M_{p}})| + |\mathcal{J}^{\textsf{hs}}(\mathcal{I}_{M_{c}}) \backslash \mathcal{J}^{\textsf{DP}}(\mathcal{I}_{M_{c}})| \geq |\mathcal{J}^{\textsf{DP}}(\mathcal{I}_{L}) \backslash \mathcal{J}^{\textsf{hs}}(\mathcal{I}_{L})|
    \end{align}
    Conditional on the event $\mathcal{E}^{c} \cap \mathcal{H} \cap \mathcal{C} \cap \mathcal{P}$, we have that
    \begin{align*}
        &\mathbb{E}\left[\Lambda(B,T; \textsf{DP}) | \mathcal{E}^{c} \cap \mathcal{H} \cap \mathcal{C} \cap \mathcal{P} \right] \\
        &\stackrel{(a)}= \left(\sum_{\mathcal{A} \in \{\mathcal{I}_{H}, \mathcal{I}_{M_{p}}, \mathcal{I}_{M_{c}}\}} \sum_{k \in \mathcal{J}^{\textsf{hs}}(\mathcal{A}) \backslash \mathcal{J}^{\textsf{DP}}(\mathcal{A})} \theta_{k} + \sum_{\mathcal{A} \in \{\mathcal{I}_{H}, \mathcal{I}_{M_{p}}, \mathcal{I}_{M_{c}}\}} \sum_{k \in \mathcal{J}^{\textsf{DP}}(\mathcal{A})} \theta_{k} + \sum_{k \in \mathcal{J}^{\textsf{hs}}(\mathcal{I}_{L})} \theta_{k} \right)  \\
        &\quad - \left(\sum_{\mathcal{A} \in \{\mathcal{I}_{H}, \mathcal{I}_{M_{p}}, \mathcal{I}_{M_{c}} \}} \sum_{k \in \mathcal{J}^{\textsf{DP}}(\mathcal{A})}  \theta_{k} + \sum_{k \in \mathcal{J}^{\textsf{DP}}(\mathcal{I}_{L}) \backslash \mathcal{J}^{\textsf{hs}}(\mathcal{I}_{L})} \theta_{k} + \sum_{k \in \mathcal{J}^{\textsf{hs}}(\mathcal{I}_{L})} \theta_{k} \right), \\
        &\stackrel{(b)}= \sum_{\mathcal{A} \in \{\mathcal{I}_{H}, \mathcal{I}_{M_{p}}, \mathcal{I}_{M_{c}}\}} \sum_{k \in \mathcal{J}^{\textsf{hs}}(\mathcal{A}) \backslash \mathcal{J}^{\textsf{DP}}(\mathcal{A})} \theta_{k} - \sum_{k \in \mathcal{J}^{\textsf{DP}}(\mathcal{I}_{L}) \backslash \mathcal{J}^{\textsf{hs}}(\mathcal{I}_{L})} \theta_{k}, \\
        &\stackrel{(c)}\geq |\mathcal{J}^{\textsf{hs}}(\mathcal{I}_{H}) \backslash \mathcal{J}^{\textsf{DP}}(\mathcal{I}_{H})| \left(u + c_{0}\Delta_{\beta} \right) + |\mathcal{J}^{\textsf{hs}}(\mathcal{I}_{M_{p}}) \backslash \mathcal{J}^{\textsf{DP}}(\mathcal{I}_{M_{p}})| \left(\ell - c_{0} \Delta_{\beta} \right)  \\ & \quad + |\mathcal{J}^{\textsf{hs}}(\mathcal{I}_{M_{c}}) \backslash \mathcal{J}^{\textsf{DP}}(\mathcal{I}_{M_{c}})| \left(\ell - \Delta_{\beta} \right) - |\mathcal{J}^{\textsf{DP}}(\mathcal{I}_{L}) \backslash \mathcal{J}^{\textsf{hs}}(\mathcal{I}_{L})| \left( \ell - c_{0} \Delta_{\beta} \right) \\
        &\stackrel{d}\geq |\mathcal{J}^{\textsf{hs}}(\mathcal{I}_{H}) \backslash \mathcal{J}^{\textsf{DP}}(\mathcal{I}_{H})|\left[\left(u + c_{0} \Delta_{\beta} \right) - \left(\ell - c_{0}\Delta_{\beta} \right)\right] + |\mathcal{J}^{\textsf{hs}}(\mathcal{I}_{M_{c}}) \backslash \mathcal{J}^{\textsf{DP}}(\mathcal{I}_{M_{c}})| \left[\left(\ell - \Delta_{\beta} \right) - \left(\ell - c_{0} \Delta_{\beta} \right) \right] \\
        & \quad + |\mathcal{J}^{\textsf{hs}}(\mathcal{I}_{M_{p}}) \backslash \mathcal{J}^{\textsf{DP}}(\mathcal{I}_{M_{p}})| \left[ \left( \ell - c_{0}\Delta_{\beta} \right) - \left( \ell - c_{0}\Delta_{\beta} \right) \right] \\
        &\stackrel{(e)}= |\mathcal{J}^{\textsf{hs}}(\mathcal{I}_{H}) \backslash \mathcal{J}^{\textsf{DP}}(\mathcal{I}_{H})|\left[\left(u + c_{0} \Delta_{\beta} \right) - \left(\ell - c_{0}\Delta_{\beta} \right)\right] + |\mathcal{J}^{\textsf{hs}}(\mathcal{I}_{M_{c}}) \backslash \mathcal{J}^{\textsf{DP}}(\mathcal{I}_{M_{c}})| \alpha_{0}\Delta_{\beta}, \\
        &\stackrel{(f)}\geq \frac{\alpha_{0} \ell}{8} T^{\frac{1}{2} - \frac{1}{2(1 + \beta)}},
    \end{align*}
    where (a) follows from \eqref{eq:offline-cost-decomposition} and \eqref{eq:dp-cost-decomposition}, (b) follows trivially, (c) follows from the fact $\theta_{k} \mathbbm{1}\{\theta_{k} \in \mathcal{I}_{H}\} \geq u + c_{0}\Delta_{\beta}, \theta_{k} \mathbbm{1}\{\theta_{k} \in \mathcal{I}_{M_{p}}\} \geq \ell - c_{0}\Delta_{\beta}, \theta_{k} \mathbbm{1}\{\theta_{k} \in \mathcal{I}_{M_{c}}\} \geq \ell - \Delta_{\beta}$ and $\theta_{k} \mathbbm{1}\{\theta_{k} \in \mathcal{I}_{L}\} \leq \ell - c_{0}\Delta_{\beta}$, (d) follows from \eqref{eq:set-size-inequality-2}, (e) follows from the definition of $\alpha_{0} = c_{0} - 1$, (f) follows from the fact that $|\mathcal{J}^{\textsf{hs}}(\mathcal{I}_{M_{c}}) \backslash \mathcal{J}^{\textsf{DP}}(\mathcal{I}_{M_{c}})| \geq \frac{1}{8}\sqrt{T}$ which is due to fact that under the event $\mathcal{H} \cap \mathcal{C} \cap \mathcal{P}$, the hindsight optimal will accept all the arrivals in the set $\mathcal{I}_{M_{c}}$ however under the event $\mathcal{E}^{c}$ will accept at most $\frac{1}{8}\sqrt{T}$ arrivals in the first $B$ time steps and this will result in incorrectly rejecting at least $\frac{1}{4}\sqrt{T} - \frac{1}{8}\sqrt{T} = \frac{1}{8}\sqrt{T}$ arrivals in the set $\mathcal{I}_{M_{c}}$. \QEDB

\subsection{Proof of Lemma \ref{lem:high-type-concentration}}
    \label{subsec:high-type-concentration}
    Recall the definition of the events $\mathcal{H}_{1}, \mathcal{H}_{2}$ and $\Tilde{\mathcal{H}}_{2}$ as defined in \eqref{event:high-type-first-half-low}, \eqref{event:high-type-second-half-low} and \eqref{event:high-type-second-half-high} respectively.
    We have that $N(\mathcal{I}_H, 1, B) \sim \text{Bin}\left(B, p_{\mathcal{I}_H}\right)$, where $p_{\mathcal{I}_H} = \frac{1}{2} - \frac{129}{256}T^{-\frac{1}{2}}$. Therefore we have that $\mu_1^B(\mathcal{I}_H) = \mathbb{E}\left[N(\mathcal{I}_H, 1, B) \right] \approx \frac{T}{4} - \frac{129}{512}\sqrt{T}$ since $B \approx T/2$. Therefore we can write the event $\mathcal{H}_1$ as 
    \begin{align*}
        \mathcal{H}_1 = \bigg\{\mu_1^B(\mathcal{I}_{H}) - \frac{127}{512}\sqrt{T} \leq N(\mathcal{I}_H, 1, B) \leq \mu_{1}^{B}(\mathcal{I}_H) + \frac{129}{512}\sqrt{T}\bigg\}
    \end{align*}
    Therefore we have that 
    \begin{align*}
        \mathbb{P}\left(\mathcal{H}_1 \right) &\stackrel{(a)}= \mathbb{P}\left(\mu_1^B(\mathcal{I}_{H}) - \frac{127}{512}\sqrt{T} \leq N(\mathcal{I}_H, 1, B) \leq \mu_{1}^{B}(\mathcal{I}_H) + \frac{129}{512}\sqrt{T} \right) \\
        &\stackrel{(b)}= \mathbb{P}\left( - \frac{127}{512}\sqrt{T} \leq N(\mathcal{I}_H, 1, B) - \mu_1^B(\mathcal{I}_H) \leq \frac{129}{512}\sqrt{T} \right) \\
        &\stackrel{(c)}= \mathbb{P}\left( - \frac{127}{256}\frac{1}{\sqrt{p_{\mathcal{I}_H}(1 - p_{\mathcal{I}_H})}} \leq \frac{N(\mathcal{I}_H, 1, B) - \mu_1^B(\mathcal{I}_H)}{\sqrt{B p_{\mathcal{I}_H}(1 - p_{\mathcal{I}_H})}} \leq \frac{129}{256}\frac{1}{\sqrt{p_{\mathcal{I}_H}(1 - p_{\mathcal{I}_H})}}\right) \\
        &\stackrel{(d)}\geq \mathbb{P}\left( 0 \leq \frac{N(\mathcal{I}_H, 1, B) - \mu_1^B(\mathcal{I}_H)}{\sqrt{B p_{\mathcal{I}_H}(1 - p_{\mathcal{I}_H})}} \leq \frac{258}{256} \right) \\
        &\stackrel{(e)}= \mathbb{P}\left( \frac{N(\mathcal{I}_H, 1, B) - \mu_1^B(\mathcal{I}_H)}{\sqrt{B p_{\mathcal{I}_H}(1 - p_{\mathcal{I}_H})}} \leq \frac{258}{256}\right) - \mathbb{P}\left( \frac{N(\mathcal{I}_H, 1, B) - \mu_1^B(\mathcal{I}_H)}{\sqrt{B p_{\mathcal{I}_H}(1 - p_{\mathcal{I}_H})}} \leq 0\right) \\
        &\stackrel{(f)}\geq \Phi\left(\frac{258}{256} \right) - \Phi\left( 0\right) - \frac{c}{\sqrt{T}} \\
        &\stackrel{(g)}\geq 0.34 - \frac{c}{\sqrt{T}}
    \end{align*}
    where (a) follows from definition of event $\mathcal{H}_1$, (b,c) follows trivially, (d) follows the fact that $\bigg\{0 \leq \frac{N(\mathcal{I}_H, 1, B) - \mu_1^B(\mathcal{I}_H)}{\sqrt{B p_{\mathcal{I}_H}(1 - p_{\mathcal{I}_H})}} \leq \frac{258}{256}\bigg\} \subseteq \bigg\{ - \frac{127}{256}\frac{1}{\sqrt{p_{\mathcal{I}_H}(1 - p_{\mathcal{I}_H})}} \leq \frac{N(\mathcal{I}_H, 1, B) - \mu_1^B(\mathcal{I}_H)}{\sqrt{B p_{\mathcal{I}_H}(1 - p_{\mathcal{I}_H})}} \leq \frac{129}{256}\frac{1}{\sqrt{p_{\mathcal{I}_H}(1 - p_{\mathcal{I}_H})}} \bigg\}$ since $p_{\mathcal{I}_H}(1 - p_{\mathcal{I}_H}) \leq 1/4$, (e) follows trivially, (f) follows from Berry Esseen Theorem and (g) follows trivially.
    Now there exists a $T_0 < \infty$ such that for all $T \geq T_0$, we have that $\mathbb{P}(\mathcal{H}_1) \geq 0.003$. An analogous proof follows for $\mathcal{H}_2$ and $\tilde{\mathcal{H}}_2$ as well, we omit it to avoid repetition. \QEDB

\subsection{Proof of Lemma \ref{lem:mid-type-concentration}}
\label{subsec:niddle-type-concentration}
Recall the definition of events $\mathcal{C}_1, \mathcal{C}_2, \mathcal{P}_1$ and $\mathcal{P}_2$ as defined in \eqref{event:middle-center-type-first-half}, \eqref{event:middle-center-type-second-half}, \eqref{event:middle-periphery-type-first-half} and \eqref{event:middle-periphery-type-second-half} respectively. We have that $N(\mathcal{I}_{M_c}, 1, B) \sim \text{Bin}(B, p_{\mathcal{I}_{M_c}})$ where $p_{\mathcal{I}_{M_c}} = T^{-\frac{1}{2}}$. Therefore we have that $\mu_1^B(\mathcal{I}_{M_c}) = \mathbb{E}\left[ N(\mathcal{I}_{M_c}, 1, B)\right] \approx \sqrt{T}/2$ since $B \approx T/2$. Now the event $\mathcal{C}_1$ can be written as
\begin{align*}
    \mathcal{C}_1 = \bigg\{\mu_1^B(\mathcal{I}_{M_c}) / 2 \leq N(\mathcal{I}_{M_c}, 1, B) \leq 2\mu_1^B(\mathcal{I}_{M_c})\bigg\}
\end{align*}
Therefore we have that 
\begin{align*}
    \mathbb{P}\left(\mathcal{C}_1^c \right) &\stackrel{(a)}= \mathbb{P}\left(\{N(\mathcal{I}_{M_c}, 1, B) \geq 2 \mu_1^B(\mathcal{I}_{M_c})\} \cup \{N(\mathcal{I}_{M_c}, 1, B) \leq \mu_1^B(\mathcal{I}_{M_c}) / 2\} \right) \\
    &\stackrel{(b)}\leq \mathbb{P}\left(\{N(\mathcal{I}_{M_c}, 1, B) \geq 2 \mu_1^B(\mathcal{I}_{M_c})\} \right) + \mathbb{P}\left(\{N(\mathcal{I}_{M_c}, 1, B) \leq \mu_1^B(\mathcal{I}_{M_c}) / 2\} \right) \\
    &\stackrel{(c)}\leq cT^{-\frac{1}{4}}
\end{align*}
where (a) follows from definition of $\mathcal{C}_1$, (b) follows from union bound, (c) follows from Berry Esseen theorem as applied before. From this it follows that there exists $T_0 < \infty$ such that $\mathbb{P}(\mathcal{C}_{1}^c) \leq 0.001$ for all $T \geq T_0$. An analogous proof follows for $\mathcal{C}_2, \mathcal{P}_1$ and $\mathcal{P}_2$, we omit it to avoid repetition. \QEDB

\section{Details and Analysis of {\sf CwG} Policy}
\label{onapp:analysis-cwg-policy}

In this section we will provide some more details about the {\sf CwG} algorithm (Algorithm \ref{alg:CwG}) and also provide the proof of Theorem \ref{thm:general-cwg-upper-bound}. In Section \ref{subsec:phase-structure-cwg}, we provide a discussion about the phase structure of Algorithm \ref{alg:CwG}. In Section \ref{onapp:hindsight-to-go-threshold}, we define the concept of {\it hindsight-to-go} ({\sf HTG}) which will aid our analysis. In Section \ref{onapp:proof-outline-cwg-analysis}, we provide a proof outline for Theorem \ref{thm:general-cwg-upper-bound}. In Section \ref{sec:helper_lemmas}, we will provide some helper lemmas to formalize our analysis with their proofs deferred to Section \ref{onapp:helper-lemma-cwg-proofs}. In Section \ref{onapp:formal-proof-cwg-upper-bound}, we provide the formal proof of Theorem \ref{thm:general-cwg-upper-bound}.

\new{
\subsection{Phase Structure of Algorithm \ref{alg:CwG}}
\label{subsec:phase-structure-cwg}
The phase structure of the {\sf CwG} policy has been devised to simplify the analysis of the {\sf CwG} policy. The key idea of the {\sf CwG} policy is that if the {\sf CE} threshold $p_{t}^{\textsf{CE}}$ at time $t$ is within a ball of radius  $\Delta_{t} \triangleq \sqrt{2\log \tau / \tau}$ (where $\tau = T - t + 1$ is the number of remaining time steps) of a gap quantile $q_{i}^{\star}$, then the CwG threshold is set to the gap quantile $q_{i}^{\star}$ itself. As $t$ increases, so does the size of the radius and hence eventually there will be more than one gap quantiles in this ball. If there are more than one gap quantiles in $\Delta_{t}$, we need a tie-breaking rule to decide which gap quantile the {\sf CwG} threshold is assigned to. This tie-breaking rule further complicates an already involved analysis and hence to simplify the technical analysis, we define the {\sf CwG} algorithm by dividing it into two phases. 
    
    In the first phase, it suffices to ensure that there will always be at most one gap quantile in the $\Delta_{t}$-neighbourhood of $p_{t}^{\textsf{CE}}$ for any value of $p_{t}^{\textsf{CE}}$ and there is no need for a tie-breaking rule. One way to ensure this, is to find $t^{\star}$ such that for all $t \leq t^{\star}$, we have that $\sqrt{2 \log \tau / \tau} \leq \varepsilon_{0} / 2$. Note that irrespective of the value of $p_{t}^{\textsf{CE}}$, there is at most one gap quantile in the $\Delta_{t}$ neighborhood of $p_{t}^{\textsf{CE}}$. Further, note that for $t \leq T - 2$, $\sqrt{2 \log \tau / \tau}$ is increasing in $t$ and hence it suffices to verify that $\sqrt{2 \log \tau^\star / \tau^{\star}} \leq \varepsilon_{0} / 2$ for $\tau^{\star} = \ceil{64 \log (1 / \varepsilon_{0}) / \varepsilon_{0}^{2}}$. Given that we are guaranteed to have at most one gap quantile in $\Delta_{t}$-neighbourhood of $p_{t}^{\textsf{CE}}$, our analysis is great simplified.

    The second phase is of length $\ceil{64 \log (1/\varepsilon_{0}) / \varepsilon_{0}^{2}}$ and we use a static allocation rule in the second phase. The contribution to regret because of the static policy is at most $C\sqrt{\log (1/\varepsilon_{0})} / \varepsilon_{0}$ for some universal constant $C < \infty$.
}

\subsection{Hindsight To Go ({\sf HTG}) and {\sf HTG} Threshold}
\label{onapp:hindsight-to-go-threshold}
Let $\bm{q}^{\theta}_{\geq t}(n)$ denote the $n$-th largest value quantile in $\bm{q}^{\theta}_{\geq t}$ for an integer $n \in \mathbb{N}$. Define the following quantile values $q_{t}^l \triangleq \bm{q}^{\theta}_{\geq t}(B_t + 1)$ and $q_t^u \triangleq \bm{q}^{\theta}_{\geq t}(B_t)$ and denote their corresponding values by $l_t = F^{-1}(q_{t}^l), u_t = F^{-1}(q_t^u)$, where $B_t$ is the remaining budget at time $t$. Note that since the principle of compensated coupling is to persuade the hindsight policy to take the same action as the online policy using sufficient compensations, the hindsight policy at time $t$ may look different from the hindsight policy initially and being adapted to the budget which evolves according to the online policy. To distinguish between the two, at any time $t$, we will instead refer to the hindsight policy as the \emph{Hindsight To Go} ({\sf HTG}) policy, which due to coupling follows the same actions as the online policy up till time $t - 1$ and then from time $t$ onwards takes the optimal hindsight decision with arrivals in $\vect{\omega}_{\geq t}$ given the remaining budget $B_t$. Given the \textsf{CwG} quantile threshold $p_t^{\sf CwG}$, we define $p_{t}^{\sf HTG} \triangleq \argmax_{x \in [q_t^l, q_t^u]} |p_{t}^{\sf CwG} - x|$ when $B_t > 0$, otherwise $p_t^{\textsf{HTG}} = 1$. The reason to adopt this particular $p_{t}^{\textsf{CwG}}$ dependent definition of $p_{t}^{\textsf{HTG}}$ is that the compensation needed at time $t$ will now be bounded above by the separation between the {\sf CwG} threshold and the {\sf HTG} threshold in value space.

\subsection{Proof Outline}
\label{onapp:proof-outline-cwg-analysis}
We first provide a proof outline.
Recall $\tilde{T} = T - \ceil{64 \log (1 / \varepsilon_{0})} / \varepsilon_{0}^{2}$ in Algorithm \ref{alg:CwG} (the {\sf CwG} policy). The algorithm operates in two phases, the first phase includes time steps $t$ such that $1 \leq t \leq \tilde{T}$ while the second phase consists of the remaining time steps $t$ such that $\Tilde{T} + 1 \leq t \leq T$. 
\paragraph{\underline{Analysis of First Phase}.} The analysis of the first phase makes use of the regret decomposition given in Lemma \ref{lem:cwg-regret-decomposition}. To bound the expected compensation term $\mathbb{E}_{B_{t}^{\pi}}\left[\partial \mathcal{R}_{t}(B_{t}^{\pi}, a_{t}^{\pi}) \right]$ in Lemma \ref{lem:cwg-regret-decomposition} for $\pi = {\sf CwG}$, we will analyse two thresholds: the \textsf{CwG} quantile threshold denoted as $p_{t}^{\textsf{CwG}}$ and Hindsight To Go (\textsf{HTG}) quantile threshold $p_{t}^{\textsf{HTG}}$. Note that given a tail sequence $\vect{\theta}_{\geq t}$ and the remaining budget $B_{t}$, the Hindsight To Go threshold is set such that on the sample path $\vect{\theta}_{\geq t}$, the top $B_{t}$ candidates are chosen. We bound the expected compensation at time $t$ for $t \in [1, \tilde{T}]$ and we do so by dividing the analysis into two events: $(a) E_{t} = \{1 - B_{t} / \tau > 4 \sqrt{\log \tau / \tau} \}$ and $(b) E_{t}^c = \{1 - B_{t} / \tau \leq 4 \sqrt{\log \tau / \tau}\}$ where $\tau = T - t + 1$. At any time either of the two events arises and we bound the expected compensation conditional on each of the two events. The analysis for both the events utilizes the same recipe. We show that with high probability the difference between \textsf{CwG} quantile threshold $p_{t}^{\textsf{CwG}}$ and the \textsf{HTG} quantile threshold $p_{t}^{\textsf{HTG}}$ is bounded above by $C \sqrt{\log \tau / \tau}$ (Lemma \ref{lem:crt-HTG-close}). As a result of this, we establish that with high probability the two thresholds $p_{t}^{\textsf{CwG}}$ and $p_{t}^{\textsf{HTG}}$ belong to the same cluster (Lemma \ref{lem:opposite-side-crt-HTG}). Now compensation is need at time $t$ only if there is a candidate ability arrival $\theta_t$ such that its quantile $F(\theta_t)$ lies between the two thresholds $p_{t}^{\textsf{CwG}}$ and $p_{t}^{\textsf{HTG}}$ and the amount of compensation is bounded by $|F^{-1}(p_{t}^{\textsf{CwG}}) - F^{-1}(p_{t}^{\textsf{HTG}})|$ (Lemma \ref{lem:compensation-upper-bound}). Using Lemmas \ref{lem:crt-HTG-close}, \ref{lem:opposite-side-crt-HTG}, \ref{lem:compensation-upper-bound} and definition of the $(\beta, \varepsilon_{0}, \delta)$-clustered distribution, we show that the expected compensation at time $t$ is bounded as follow.
\begin{lemma}
    \label{lem:expected-compensation-bound}
    There is a universal constant $C < \infty$ such that the following occurs. For any $\beta \in [0, \infty), \varepsilon_{0} \in (0,1]$ and $\delta \in (0,1]$, suppose the candidate-ability distribution $F$ with associated gaps is $(\beta, \varepsilon_{0}, \delta)$-clustered. Then for $t \in \{1,2,\dots, \tilde{T}\}$, for the \textsf{CwG} policy we have that the expected compensation at time $t$ is bounded above as
    \begin{align*}
        \sup_{B_{t} \geq 0} \partial \mathcal{R}_{t}\left( B_{t}, a_{t}^{\textsf{CwG}}\right) \leq C\left(\left( \log \tau / \tau \right)^{\frac{1}{2} + \frac{1}{2(1 + \beta)}} + \delta \sqrt{\log \tau / \tau} \right),
    \end{align*}
    where $\tau = T - t + 1$. Note that the above implies that 
    \begin{align*}
\mathbb{E}_{B_{t}}\left[\partial \mathcal{R}_{t}\left(B_{t}, a_{t}^{\textsf{CwG}} \right) \right] \leq C\left(\left( \log \tau / \tau \right)^{\frac{1}{2} + \frac{1}{2(1 + \beta)}} + \delta \sqrt{\log \tau / \tau} \right).
    \end{align*}
\end{lemma}
Using Lemma \ref{lem:cwg-regret-decomposition} and Lemma \ref{lem:expected-compensation-bound}, the cummulative regret accrued up till time $\tilde{T}$ is upper bounded by
    $C\left((\log T)^{\frac{1}{2} + \frac{1}{2(1 + \beta)}} T^{\frac{1}{2} - \frac{1}{2(1 + \beta)}} \mathbbm{1}\{\beta > 0\} + \log^2 T \mathbbm{1}\{\beta = 0\} + \delta \sqrt{T \log T} \right).$

\paragraph{\underline{Analysis of Second Phase}.}
Recall that the {\sf CwG} policy (Algorithm 2) in the last $64 \log (1 / \varepsilon_{0}) / \varepsilon_{0}^{2}$ time steps, makes use of the static allocation policy where we solve for the CE quantile threshold $p_{\tilde{T}}^{\textsf{CE}}$ and thereafter use the time invariant quantile threshold $p_{\tilde{T}}^{\textsf{CE}}$. Using a well known fact in the network revenue management literature, we know that the regret accrued under a static allocation policy is upper bounded as $C \sqrt{\text{horizon length}}$ for some universal constant $C < \infty$. Since the \textsf{CwG} policy (Algorithm \ref{alg:CwG}) employs the static allocation policy for the last $\ceil{64 \log (1 / \varepsilon_{0}) / \varepsilon_{0}^{2}}$, the regret accrued over the last $T - \tilde{T}$ time steps is upper bounded as $C \sqrt{\log (1 / \varepsilon_{0})} / \varepsilon_{0}$. Adding up the regret over the two phases results in the regret scaling in Theorem \ref{thm:general-cwg-upper-bound}.

\subsection{Preliminaries and Helper Lemmas}
\label{sec:helper_lemmas}
We introduce some helper lemmas which we will use to prove the regret bound. We defer the proof of these lemmas to Appendix \ref{onapp:helper-lemma-cwg-proofs}. Let $t$ denote the current time step and $\tau = T - t + 1$ denote the remaining number of times steps. Assume that $T \geq \ceil{64\log(1 / \varepsilon_{0})/ \varepsilon_{0}^{2}}$ and define $\tilde{T} \triangleq T - \floor{ 64 \log(1 / \varepsilon_{0}) / \varepsilon_{0}^{2}}$.
Define the following events for $t \leq \tilde{T}$:
\begin{align}
\label{eq:event_ce_HTG_close}
    \mathcal{A}_{1,t} &= \{|p_{t}^{\sf CE} - p_{t}^{\sf HTG}| \leq \sqrt{2 \log \tau / \tau}\}, \\
    \label{eq:event_cwg_HTG_close}
    \mathcal{A}_{2,t} &= \{|p_{t}^{\sf CwG} - p_{t}^{\sf HTG}| \leq 3 \sqrt{\log \tau / \tau}\}, \\
    \label{eq:event_cwg_HTG_same_quantile}
    \mathcal{A}_{3,t} &= \cup_{i = 1}^{n + 1}\{p_{t}^{\sf CwG} \in \bar{Q}_{i}, p_{t}^{\sf HTG} \in \bar{Q}_{i}\},
\end{align}
where $\bar{Q}_{i} = [q_{i - 1}^{\star}, q_{i}^{\star}]$ and $n$ denotes the number of gaps. The interpretation of $A_{3,t}$ is that the {\sf CwG} policy threshold and the {\sf HTG} policy threshold  belong (weakly) to the same mass cluster. 
The following lemmas show that these three events are very likely to occur for $t \leq \tilde{T}$:
\begin{lemma}
\label{lem:rt-HTG-close}
Consider the event ${\cal A}_{1,t}$ defined in \eqref{eq:event_ce_HTG_close}. We have that $\mathbb{P}( {\cal A}_{1,t}^c) \leq 2/ \tau^{4}$.
\end{lemma}

\begin{lemma}
\label{lem:crt-HTG-close}
Consider the event ${\cal A}_{2,t}$ defined in \eqref{eq:event_cwg_HTG_close}. We have that $\mathbb{P}({\cal A}_{2,t}^c) \leq {2} / {\tau^{4}}$.
\end{lemma}

\begin{lemma}
\label{lem:opposite-side-crt-HTG}
Consider the event ${\cal A}_{3,t}$ defined in \eqref{eq:event_cwg_HTG_same_quantile}. We have that $\mathbb{P}({\cal A}_{3,t}^{c}) \leq 2n(n + 1) / \tau^{4}$, where $n$ is the number of gaps.
\end{lemma}

Let  $q^{\theta}_t = F(\theta_{t})$ be the quantile of the candidate ability $\theta_{t}$ at time $t$. 
If $p_{t}^{\textsf{CwG}} < q_{t}^{l}$ then we have that $p_{t}^{\textsf{HTG}} = q_{t}^{u}$ and compensation is needed only if $q_{t}^{\theta} \in [p_{t}^{\textsf{CwG}}, p_{t}^{\textsf{HTG}}]$. If $p_{t}^{\textsf{CwG}} > q_{t}^{u}$ then we have that $p_{t}^{\textsf{HTG}} = q_{t}^{l}$ and compensation is needed only if $q^{\theta}_t \in [p_{t}^{\textsf{HTG}}, p_{t}^{\textsf{CwG}}]$. If $p_{t}^{\textsf{CwG}} \in (q_{t}^{l}, q_{t}^{u})$, then no compensation is required.
\begin{lemma}
\label{lem:compensation-upper-bound}
Let $q_{t}^{\theta} = F(\theta_{t})$ denote the quantile corresponding to $\theta_{t}$.  Compensation needs to be provided only if $q_{t}^{\theta} \in (\min\{p_{t}^{\textsf{CwG}}, p_{t}^{\textsf{HTG}}\}, \max\{p_{t}^{\textsf{CwG}}, p_{t}^{\textsf{HTG}}\})$; let $\partial R_{t}(B_{t}, a_{t}^{\textsf{CwG}})$ denote the compensation. Then we have that $\partial R_{t}(B_{t}, a_{t}^{\textsf{CwG}}) \leq \max\{ F^{-1}(p_{t}^{\textsf{CwG}}) - F^{-1}(p_{t}^{\textsf{HTG}}), F^{-1}(p_{t}^{\textsf{HTG}}) - F^{-1}((p_{t}^{\textsf{CwG}})^{+}) \}$.
\end{lemma}

\subsection{Formal Proof of Theorem \ref{thm:general-cwg-upper-bound}}
\label{onapp:formal-proof-cwg-upper-bound}
\proof{\underline{Proof of Theorem \ref{thm:general-cwg-upper-bound}}.}
Define $\tilde{T} = T - \floor{64 \log (1 / \varepsilon_{0}) / \varepsilon_{0}^2}$ and define $\tau_0 = \floor{64 \log (1 / \varepsilon_{0}) / \varepsilon_{0}^2}$. Consider some time $t \leq \tilde{T}$ and let $\tau = T - t + 1$ denote the remaining time. Recall that $p_{t}^{\sf CwG} \in \mathcal{F}_t$ and $p_{t}^{\sf HTG}$ depends on the candidate abilities $\vect{\theta}_{\geq t}$ but only via the $B_t$-th largest quantile $q_t^u$ and $B_t + 1$-th largest quantile $q_t^l$. To facilitate our analysis, we employ the so-called principle of deferred decisions, and only reveal $q_t^u$ and $q_t^l$ (in addition to the history up to time $t$ i.e. $\mathcal{F}_t$), which uniquely determines $p_t^{\sf HTG}$. Define the event $\mathcal{L}_t \triangleq \{p_{t}^{\sf CwG} \leq q_t^{l}\}$ and $\mathcal{H}_t \triangleq \{p_{t}^{\sf CwG} \geq q_t^u\} $. For the rest of the proof, we will condition on the event $\mathcal{L}_t$ and prove an upper bound on the expected compensation $\partial \mathcal{R}_{t}(B_t, a_t)$ (conditional on $\mathcal{L}_t$). A similar bound can be analogously shown under the event $\mathcal{H}_t$ and we omit the details to avoid repetition. Let $q^{\theta}_{t}$ denote the quantile corresponding to the candidate ability $\theta_{t}$. Now compensation is needed only if $q^{\theta}_t \in \left[p_{t}^{\sf CwG}, q_t^u\right]$. Let $\mathcal{C}_t$ denote the event that compensation is needed i.e., the action under the \textsf{CwG} threshold is different from the action under the \textsf{HTG} threshold. 
Given $q_t^l$ and $q_t^u$, we know that the $\tau$ periods to go include a random subset of $B_t$ quantiles located above $q_t^u$ (these quantiles are i.i.d uniform in $[q_t^u,1]$) and the remaining $\tau - B_t$ quantiles are below $q_t^l$ (these quantiles are i.i.d uniform in $[0,q_t^l]$). If $p_t^{\sf CwG} \in (q_t^l, q_t^u)$, no compensation is needed. Compensation is needed only if $q^{\theta}_t \in [p_t^{\sf CwG}, q_t^l]$ and this event occurs if (a) the realized quantile $q^{\theta}_t = q_t^l$ or (b) $q^{\theta}_t \in [p_t^{\sf CwG}, q_t^l]$ is one of the $\tau - B_t - 1$ lower quantiles. The probability of case (a) is $1/\tau$ and probabilty of (b) is $(\tau - B_t - 1)(q_t^l - p_{t}^{\sf CwG}) / q_t^l \tau$. Combining the two we have that
\begin{align}
        \label{eq:compensation-probability}
        \mathbb{P}\left( \mathcal{C}_t | \mathcal{F}_{t}, q_t^l, q_t^u, \mathcal{L}_t \right) = \frac{\mathbbm{1}_{\{p_{t}^{\sf CwG} \leq q_t^l\}} }{\tau} + \frac{(\tau - B_t - 1)(q^l_t - p_t^{\sf CwG})_{+}}{q_t^l \tau}
    \end{align}
    where $\tau = T - t + 1$ and $(x)_+ = \max\{x, 0\}$. Using Lemma \ref{lem:compensation-upper-bound} and \eqref{eq:compensation-probability}, we have the following bound the expected compensation
    \begin{align}
        \nonumber
        \mathbb{E}\left[\partial \mathcal{R}_t(B_t^{\textsf{CwG}}, a_t^{\sf CwG}, \vect{\theta}_{\geq t}) | \mathcal{F}_t, q_t^l, q_t^u, {\cal L}_t \right] & \leq \frac{|F^{-1}(q_t^u) - F^{-1}((p_t^{\sf CwG})^+)|}{\tau}  \\
        \label{eq:compensation-general-bound}
        & + \frac{(q_t^l - p_t^{\sf CwG}) |F^{-1}(q_t^u) - F^{-1}((p_t^{\sf CwG})^+)|(\tau - B_t - 1)}{q_t^l \tau} 
    \end{align}

Next we need to bound the ratio $(\tau - B_t - 1) / (q_t^l \tau)$ and at any time $t \leq T - \tau_0$, exactly of the following complementary events occurs: (a) ${\cal E}_t = \{1 - B_t / \tau > 4 \sqrt{\log \tau / \tau }\}$ and (b) ${\cal E}_t^c = \{1 - B_t / \tau \leq 4 \sqrt{\log \tau / \tau}\}$. Recall the event ${\cal A}_{3,t}$ defined in \ref{eq:event_cwg_HTG_same_quantile}, which states that 
\begin{align*}
    {\cal A}_{3,t} &= \{q_t^l, q_t^u, p_t^{\sf CwG} \text{ are quantiles belonging (weakly) to the same cluster}\}
\end{align*}
Next we will establish an upper bound on \eqref{eq:compensation-general-bound} for each of the events $\mathcal{E}_t$ and $\mathcal{E}_{t}^c$.
\paragraph{\underline{Case (a): $\mathcal{E}_{t} = \{1 - B_{t} / \tau > 4 \sqrt{\log \tau / \tau}\}$}.} 
Define the following events:
\begin{align*}
\mathcal{A}_{4,t} &\triangleq \{q^{l}_{t} \geq (1/2)(1 - B_{t} / \tau)\} \, ,\\
\mathcal{A}_{5,t} &\triangleq {\cal A}_{3,t} \cap {\cal A}_{4,t} \, .
\end{align*}
Under the event $\mathcal{A}_{3,t}$, from Definition \ref{def:beta-epsilon-delta-clustered} (a) it follows that $|F^{-1}(q^{u}_{t}) - F^{-1}((p_{t}^{\textsf{CwG}})^{+})| \leq |q^{u}_{t} - p_{t}^{\textsf{CwG}}|^{\frac{1}{1+\beta}} + \delta$. Now, on the event ${\cal A}_{4,t}$, we have that $(\tau - B_{t} - 1) / (\tau q_{t}^{l}) \leq 2$. 
We have that
\begin{align*}
    \mathbb{P}({\cal A}_{4,t}^{c}|B_t, \mathcal{E}_{t}) &= \mathbb{P}(q_{t}^{l} < (1/2) (1 - B_{t}/ \tau) ) 
    \leq \mathbb{P}(\textup{Binomial}(\tau, (1/2) (1 - B_{t}/\tau)^{-}) \geq \tau - B_t - 1 ) \\
    &\leq \exp(-\Omega(\tau - B_{t})) \leq C / (\tau - B_{t})^{4}  \leq C / \tau^{2} \, .
\end{align*}
where the last inequality follows from the case assumption that $\tau - B_{t} \geq 4\sqrt{\tau \log \tau}$ and the inequality is true for some appropriately defined constant $C < \infty$.
It follows that 
\begin{align}
    \mathbb{P}({\cal A}_{4,t}^{c}|\mathcal{E}_{t}) \leq C / \tau^{2} \, .
      \label{eq:bound-A-4}
\end{align}

Using \eqref{eq:compensation-general-bound}, and the definitions of the events ${\cal A}_{3,t}$ and ${\cal A}_{4,t}$, we have that 
\begin{align}
    \nonumber
    \mathbb{E}&\left[\partial \mathcal{R}_{t}(B_{t}, a_t^{\sf CwG}, \vect{\theta}_{\geq t}) | \mathcal{F}_{t},q^{l}_{t}, q^{u}_{t}, \mathcal{L}_{t}, \mathcal{E}_{t} \right] \\
    &\leq  \mathbbm{1}_{{\cal A}_{3,t}} \cdot \big[|q^{u}_{t} - p_{t}^{\textsf{CwG}}|^{\frac{1}{1+\beta}} + \delta\big ] / \tau + 2\mathbbm{1}_{{\cal A}_{3,t}} \mathbbm{1}_{{\cal A}_{4,t}}\left[|q^{u}_{t} - p_{t}^{\textsf{CwG}}|^{1+\frac{1}{1+\beta}} + |q^{u}_{t} - p_{t}^{\textsf{CwG}}| \delta \right] + \mathbbm{1}_{{\cal A}_{3,t}^c} + \mathbbm{1}_{{\cal A}_{4,t}^c}, \label{eq:reg-conditional-ub-1}\\
    &\leq |q^{u}_{t} - p_{t}^{\textsf{CwG}}|^{\frac{1}{1+\beta}}/\tau + \delta / \tau + 2|q^{u}_{t} - p_{t}^{\textsf{CwG}}|^{1+\frac{1}{1+\beta}} + 2|q^{u}_{t} - p_{t}^{\textsf{CwG}}| \delta +  \mathbbm{1}_{{\cal A}_{4,t}^c} + \mathbbm{1}_{{\cal A}_{3,t}^{c}},
    \label{eq:Reg-conditional-upperbound}
\end{align}
where the first inequality follows from $q_{t}^{l} - p_{t}^{\textsf{CwG}} \leq q_{t}^{u} - p_{t}^{\textsf{CwG}}$, and the second inequality follows from the fact that $\mathbbm{1}_{{\cal A}_{3,t}}, \mathbbm{1}_{{\cal A}_{4,t}} \leq 1$. Using the definition of the event $A_{2,t}$ in \eqref{eq:event_cwg_HTG_close} and Lemma \ref{lem:crt-HTG-close}, we have that for all $\alpha \in (0,2]$, we have
\begin{align}
    \label{eq:bound-expected-HTG-cwg-distance}
    \mathbb{E}\left[|q_{t}^{u} - p_{t}^{\textsf{CwG}}|^{\alpha} \right] &\leq \mathbb{E}\left[\mathbbm{1}_{{\cal A}_{2,t}}|q_{t}^{u} - p_{t}^{\textsf{CwG}}|^{\alpha}  + \mathbbm{1}_{{\cal A}_{2,t}^{c}}\right] \leq 3^{\alpha} (\log \tau / \tau)^{\alpha /2} + 2 / \tau^{4} \leq C \left(\log \tau / \tau  \right)^{\alpha / 2}\, .
\end{align}
Taking expectations on both sides of \eqref{eq:Reg-conditional-upperbound}, we obtain that 
\begin{align}
    \nonumber
    \mathbb{E}&[\partial \mathcal{R}_{t}(B_{t}, a_t^{\sf CwG}, \vect{\theta}_{\geq t}) | \mathcal{E}_{t}, \mathcal{L}_{t}] \\ \nonumber &\stackrel{(i)}\leq \left(\mathbb{E}\left[|q_{t}^{u} - p_{t}^{\textsf{CwG}}|^{\frac{1}{\beta + 1}} \right] + \delta \right)/ \tau +  2 \mathbb{E}\left[| q_{t}^{u} - p_{t}^{\textsf{CwG}}|^{1 + \frac{1}{\beta + 1}} \right] + 2 \mathbb{E}\left[|q_{t}^{u} - p_{t}^{\textsf{CwG}}| \right] \delta + \mathbb{P}({\cal A}_{3,t}^{c}|\mathcal{E}_{t}) + \mathbb{P}({\cal A}_{4,t}^{c}|\mathcal{E}_{t}), \\
    \label{eq:case-a-regret-decomposition}
    &\stackrel{(ii)}\leq 6 \left(\log \tau / \tau  \right)^{\frac{1}{2(\beta + 1)}} / \tau + \delta / \tau + 36 \left( \log \tau / \tau \right)^{\frac{1}{2} + \frac{1}{2(\beta + 1)}} +  6 \delta \sqrt{\log \tau / \tau} + \mathbb{P}({\cal A}_{3,t}^{c}|\mathcal{E}_t)+  C / \tau^{2},
\end{align}
where inequality (i) follows from the taking expectation on both sides, and inequality (ii) follows from using \eqref{eq:bound-expected-HTG-cwg-distance} for the first, third and fourth summands, and the sixth summand follows from \eqref{eq:bound-A-4}. 

\paragraph{\underline{Case (b): $\mathcal{E}_{t}^c = \{1 - B_{t} / \tau \leq 4 \sqrt{\log \tau / \tau}\}$.}}
The event $1 - B_{t} / \tau \leq 4 \sqrt{\log \tau / \tau}$ implies that $(\tau - B_{t} - 1) / \tau \leq 4 \sqrt{\log \tau / \tau}$, and obviously we have $(q_{t}^{l} - p_{t}^{\textsf{CwG}})  / q_{t}^{l} \leq 1$. Therefore the second term in the RHS of \eqref{eq:compensation-general-bound} is bounded above as
\begin{align*}
    |q^{l}_{t} - p_{t}^{\textsf{CwG}}||F^{-1}(q^{u}_{t}) - F^{-1}((p_{t}^{\textsf{CwG}})^{+})|(\tau-B_{t}-1)/(q^{l}_{t} \tau) \leq  4 \sqrt{\log \tau / \tau} \left| F^{-1}(q^{u}_{t}) - F^{-1}((p_{t}^{\textsf{CwG}})^{+})\right|
\end{align*}
Therefore we can upper bound $\mathbb{E}\left[\partial \mathcal{R}_{t}(B_{t}, \mu_{t}^{\sf CwG}, \vect{\theta}_{\geq t}) | \mathcal{F}_{t},q^{l}_{t}, q^{u}_{t}, \mathcal{E}_{t}^{c}, \mathcal{L}_t \right]$ as 
\begin{align*}
    &\mathbb{E}\left[\partial \mathcal{R}_{t}(B_{t}, a_{t}^{\sf CwG}, \vect{\theta}_{\geq t}) | \mathcal{F}_{t},q^{l}_{t}, q^{u}_{t}, \mathcal{E}_{t}^{c}, \mathcal{L}_t \right] \nonumber\\ 
    &\leq \mathbbm{1}_{A_{3,t}} \left[\left|q_{t}^{u} - p_{t}^{\textsf{CwG}}\right|^{\frac{1}{1 + \beta}} / \tau + \delta / \tau \right] + \mathbbm{1}_{A_{3,t}} \left[ 4 \sqrt{\log \tau / \tau} |q_{t}^{u} - p_{t}^{\textsf{CwG}}|^{\frac{1}{\beta + 1}} + 4\delta \sqrt{\log \tau / \tau}  \right] + \mathbbm{1}_{A_{3,t}^{c}}, \\
    &\leq \left|q_{t}^{u} - p_{t}^{\textsf{CwG}}\right|^{\frac{1}{1 + \beta}} / \tau + \delta / \tau + 4 \sqrt{\log \tau / \tau} |q_{t}^{u} - p_{t}^{\textsf{CwG}}|^{\frac{1}{\beta + 1}} + 4 \delta \sqrt{\log \tau / \tau} +  \mathbbm{1}_{A_{3,t}^{c}},
\end{align*}
Taking expectations on both sides we get that 
\begin{align}
    \nonumber
    &\mathbb{E}\left[\partial \mathcal{R}_{t}(B_{t}, a_{t}^{\sf CwG}, \vect{\theta}_{\geq t}) | \mathcal{E}_{t}^{c}, \mathcal{L}_t \right] \\
    &\leq \mathbb{E}\left[|q_{t}^{u} - p_{t}^{\texttt{cwg}}|^{\frac{1}{1 + \beta}} \right] / \tau + \delta / \tau + 4 \sqrt{\log \tau / \tau} \mathbb{E}\left[\left|q_{t}^{u} - p_{t}^{\texttt{cwg}} \right|^{\frac{1}{1 + \beta}} \right] + 4 \delta \sqrt{\log \tau / \tau} + \mathbb{P}(\mathcal{A}_{3,t}^{c} | \mathcal{E}_{t}^c), \nonumber\\
    \label{eq:case-b-regret-decomposition}
    &\leq  6\left( \log \tau / \tau \right)^{\frac{1}{2(1 + \beta)}} / \tau + \delta / \tau + 24 \left( \log \tau / \tau \right)^{\frac{1}{2} + \frac{1}{2(1 + \beta)}} + 4 \delta \sqrt{\log \tau / \tau} + \mathbb{P}(\mathcal{A}_{3,t}^{c} | \mathcal{E}_{t}^c),
\end{align}
where the second inequality follows from the fact that the first and the second term are bounded by \eqref{eq:bound-expected-HTG-cwg-distance}. This completes for event $\mathcal{E}_t^c$. From Lemma \ref{lem:opposite-side-crt-HTG}, it follows that $\mathbb{P}(\mathcal{A}_{3,t}^{c}) \leq 2n(n +1) / \tau^{4}$.
\begin{align}
    \label{eq:bound-A-3-complement-sum}
    \sum_{t = 1}^{\tilde{T}} \mathbb{P}({\cal A}_{3,t}^{c}) \leq \sum_{\tau = \tau_{0}}^{T} 2n(n + 1) / \tau^{4} \leq n(n + 1) /  \tau_{0}^{3} \leq n(n + 1) \varepsilon_{0}^{6} \stackrel{(\star)}\leq 1
\end{align}
where ($\star$) follows since if there are $n$ gaps, there are $n+1$ clusters, and hence $\varepsilon_{0} \leq 1/ (n+1)$.
Combining  \eqref{eq:case-a-regret-decomposition} and \eqref{eq:case-b-regret-decomposition}, for a constant $C < \infty$ we have that
\begin{align*}
    \mathbb{E}\left[ \partial \mathcal{R}_{t}(B_{t}, a_{t}^{\sf CwG}, \vect{\theta}_{\geq t}) | \mathcal{L}_t\right] \leq C (\log \tau / \tau)^{\frac{1}{2} + \frac{1}{2(\beta + 1)}}   +  C \delta \sqrt{\log \tau / \tau} + \mathbb{P}({\cal A}_{3,t}^{c}) \, .
\end{align*}
An identical bound holds for the regret contribution from the event $\mathcal{H}_t$ where $\mathcal{H}_t = \{p_t^{\textup{cwg}}\geq q_t^u\}$, by a symmetric argument. As a result, we can bound the expected total regret at time $t$ as per $\mathbb{E}\left[ \partial \mathcal{R}_{t}(B_t, a_t^{\sf CwG}, \theta_{\geq t})\right] \leq 2\mathbb{E}\left[ \partial \mathcal{R}_{t}(B_t, a_t^{\sf CwG}, \theta_{\geq t}) | \mathcal{L}_{t} \right]$. Therefore we have that there exists a constant $C < \infty$ such that 
\begin{align}
    \label{eq:compensation-cwg-bound}
    \partial \mathcal{R}_{t}(B_t, a_t^{\sf CwG}) \leq C (\log \tau / \tau)^{\frac{1}{2} + \frac{1}{2(\beta + 1)}}   +  C \delta \sqrt{\log \tau / \tau} + \mathbb{P}(A_{3,t}^{c})
\end{align}
Note that the RHS for \eqref{eq:compensation-cwg-bound} does not depend on the remaining budget $B_t$ and hence we have a uniform bound on the expected compensation given below.
\begin{align}
    \label{eq:compensation-cwg-bound-uniform}
    \sup_{B_{t} \geq 0} \partial \mathcal{R}_{t}(B_t, a_t^{\sf CwG}) \leq C (\log \tau / \tau)^{\frac{1}{2} + \frac{1}{2(\beta + 1)}}   +  C \delta \sqrt{\log \tau / \tau} + \mathbb{P}({\cal A}_{3,t}^{c})
\end{align}
This further implies that 
\begin{align}
    \label{eq:compensation-cwg-bound-expected}
    \mathbb{E}_{B_{t}}\left[\partial \mathcal{R}_{t}(B_t, a_t^{\sf CwG})\right] \leq C (\log \tau / \tau)^{\frac{1}{2} + \frac{1}{2(\beta + 1)}}   +  C \delta \sqrt{\log \tau / \tau} + \mathbb{P}({\cal A}_{3,t}^{c})
\end{align}
Now summing this bound from $t = 1$ to $t = \tilde{T}$, we have, using \eqref{eq:bound-A-3-complement-sum} that for a constant $C < \infty$, we have that
\begin{align*}
    \sum_{t = 1}^{\tilde{T}} \mathbb{E}_{B_t}[\partial \mathcal{R}_{t}(B_t, a_t^{\sf CwG})] &\leq C\left[(1 + 1 / \beta) (\log T)^{\frac{1}{2} + \frac{1}{2(\beta + 1)}} T^{\frac{1}{2} - \frac{1}{2(1 + \beta)}} \cdot \mathbbm{1}\{\beta > 0\} + (\log T)^{2} \mathbbm{1}\{\beta = 0\} \right] \\ & \quad  + C \delta \sqrt{T \log T}\, .
\end{align*}
Finally, consider time steps $t$ such that $\tilde{T} + 1 \leq t \leq T$. In the last $64\log(1 / \varepsilon_{0}) / \varepsilon_{0}^{2}$ time steps, we make use of the static allocation policy and as noted before the regret accrued during the static allocation policy is upper bounded by $C\sqrt{\tau_{0}} = C\sqrt{\log (1 / \varepsilon_{0})} / \varepsilon_{0}$ for some universal constant $C < \infty$. Combining the two parts completes the proof. \QEDB
\endproof

\subsection{Proof of Helper Lemmas}
\label{onapp:helper-lemma-cwg-proofs}
\proof{\underline{Proof of Lemma \ref{lem:rt-HTG-close}.}}
\label{subsec:ce-HTG-close}
Let us assume that $p_{t}^{\textsf{CE}} \geq p_{t}^{\textsf{HTG}} + \sqrt{2 \log \tau / \tau}$. Now conditional on $B_{t}$ and given the knowledge of $p_{t}^{\textsf{HTG}}$, we know that there are $B_{t}$ arrivals with quantile larger than $p_{t}^{\textsf{HTG}}$ and $\tau - B_{t}$ arrivals with quantiles less than $p_{t}^{\textsf{HTG}}$. Let $X_{t} \triangleq \text{Ber}(\tau, (p_{t}^{\textsf{CE}} - \sqrt{2\log \tau / \tau})^{+})$ with $\mathbb{E}\left[X_{t} | B_{t}\right] = (\tau - B_{t} - \sqrt{2 \tau \log \tau})^{+}$. Then we have that
\begin{align*}
    \mathbb{P}\left(p_{t}^{\textsf{CE}} \geq p_{t}^{\textsf{HTG}} + \sqrt{2 \log \tau / \tau} \bigg| B_{t}\right) &\leq \mathbb{P}\left( X_{t} \geq \tau - B_{t} \bigg| B_{t}\right) \leq \mathbb{P}\left(X_{t} - \mathbb{E}\left[ X_{t} | B_{t}\right]  \geq \sqrt{2 \log \tau / \tau} \bigg| B_{t} \right) \leq {1}/{\tau^{4}}
\end{align*}
where the last inequality follows from the Hoeffding inequality. It follows that $\mathbb{P}(p_{t}^{\textsf{CE}} \geq p_{t}^{\textsf{HTG}} + \sqrt{2 \log \tau / \tau}) \leq 1 / \tau^{4}$. Analogously, we can show the same for the case of $p_{t}^{\textsf{HTG}} \geq p_{t}^{\textsf{CE}} + \sqrt{2 \log \tau / \tau}$. \QEDB
\endproof

\proof{\underline{Proof of Lemma \ref{lem:crt-HTG-close}.}}
We have that $| p_{t}^{\textsf{CwG}} - p_{t}^{\textsf{HTG}} | = | p_{t}^{\textsf{CwG}} - p_{t}^{\textsf{CE}} + p_{t}^{\textsf{CE}} - p_{t}^{\textsf{HTG}} | \leq |p_{t}^{\textsf{CwG}} - p_{t}^{\textsf{CE}} | + |p_{t}^{\textsf{CE}} - p_{t}^{\textsf{HTG}} |$. By the definition of the algorithm we have that $|p_{t}^{\textsf{CwG}} - p_{t}^{\textsf{CE}}| \leq \sqrt{2\log\tau / \tau}$. Now conditional on $B_{t}$, we have that event ${\cal A}_{1,t}$ implies the event ${\cal A}_{2,t}$ and hence we have that ${\cal A}_{2,t}^c$ implies ${\cal A}_{1,t}^c$ which implies that $\mathbb{P}({\cal A}_{2,t}^c | B_{t}) \leq \mathbb{P}({\cal A}_{1,t}^c | B_{t})$. Using the proof of Lemma \ref{lem:rt-HTG-close}, we have that $\mathbb{P}({\cal A}_{2,t}^{c} | B_{t}) \leq 2/\tau^{4}$ and the claim of the lemma follows. \QEDB
\endproof

\proof{\underline{Proof of Lemma \ref{lem:opposite-side-crt-HTG}.}}
We have that ${\cal A}_{3,t}^{c} = \cup_{i,j: Q_{i}^{\circ} \cap Q_{j}^{\circ} = \emptyset} \{p_{t}^{\sf CwG} \in Q_{i}, p_{t}^{\sf HTG} \in Q_{j} \}$ where $A^{\circ}$ denotes the interior of the set $A$. Consider the event $\{p_{t}^{\sf CwG} \in Q_{i}, p_{t}^{\sf HTG} \in Q_{j} \}$ such that $Q_{i}^{\circ} \cap Q_{j}^{\circ} = \emptyset$. From the definition of $p_{t}^{\sf CwG}$ in Algorithm \ref{alg:CwG} and the fact that $Q_{i}^{\circ} \cap Q_{j}^{\circ} = \emptyset$ implies that $|p_{t}^{\sf CE} - p_{t}^{\sf HTG}| \geq \sqrt{2\log \tau / \tau}$. Using Lemma \ref{lem:rt-HTG-close} and the union bound completes the proof. \QEDB
\endproof

\proof{\underline{Proof of Lemma \ref{lem:compensation-upper-bound}.}}
Assume that $p_{t}^{\sf CwG} \leq q_{t}^{l}$, then according to the definition of $p_{t}^{\sf HTG}$, we have that $p_{t}^{\sf HTG} = q_{t}^{u}$. Compensation is provided only if $q^{\theta}_{t} \in [p_{t}^{\sf CwG}, p_{t}^{\sf HTG}]$. Suppose that is the case, then we have that $F^{-1}((p_{t}^{\sf CwG})^{+}) \leq \theta_{t} \leq F^{-1}(p_{t}^{\sf CwG}) = F^{-1}(q_{t}^{u}) = u_{t}$. The {\sf CwG} policy would accept the candidate with ability $\theta_{t}$ since $\theta_{t} \geq F^{-1}((p_{t}^{\sf CwG})^{+})$ where as the HTG would want to reject the candidate with ability $\theta_{t}$, because in the future it knows that it can select a candidate with ability at least $u_{t} \geq \theta_{t}$. Hence to persuade the HTG, we need to compensate it $u_{t} - \theta_{t} = F^{-1}(p_{t}^{\sf HTG}) - \theta_{t}$ and maximum compensation can hence be $F^{-1}(p_{t}^{\sf HTG}) - F^{-1}((p_{t}^{\sf CwG})^{+})$. An analogous analysis can be done for the case when $p_{t}^{\sf CwG} \geq q_{t}^{u}$ which follows similarly. \QEDB
\endproof

\section{Proof of Corollary \ref{cor:discrete-distribution-bound}}
\label{sec:proof-cwg-discrete-type}
\proof{\underline{Proof of Corollary \ref{cor:discrete-distribution-bound}.}} The discrete distribution considered is a $(\beta=0, \varepsilon_0, \delta=0)$-clustered distribution for $\varepsilon_{0} = \min_{1 \leq i \leq m}\{f_{i}\}$.
As done for the general case above, our analysis for the case of discrete distributions as considered in the Example \ref{ex:discret-distributions} also follows in two parts. The regret accrued during the second part due to the static allocation policy is upper bounded by $C\sqrt{\log(1 / \varepsilon_{0})} / \varepsilon_{0}$ for some universal constant $C < \infty$. Next we will consider the first part. The argument for the first part will mirror the analysis presented in the proof of Theorem \ref{thm:general-cwg-upper-bound} except for one important improvement we make for this special case. Consider the regret contribution of sample paths satisfying $\mathcal{L}_t \triangleq \{p_{t}^{\sf CwG} \leq q_{t}^{l}\}$ as we did previously. (Again, there is a analogous analysis for the symmetric event $\mathcal{H}_t \triangleq \{p_{t}^{\sf CwG} > q_{t}^{u}\}$, which we omit to avoid repetition.) The only but important distinction in the case of discrete distributions is that on the event ${\cal A}_{3,t}$, which is that $q_{t}^{l}, q_{t}^{u}$ and $p_{t}^{\sf CwG}$ are quantiles belonging to the same cluster, the compensation is given as $F^{-1}(q_{t}^{u}) - F^{-1}((p_{t}^{\sf CwG})^{+})$, however for discrete distributions, we have that $F^{-1}(q_{t}^{u}) - F^{-1}((p_{t}^{\sf CwG})^{+}) = 0$. Previously, in the general case, we had upper bounded $F^{-1}(q_{t}^{u}) - F^{-1}((p_{t}^{\sf CwG})^{+})$ by $|q^{u}_{t} - p_{t}^{\sf CwG}|^{{1} / {(1+\beta)}} + \delta$ using Definition \ref{def:beta-epsilon-delta-clustered}. Because $F^{-1}(q_{t}^{u}) - F^{-1}((p_{t}^{\sf CwG})^{+}) = 0$ on the event $\mathcal{A}_{3,t}$, from \eqref{eq:compensation-general-bound}, we have that 
\begin{align*}
    \mathbb{E}\left[\partial \mathcal{R}_{t}(B_{t}, a_t^{\sf CwG}, \vect{\theta}_{\geq t}) | \mathcal{F}_{t},q^{l}_{t}, q^{u}_{t}, \mathcal{L}_{t}, \mathcal{E}_{t} \right] 
    \leq \mathbbm{1}_{{\cal A}_{3,t}^c} + \mathbbm{1}_{{\cal A}_{4,t}^c},
\end{align*}
because $\sup_{B_t, a_t, \bm{\theta}_{\geq t}}\partial \mathcal{R}_t(B_t, a_t; \bm{\theta}_{\geq t}) \leq 1$. Taking expectations on both sides, we have that 
\begin{align*}
    \mathbb{E}\left[\partial \mathcal{R}_{t}(B_{t}, a_t^{\sf CwG}, \vect{\theta}_{\geq t}) | \mathcal{L}_{t} \right] 
    \leq \mathbb{P}({\cal A}_{3,t}^c) + \mathbb{P}({\cal A}_{4,t}^c) \leq \mathbb{P}({\cal A}_{3,t}^{c}) + C / \tau^{2},
\end{align*}
 Summing this upper bound from $\tau = \tau_{0}$ to $\tau = T$, we get that the summation is upper bounded by a universal constant $C < \infty$ using \eqref{eq:bound-A-3-complement-sum}. Combining the regret accrued in the two parts, we get the required result. \QEDB
\endproof

\section{Proof of Corollary \ref{cor:beta-one-clustered-regret}}
\label{app:proof-corollary-beta-one-clustered}
\proof{\underline{Proof of Corollary \ref{cor:beta-one-clustered-regret}. }}
Since by assumption, there are no gaps in the distribution, we have that $p_{t}^{\sf CwG} = p_{t}^{\sf CE}$ for all $t$. Our analysis will follow along the same lines as the analysis for Theorem \ref{thm:general-cwg-upper-bound} with $\delta = 0$.
From \eqref{eq:compensation-general-bound}, we have that 
\begin{align}
        \nonumber
        \mathbb{E}\left[\partial \mathcal{R}_t(B_t, a_t^{\sf CwG}, \vect{\theta}_{\geq t}) | \mathcal{F}_t, q_t^l, q_t^u, {\cal L}_t \right] & \leq \frac{|F^{-1}(q_t^u) - F^{-1}((p_t^{\sf CwG})^+)|}{\tau}  \\
        \nonumber
        & + \frac{(q_t^l - p_t^{\sf CwG}) |F^{-1}(q_t^u) - F^{-1}((p_t^{\sf CwG})^+)|(\tau - B_t - 1)}{q_t^l \tau} 
\end{align}
This is where our proof departs from the proof of Theorem \ref{thm:general-cwg-upper-bound}. The fact that the {\sf CwG} policy boils down to the {\sf CE} policy when there are no non-trivial gaps simplifies the analysis to a great extent. Instead of considering two cases to bound the ratio $(\tau - B_{t} - 1) / (q_{t}^{l} \tau)$, we can bound it much simply. From the definition of $p_{t}^{\sf CE}$, we have that $p_{t}^{\sf CE} = 1 - B_{t} / \tau$, which implies that $(\tau - B_{t} - 1) / (q_{t}^{l} \tau) \leq p_{t}^{\sf CE} \tau / (q_{t}^{l} \tau) \leq 1$ since $p_{t}^{\sf CwG} = p_{t}^{\sf CE}$ and we are considering the sample paths on which $p_{t}^{\sf CwG} \leq q_{t}^{l}$. Since $F$ is a $(\beta, \varepsilon_0 = 1, \delta = 0)$-clustered distribution, we have that 
\begin{align*}
    \mathbb{E}\left[\partial \mathcal{R}_t(B_t, a_t^{\sf CwG}, \vect{\theta}_{\geq t}) | \mathcal{F}_t, q_t^l, q_t^u, {\cal L}_t \right]  &\leq \frac{|q_{t}^{u} - p_{t}^{\sf CwG}|^{\frac{1}{1 + \beta}}}{\tau} + |q_{t}^{u} - p_{t}^{\sf CwG}|^{1 + \frac{1}{1 + \beta}} 
\end{align*}

Taking expectations, we have that
\begin{align*}
    \mathbb{E}\left[\partial \mathcal{R}_t(B_t, a_t^{\sf CwG}, \vect{\theta}_{\geq t}) |  L_t \right] &\leq \mathbb{E}\left[ |q_{t}^{u} - p_{t}^{\sf CwG}|^{\frac{1}{1 + \beta}} \right] / \tau + \mathbb{E}\left[|q_{t}^{u} - p_{t}^{\sf CwG}|^{1 + \frac{1}{1 + \beta}} \right] \\ 
    &\leq C \left(\tau^{-\frac{1}{2(1 + \beta)} - 1} + \tau^{-\frac{1}{2(1 + \beta)} - \frac{1}{2}}\right),
\end{align*}
where the second inequality follows from the fact that $\mathbb{E}\left[|q_{t}^{u} - p_{t}^{\sf CwG}|^{\alpha} \right] \leq C(T - t + 1)^{-\alpha / 2}$ for any $\alpha \in (0, 2]$ and the fact that $1/(1 + \beta) \in (0,1]$ and $1 + 1/(1 + \beta) \in (0,2]$. Recall that $\mathbb{E}\left[\partial \mathcal{R}_t(B_t, a_t^{\sf CwG}, \vect{\theta}_{\geq t})\right] \leq 2 \mathbb{E}\left[\partial \mathcal{R}_t(B_t, a_t^{\sf CwG}, \vect{\theta}_{\geq t}) |  {\cal L}_t \right] $. Summing this over $T$ time steps gives us the regret scaling in Corollary \ref{cor:beta-one-clustered-regret}. To complete the proof, we will prove that for all $t \leq T - 1$ and $\alpha \in (0,2]$, we have that $\mathbb{E}\left[ |q_{t}^{u} - p_{t}^{\sf CwG}|^{\alpha} \right] \leq C(T - t + 1)^{-\alpha / 2}$ . This inequality follows from the Hoeffding inequality as shown below.
\begin{align*}
    \mathbb{E}\left[|q_{t}^{u} - p_{t}^{\sf CwG}|^{\alpha} \right] = \int_{0}^{\infty} \mathbb{P}\left(|q_{t}^{u} - p_{t}^{\sf CwG}|^{\alpha} \geq x \right) dx 
    \leq 2 \int_{0}^{\infty} \exp\left(-2\tau x^{2/\alpha} \right) dx 
    = 2^{-\alpha / 2} \alpha \Gamma(1 + 2 / \alpha) \tau^{-\alpha / 2},
\end{align*}
This completes the proof of Corollary \ref{cor:beta-one-clustered-regret}.
\QEDB

\section{Recovering existing regret guarantees for {\sf RAMS}}
\label{onapp:recovering-existing-guarantees-rams}
\simulation{
In this section, we provide corollaries which show that {\sf RAMS} attains near-optimal regret scaling for a variety of online resource allocation problems under different assumptions. As a first application of Theorem \ref{thm:meta-performance-rams}, we consider the multisecretary problem. For analytical simplicity, we consider a minor variant of {\sf RAMS} where for the first $\tilde{T} = T - \floor{64 \log(1/\varepsilon_0) / \varepsilon_{0}^2}$ time steps, we implement {\sf RAMS} as stated in Algorithm \ref{alg:rams} and in the final $\ceil{64 \log(1/\varepsilon_0) / \varepsilon_{0}^2}$ time steps, we implement a static threshold policy as done in the case of Algorithm \ref{alg:CwG}. This minor variant of {\sf RAMS} inherits the guarantees established in Theorem \ref{thm:general-cwg-upper-bound}.
\begin{corollary}[$\beta$-dependent regret for multisecretary]
    \label{cor:cwg-rams-guarantee} Consider the multisecretary problem with the candidate ability distribution $F$ being $(\beta, \varepsilon_{0}, \delta)$-clustered for some fixed $\beta \in [0, \infty),\varepsilon_{0} \in (0,1]$ and $\delta \in [0,1]$. Fix the parameter $\eta > 2$ in Theorem \ref{thm:meta-performance-rams}. Assume that the number of sample paths drawn at time $t$ is sufficiently large, specifically $K_t \geq (T - t + 1)^{\eta + \nu}$ for some  $\nu > 0$. 
    Then there exists a constant $C \equiv C(F, \eta, \nu) < \infty$, such that for all $T \in \mathbb{N}$ and $B \in \mathbb{N}$, the regret for \textup{\sf RAMS} is bounded above as 
    \begin{align*}
        \text{Regret}(B,T; \textup{\sf RAMS}) &\leq C(1 + 1 / \beta) (\log T)^{\frac{1}{2} + \frac{1}{2(\beta + 1)}} T^{\frac{1}{2} - \frac{1}{2(\beta + 1)}} \cdot \mathbbm{1}\{\beta > 0\} + C(\log T)^{2} \mathbbm{1}\{\beta = 0\} \\
    &\quad  + C \delta \sqrt{T \log T} + C\sqrt{\log(1 / \varepsilon_{0})} / \varepsilon_{0}. 
    \end{align*}
\end{corollary}
}

\simulation{
Next we zoom out from the multisecretary problem and consider the more general network revenue management and online matching problems. We present four assumptions under which these problems have been studied. These assumptions are stated in the notation introduced in this paper.
}

\simulation{
\begin{assumption}[Small number of types for NRM]
    \label{ass:nrm-finite-types}
    The type distribution $F$ is supported on a discrete set
    $\{(r_1, \bm{c}_{1}), (r_2, \bm{c}_2), \dots, (r_n, \bm{c}_n)\}$ with $c_{\theta,k} \in \{0,1\}$ for all $\theta \in \{1,\dots, n\},k \in \{1, \dots, d\}$. 
\end{assumption}
}

\simulation{
\begin{assumption}[Infinitely many types for NRM with density bounded below]
    \label{ass:nrm-infinite-types}
    The consumption random vector $\bm{c}_\theta$ is bounded i.e. $\underline{\nu} \leq \|\bm{c}_{\theta}\|_\infty \leq \bar{\nu}$ for $0 < \underline{\nu} \leq \bar{\nu} < \infty$ for all $\theta \in \Theta$.
    Conditional on the consumption vector $\bm{c}_\theta$, the reward distribution $F_{\theta}$ is assumed to be $(\beta = 0, \varepsilon_{0} = 1)$-clustered with reward random variable $r_\theta$ being bounded in $[0,1]$.
\end{assumption}
}

\simulation{
\begin{assumption}[Infinitely many types for NRM]
    \label{ass:nrm-semi-infinite-types}
    The consumption random vector $\bm{c}_\theta$ is supported on a small discrete set 
    $\{\bm{c}_1, \dots, \bm{c}_{n}\}$ with $c_{\theta,k} \in \{0,1\}$ for all $\theta \in \{1,\dots, n\}$ and $k \in \{1,\dots,d\}$. Conditional on the consumption vector $\bm{c}_\theta$, the reward distribution $F_{\theta}$ is assumed to be $(\beta = 0, \varepsilon_{0})$-clustered distribution with $\varepsilon_0 \in (0,1]$ and the reward random variable $r_\theta$ being bounded in $[0,1]$.
\end{assumption}
}

\simulation{
\begin{assumption}[Small number of types for Online Matching]
    \label{ass:online-matching-finite-types}
    The type distribution $F$ is supported on a discrete set of reward vectors $\{\bm{r}_{1}, \dots, \bm{r}_n\}$ where $\bm{r}_\theta \in [0,1]^d$ for all $\theta \in \{1, \dots, n\}$.
\end{assumption}
}

\simulation{
{\it Discussion of the assumptions.} Recall Assumptions \ref{ass:ms-finite-types} (a few discrete types) and \ref{ass:ms-infinite-types} (continuous types) for the multisecretary problem. Assumptions \ref{ass:nrm-finite-types} and \ref{ass:online-matching-finite-types} are a natural generalization of Assumption \ref{ass:ms-finite-types} in the context of network revenue management and online matching respectively and is often a standard assumption in this literature \citep{vera2021bayesian, bumpensanti2020re, jasin2012re}. Assumptions \ref{ass:nrm-infinite-types} and \ref{ass:nrm-semi-infinite-types} are a generalization of Assumption \ref{ass:ms-infinite-types} for the NRM problem with multiple resources. Assumption \ref{ass:nrm-infinite-types} resembles the assumption studied in \cite{bray2022logarithmic}, however Assumption \ref{ass:nrm-infinite-types} is stronger than the one in \cite{bray2022logarithmic} in the sense that \cite{bray2022logarithmic} allows for arbitrarily small consumption vectors (i.e., $\underline{\nu} = 0$) while Assumption \ref{ass:nrm-infinite-types} assumes that consumption vectors are bounded below. 
Additionally \cite{bray2022logarithmic} allows for unbounded rewards while Assumption \ref{ass:nrm-infinite-types} assumes that the rewards are bounded in the interval $[0,1]$.
Note that we study a stronger version of the assumptions in \cite{bray2022logarithmic} for the sake of technical simplicity and conjecture that {\sf RAMS} will achieve the same logarithmic regret scaling under the assumptions studied in \cite{bray2022logarithmic}.
The key similarity between Assumption \ref{ass:nrm-infinite-types} and the assumption studied in \cite{bray2022logarithmic} is that both assumptions imply that the fluid problem is non-degenerate which enables the logarithmic regret scaling. Assumption \ref{ass:nrm-semi-infinite-types}, while being similar to Assumption \ref{ass:nrm-infinite-types}, allows for degeneracy in the fluid problem and was recently studied by \cite{jiang2022degeneracy}. There are two key distinctions between Assumptions \ref{ass:nrm-infinite-types} and \ref{ass:nrm-semi-infinite-types}: (i) Assumption \ref{ass:nrm-semi-infinite-types} only permits a few consumption types and (ii) Assumption \ref{ass:nrm-semi-infinite-types} allows for gaps in the (conditional) reward distributions which in turn permits degeneracy in the fluid problem.
}

\simulation{
Theorem~\ref{thm:meta-performance-rams} tells us that {\sf RAMS} inherits the regret guarantees previously established for other algorithms, under Assumptions \ref{ass:nrm-finite-types}-\ref{ass:online-matching-finite-types}. This is formalized in the following corollaries. Note that in all our regret guarantees provided below, the only scaling parameters are the time horizon $T$ and the budget $B$ and all other parameters are considered constant. Moreover, we emphasis that the distribution $F$ is initially fixed and its parameters do not scale with the scaling parameter $T$ and $B$. Therefore, the minimum probability parameter $\varepsilon_0$ for the distributions considered in Assumptions \ref{ass:nrm-infinite-types} and \ref{ass:nrm-semi-infinite-types} is also fixed and subsumed in the constants presented below.
}

\simulation{
\begin{corollary}[Regret for NRM]
    \label{cor:nrm-regret-profiles}
    Consider the network revenue management problem with request distribution $F$. Fix the parameter $\eta > 2$ in Theorem \ref{thm:meta-performance-rams}. Assume that the number of sample paths drawn at time $t$ is large enough as per $K_t \geq (T - t + 1)^{\eta + \nu}$ for some $\nu > 0$. 
    Then there exists a constant $C \equiv C(F, \eta, \nu) < \infty$,  
    such that for all $T \in \mathbb{N}$ and $B \in \mathbb{R}^d$, we have that 
    \begin{itemize}
        \item[(a)] \textup{(Constant Regret with few types)} If $F$ satisfies Assumption \ref{ass:nrm-finite-types}, $\text{Regret}(B,T; \textup{\sf RAMS}) \leq C$.
        \item[(b)] \textup{(Logarithmic Regret)} If $F$ satisfies Assumption \ref{ass:nrm-infinite-types}, $\text{Regret}(B,T; \textup{\sf RAMS}) \leq C \log T$.
        \item[(c)] \textup{(Log-Squared Regret)} If $F$ satisfies Assumption \ref{ass:nrm-semi-infinite-types}, $\text{Regret}(B,T; \textup{\sf RAMS}) \leq C \log^2 T$.
    \end{itemize}
\end{corollary}
}

\simulation{
\begin{corollary}[Constant Regret for Online Matching]
    \label{cor:online-matching-constant}
    Consider the online matching setting with request distribution $F$ satisfying Assumption \ref{ass:online-matching-finite-types}. Fix the parameter $\eta > 2$ in Theorem \ref{thm:meta-performance-rams}. Assume that the number of sample paths drawn at time $t$  is large enough as per $K_t \geq (T - t + 1)^{\eta + \nu}$ for some $\nu > 0$. 
    Then there exists a constant $C \equiv C(F, \eta, \nu) < \infty$ 
    such that for all $T \in \mathbb{N}$ and $B \in \mathbb{R}^d$, the regret for \textup{\sf RAMS} is bounded above as
        $\text{Regret}(B,T; \textup{\textsf{RAMS}}) \leq C .$
\end{corollary}
}

\section{Proofs Related to \textsf{RAMS}}
\label{onapp:details-rams}

\new{
\subsection{Proof of Claim \ref{claim:sufficient-condition-compensation-uniform-bound}}
\label{onapp:sufficient-condition-compensation-proof}
\proof{\underline{Proof of Claim \ref{claim:sufficient-condition-compensation-uniform-bound}.}}Given any budget $B_t \geq 0$ and any sample path $\bm{\theta}_{\geq t + 1}$, if the hindsight to go (HTG) policy decides to accept the request $\theta_t$, we can make it reject the request $\theta_t$ by paying a maximum compensation of $r_{\max}$. On the flip side, the hindsight to go policy can extract at most $r_{\max} \bar{\nu} / \underline{\nu}$ in the future for every resource $\theta_t$ makes use of, hence if the hindsight to go (HTG) policy wants to reject $\theta_t$, we can make it accept the request $\theta_t$ by paying a compensation of $d r_{\max} \bar{\nu} / \underline{\nu}$ since the request $\theta_t$ can make use of at most $d$ resources. \QEDB

\subsection{Proof of Lemma \ref{lem:rams-equivalent-marginal-compensation-minimization}}
\label{onapp:proof-rams-marginal-compensation-equivalent}
\proof{\underline{Proof of Lemma \ref{lem:rams-equivalent-marginal-compensation-minimization}.}} Using \eqref{eq:marginal-compensation-sample-path}, for a simulated sample path ${\vect{\theta}}_{\geq t }^{(i)} \triangleq \{\theta_t, \vect{\theta}_{\geq t + 1}^{(i)} \}$, we have that
\begin{align*}
    \partial \mathcal{R}_t(B_t, a, \vect{\theta}_{\geq t}^{(i)}) &= V_t^{\sf hs}(B_t; \vect{\theta}_{\geq t}^{(i)}) - \left[ V_{t + 1}^{\sf hs}(B_t - c(\theta_t, a); \vect{\theta}_{\geq t + 1}^{(i)}) + r(\theta_t, a)\right] \nonumber \\
    &= V_t^{\sf hs}(B_t; \vect{\theta}_{\geq t}^{(i)}) - Q_{t}^{\sf hs}(B_t,a; {\vect{\theta}}_{\geq t}^{(i)}).
\end{align*}
Note that the term $V_t^{\sf hs}(B_t; \vect{\theta}_{\geq t}^{(i)})$ does not depend on the action $a \in \mathcal{A}(B_t,\theta_t)$ and hence we have 
\begin{align}
    \label{eq:rams-equiv-marginal-compensation-minimization}
    \argmax_{a \in \mathcal{A}{(B_t, \theta_t)}} K_t^{-1} \sum_{i = 1}^{K_t} Q_t^{\sf hs}(B_t, a; \tilde{\vect{\theta}}_{\geq t}^{(i)}) = \argmin_{a \in \mathcal{A}{(B_t, \theta_t)}} K_t^{-1} \sum_{i = 1}^{K_t} \partial \mathcal{R}_t(B_t, a; \tilde{\vect{\theta}}_{\geq t}^{(i)})\, ,
\end{align}
i.e., {\sf RAMS} takes an action $a \in \mathcal{A}{(B_t, \theta_t)}$ which minimizes the simulation-based estimate of the expected marginal compensation. \QEDB

\subsection{Proof of Theorem \ref{thm:meta-performance-rams}}
\label{onapp:meta-theorem-rams-proof}

\proof{\underline{Proof of Theorem \ref{thm:meta-performance-rams}.}}Given a budget $B_t$ and a request $\theta_t$ at time $t$, from Algorithm \ref{alg:rams} it follows that the action under the {\sf RAMS} policy is given by:
\begin{align}
    \label{eq:rams-action-max-Q}
    a_t^{\sf RAMS} = \argmax_{a \in \mathcal{A}} \frac{1}{K_t} \sum_{i = 1}^{K_t} Q_t^{\sf hs}\left(B_t, a; \Tilde{\bm{\theta}}_{\geq t}^{(i)} \right),
\end{align}
where $K_t$ denotes the number of simulated sample paths used at time $t$, $Q_t^{\sf hs}\left(B_t, a; \Tilde{\bm{\theta}}_{\geq t}^{(i)} \right)$ is defined in \eqref{eq:rams-opt-Q}, $\Tilde{\bm{\theta}}_{\geq t}^{(i)} \triangleq \{\theta_{t}, \Tilde{\theta}_{t + 1}^{(i)}, \Tilde{\theta}_{t + 2}^{(i)}, \dots, \Tilde{\theta}_{T}^{(i)}\}$ and $\{\Tilde{\theta}_{t + 1}^{(i)}, \Tilde{\theta}_{t + 2}^{(i)}, \dots, \Tilde{\theta}_{T}^{(i)}\}$ denote the $i$-th sequence of simulated sample paths. From Lemma \ref{lem:rams-equivalent-marginal-compensation-minimization}, it follows that the action under the RAMS policy can be equivalently written as
\begin{align}
    \label{eq:rams-action-min-compensation}
    a_t^{\sf RAMS} = \argmin_{a \in \mathcal{A}} \frac{1}{K_t} \sum_{i = 1}^{K_t} \partial \mathcal{R}_t\left(B_{t}, a; \Tilde{\bm{\theta}}_{\geq t}^{(i)} \right),
\end{align}
where $\partial \mathcal{R}_{t}\left(B_t, a; \Tilde{\bm{\theta}}_{\geq t}^{(i)} \right) = \max_{\hat{a} \in \mathcal{A}} Q_{t}^{\sf hs}\left(B_t, \hat{a}; \Tilde{\bm{\theta}}_{\geq t}^{(i)} \right) - Q_{t}^{\sf hs}\left(B_t, a; \Tilde{\bm{\theta}}_{\geq t}^{(i)} \right)$.

From the regret decomposition lemma of \cite{vera2021bayesian}, it follows that
\begin{align}
    \text{Regret}(B,T; {\sf RAMS}) &= \sum_{t = 1}^{T} \mathbb{E}_{B_{t}^{\sf RAMS}}\left[ \partial \mathcal{R}_t\left( B_{t}^{\sf RAMS}, a_{t}^{\sf RAMS}\right)\right].
\end{align}

We note that $\mathbb{E}_{B_{t}^{\sf RAMS}}\left[ \partial \mathcal{R}_{t}\left(B_{t}^{\sf RAMS}, a_{t}^{\sf RAMS} \right) \right] \leq \sup_{B_{t} \geq 0} \partial \mathcal{R}_{t}(B_{t}, a_{t}^{\sf RAMS})$ for all $t \in \{1, 2, \dots, T\}$. Now to prove to Theorem \ref{thm:meta-performance-rams}, it suffices for us to show that 
\begin{align*}
    \sup_{B_{t} \geq 0} \partial \mathcal{R}_t\left( B_t, a_t^{\sf RAMS} \right) \leq \Delta_{t}(\textup{\sf ALG}) + C(\eta, |{\cal A}|, {\cal C}) K_{t}^{-\frac{1}{\eta}},
\end{align*}
where $\Delta_t({\sf ALG})$ is the uniform upper bound assumed in condition $(i)$ at time $t$ under {\sf ALG}. We will begin by upper bounding the quantity $\partial \mathcal{R}_{t}(B_t, a_{t}^{\sf RAMS})$. For some fixed parameter $\eta > 2$, conditional on the budget $B_t$, and request $\theta_t$, define the following ``good" event $\mathcal{G}_t$
\begin{align}
    \label{def:good-event}
    \mathcal{G}_t = \cap_{a \in \mathcal{A}} \Big\{\Big| K_{t}^{-1} \sum_{i = 1}^{K_t} \partial \mathcal{R}_{t}(B_t, a; \tilde{\bm{\theta}}_{\geq t}^{(i)}) - \mathbb{E}\left[ \partial \mathcal{R}_{t}(B_t, a; \bm{\theta}_{\geq t} ) | \theta_{t}, B_t\right] \Big| \leq K_t^{-\frac{1}{\eta}} \Big\}.
\end{align} 
Using the definition of $\partial \mathcal{R}_t(B_t, a_{t}^{\sf RAMS})$ and the tower property we have, 
\begin{align}
    \label{eq:marginal-compensation-iterated-expectation}
    \partial \mathcal{R}_{t}(B_t, a_{t}^{\sf RAMS} ) = \mathbb{E}\left[ \partial \mathcal{R}_t
    (B_t, a_t^{\sf RAMS}; \bm{\theta}_{\geq t} 
    ) \big| B_{t} \right] 
   = \mathbb{E}\left[ \mathbb{E}\left[ \partial \mathcal{R}_t(B_{t}, a_t^{\sf RAMS}; \bm{\theta}_{\geq t}) \big| \theta_{t}, B_{t}\right] \big| B_{t}\right]. 
\end{align}
Now we further write the inner (conditional) expectation $\mathbb{E}\left[ \partial \mathcal{R}_t(B_{t}, a_t^{\sf RAMS}; \bm{\theta}_{\geq t}) \big| \theta_{t}, B_{t}\right]$ as
\begin{align*}
    \mathbb{E}\left[ \partial \mathcal{R}_t(B_{t}, a_t^{\sf RAMS}; \bm{\theta}_{\geq t}) \big| \theta_{t}, B_{t}\right] &= \underbrace{\mathbb{E}\left[ \partial \mathcal{R}_t(B_{t}, a_t^{\sf RAMS}; \bm{\theta}_{\geq t}) \mathbbm{1}_{\mathcal{G}_t} \big| \theta_{t}, B_{t}\right]}_{(\spadesuit)} 
     + \underbrace{\mathbb{E}\left[ \partial \mathcal{R}_t(B_{t}, a_t^{\sf RAMS}; \bm{\theta}_{\geq t}) \mathbbm{1}_{\mathcal{G}_t^c}  \big| \theta_{t}, B_{t}\right]}_{(\clubsuit)}.
\end{align*}
Now we have two terms $(\spadesuit)$ and $(\clubsuit)$ to bound. We begin by bounding the term $(\spadesuit)$. We have that 
\begin{align*}
    \mathbb{E}\left[\partial \mathcal{R}_{t}(B_t, a_t^{\sf RAMS}; \bm{\theta}_{\geq t}) \mathbbm{1}_{\mathcal{G}_t} | \theta_t, B_t \right]  &\stackrel{(a)}\leq \left(K_{t}^{-1} \sum_{i = 1}^{K_t} \partial \mathcal{R}_{t}(B_t, a_t^{\sf RAMS}; \Tilde{\bm{\theta}}_{\geq t}^{(i)}) + K_{t}^{-\frac{1}{\eta}} \right) \mathbbm{1}_{\mathcal{G}_t}, \\
    &\stackrel{(b)}\leq \left(K_t^{-1} \sum_{i = 1}^{K_t} \partial \mathcal{R}_t(B_t, a_{t}^{\sf ALG}; \tilde{\bm{\theta}}_{\geq t}^{(i)}) + K_t^{-\frac{1}{\eta}} \right) \mathbbm{1}_{\mathcal{G}_t} \\
    &\stackrel{(c)}\leq \left(\mathbb{E}\left[\partial \mathcal{R}_t(B_t, a_{t}^{\sf ALG}; {\bm{\theta}}_{\geq t}) | \theta_t, B_t \right] + 2K_t^{-\frac{1}{\eta}} \right) \mathbbm{1}_{\mathcal{G}_t} \\
    &\stackrel{(d)}\leq \mathbb{E}\left[\partial \mathcal{R}_t(B_t, a_{t}^{\sf ALG}; {\bm{\theta}}_{\geq t}) | \theta_t, B_t \right] + 2K_t^{-\frac{1}{\eta}}
\end{align*}
where (a) follows from definition of event $\mathcal{G}_t$ applied to the action $a_t^{\sf RAMS}$, (b) follows from the fact that {\sf RAMS} takes the action according to \eqref{eq:rams-action-min-compensation}, (c) follows from definition of event $\mathcal{G}_t$ applied to the action $a_t^{\sf ALG}$ and (d) follows from that fact that $\mathbbm{1}_{\mathcal{G}_t} \leq 1$. Using this it follows that 
\begin{align*}
    (\spadesuit) = \mathbb{E}\left[\partial \mathcal{R}_{t}(B_t, a_t^{\sf RAMS}; \bm{\theta}_{\geq t}) \mathbbm{1}_{\mathcal{G}_t} | \theta_t, B_t \right] \leq \mathbb{E}\left[\partial \mathcal{R}_{t}(B_{t}, a_t^{\sf ALG}; \bm{\theta}_{\geq t}) | \theta_t, B_t\right] + 2K_t^{-\frac{1}{\eta}}
\end{align*}

Next we bound the term $(\clubsuit)$. Define $Y_i(a) \triangleq \partial \mathcal{R}_{t}(B_t, a; \bm{\theta}_{\geq t}^{(i)})$. From Assumption (ii) in Theorem \ref{thm:meta-performance-rams}, we have that $\partial \mathcal{R}_t(B_t, a; \bm{\theta}_{\geq t}) \leq \mathcal{C}$ almost surely for all $B_t \geq 0, a \in \mathcal{A}$ and $\bm{\theta}_{\geq t}$. Therefore we have that $\{Y_i(a)\}_{i = 1}^{K_t}$ are i.i.d random variables with $|Y_i(a)| \leq \mathcal{C}$ almost surely. Hence we have that
\begin{align*}
    (\clubsuit) \leq \mathcal{C}\mathbb{E}[\mathbbm{1}_{\mathcal{G}_t^c} | \theta_t, B_t] = \mathcal{C} \mathbb{P}\left( \mathcal{G}_t^c | \theta_t, B_t\right) \leq C(\eta, |\mathcal{A}|, \mathcal{C}) K_t^{-\frac{1}{\eta}}
\end{align*}
where the last inequality follows from union bound and Hoeffding's inequality as described below.
\begin{align*}
    \mathbb{P}\left(\mathcal{G}_t^c | \theta_t, B_t \right) &\stackrel{(a)}\leq \sum_{a \in \mathcal{A}} \mathbb{P}\left( \Big|K_t^{-1} \sum_{i = 1}^{K_t} Y_i(a) - \mathbb{E}\left[ Y_i(a) | \theta_t, B_t\right]\Big| > K_t^{-\frac{1}{\eta}}\right), \\
    &\stackrel{(b)}\leq 2 |\mathcal{A}| \exp\left(-\frac{2K_t^{2-\frac{2}{\eta}}}{K_t \mathcal{C}^2} \right), \\
    &\stackrel{(c)}\leq 2 |\mathcal{A}| \exp\left(2 K_t^{\frac{\eta - 2}{\eta}}/{\cal C}^2 \right), \\
    &\stackrel{(d)}\leq C(\eta, |\mathcal{A}|, \mathcal{C})K_t^{-\frac{1}{\eta}},
\end{align*}
where (a) follows from union bound, (b) follows from Hoeffding's inequality, (c) follows trivially, (d) follows for some appropriate constant $C(\eta, |\mathcal{A}|, \cal C)$ since $\exp(-x) \leq C(p)x^{-p}$ for some $p > 0$.
Given the bound on $(\spadesuit)$ and $(\clubsuit)$, we have that 
\begin{align*}
    \mathbb{E}\left[ \partial \mathcal{R}_t(B_t, a_t^{\sf RAMS}; \bm{\theta}_{\geq t}) | \theta_t, B_t \right] \leq \mathbb{E}\left[\partial \mathcal{R}_t(B_t, a_t^{\sf ALG}; \bm{\theta}_{\geq t}) | \theta_t, B_t \right] + C(\eta, |\mathcal{A}|, {\cal C}) K_t^{-\frac{1}{\eta}}.
\end{align*}
Using \eqref{eq:marginal-compensation-iterated-expectation}, we have that 
\begin{align*}
    \partial \mathcal{R}_t(B_t, a_t^{\sf RAMS}) \leq \partial \mathcal{R}_t(B_t, a_t^{\sf ALG}) + C(\eta, |\mathcal{A}|, {\cal C}) K_t^{-\frac{1}{\eta}}.
\end{align*}
Taking a supremum over the budget $B_t \geq \bm{0}$, we have that
\begin{align*}
    \sup_{B_t \geq \bm{0}} \partial \mathcal{R}_{t}(B_t, a_t^{\sf RAMS}) \leq \sup_{B_t \geq \bm{0}} \partial \mathcal{R}_t(B_t, a_t^{\sf ALG}) + C(\eta, |{\cal A}|, {\cal C}(F)) K_t^{-\frac{1}{\eta}} \leq \Delta_t({\sf ALG}) + C(\eta, |{\cal A}|, {\cal C}(F)) K_t^{-\frac{1}{\eta}},
\end{align*}
where the last inequality follows from Assumption $(i)$. This completes the proof. \QEDB
\endproof

\subsection{Proof of Corollaries \ref{cor:cwg-rams-guarantee}, \ref{cor:nrm-regret-profiles} and \ref{cor:online-matching-constant}}
\label{onapp:proof-corollaries}
    From Theorem \ref{thm:meta-performance-rams}, recall that the regret upper bound for {\sf RAMS} (or minor variants of {\sf RAMS}) can be decomposed as a sum of the following two terms
    \begin{align*}
        \text{Regret}(B,T; {\sf RAMS}) \leq \underbrace{\sum_{t = 1}^{T} \Delta_{t}({\sf ALG})}_{(\diamondsuit)} + \underbrace{C \sum_{t = 1}^T K_t^{-\frac{1}{\eta}}}_{(\heartsuit)}
    \end{align*}
    where the constant $C < \infty$ is a function of the parameter $\eta$, size of the action set $|\mathcal{A}|$ and the distribution $F$. Note that $(\heartsuit)$ is common across different problem settings and assumptions while $(\diamondsuit)$ needs to be dealt with separately. For each Corollary \ref{cor:cwg-rams-guarantee}, \ref{cor:nrm-regret-profiles} and \ref{cor:online-matching-constant}, we have that $K_t \geq (T - t + 1)^{\eta + \nu}$ where $\eta > 2$ is a fixed parameter from Theorem \ref{thm:meta-performance-rams} and $\nu > 0$ is a chosen parameter. This implies that $K_t^{-\frac{1}{\eta}} \leq (T - t)^{-1 - \frac{\nu}{\eta}}$ and hence we have that $C \sum_{t = 1}^T K_t^{-\frac{1}{\eta}} \leq C \sum_{t = 1}^T (T - t + 1)^{-1 - \frac{\nu}{\eta}} \leq C \int_{1}^{T} x^{-1-\frac{\nu}{\eta}} dx \leq C(\eta,|\mathcal{A}|, F, \nu)$ since $\nu / \eta > 0$. Since this is common across all the corollaries, we have that the contribution to regret due to the number of simulated sample paths $K_t$ is a constant (depending on $\eta, \nu, F$ and $|\mathcal{A}|$). The only thing remaining to bound is $(\diamondsuit)$ under different assumptions and problem settings.

    \proof{\underline{Proof of Corollary \ref{cor:cwg-rams-guarantee}}.} 
        From Lemma \ref{lem:expected-compensation-bound}, it follows that for $t \leq \tilde{T} = T - \floor{64 \log (1/ \varepsilon_0) / \varepsilon_{0}^2}$, we have that $\sup_{B_t \geq 0} \partial \mathcal{R}_t(B_t, a_t^{\sf CwG}) \leq C\left( (\log \tau / \tau)^{\frac{1}{2} + \frac{1}{2(1 + \beta)}} + \delta \sqrt{\log \tau / \tau} \right)$ which implies that for $t \leq \Tilde{T}, \Delta_t({\sf CwG}) = C\left( (\log \tau / \tau)^{\frac{1}{2} + \frac{1}{2(1 + \beta)}} + \delta \sqrt{\log \tau / \tau} \right)$. Summing $\Delta_{t}({\sf CwG})$ from $t = 1$ to $t = \tilde{T}$, implies that the regret contribution is at most $C((\log T)^{\frac{1}{2} + \frac{1}{2(1 + \beta)}} T^{\frac{1}{2} - \frac{1}{2(1 + \beta)}} \mathbbm{1}\{\beta > 0\} + \log^2 T \mathbbm{1}\{\beta = 0\} + \delta \sqrt{T \log T})$.
        Since we are considering a variant of {\sf RAMS} which employs a static allocation policy (same as the one deployed in Algorithm \ref{alg:CwG}) for the last $\ceil{64 \log (1 / \varepsilon_0) / \varepsilon_0^2}$, the regret accrued over the last $\ceil{64 \log (1 / \varepsilon_0) / \varepsilon_0^2}$ time steps is upper bounded as $C\sqrt{\log (1/\varepsilon_0)} / \varepsilon_0$. Adding up all the contributions (including due to $(\heartsuit)$), we attain the same regret scaling as in Theorem \ref{thm:general-cwg-upper-bound}. \QEDB
    \endproof
    
    \proof{\underline{Proof of Corollary \ref{cor:nrm-regret-profiles}.}} 
    For each of the Assumptions \ref{ass:nrm-finite-types}, \ref{ass:nrm-infinite-types} and \ref{ass:nrm-semi-infinite-types}, we have that $\sup_{B \geq \bm{0}, a \in \mathcal{A}, \bm{\theta}_{\geq t}} \partial \mathcal{R}_t(B, a; \bm{\theta}_{\geq t}) \leq d$ for all $t \in \{1, \dots, T\}$ since the offline will need a compensation of atmost $r_{\max} = 1$ per resource in the future for accepting or rejecting the request $\theta_t$. Since there are $d$ {\it fixed} resources, the compensation is atmost $d$. Now under different assumptions, we have different algorithms with different values for $\Delta_t(\sf ALG)$. 
    \begin{itemize}
        \item[$(a)$] \underline{\it Under Assumption \ref{ass:nrm-finite-types}.} From (9) in \cite{vera2021bayesian}, we have that for the {\sf Bayes Selector} algorithm described in Algorithm 2 of \cite{vera2021bayesian}, $\sup_{B_t \geq \bm{0}} \partial \mathcal{R}_t(B_t, a_t^{{\sf Bayes Selector}}) \leq d \exp(-c\tau) : = \Delta_t({\sf Bayes Selector})$ for $t \leq T - T_0$ where $\tau = T - t + 1$, and $c, T_0$ are constants which depend only on the distribution $F$. Using the fact that in the last constant $T_0$, the regret accrued is atmost $d T_0$ and $\int_1^T d \exp(-c\tau) d\tau \leq C$, we have that the total regret accrued by {\sf RAMS} under Assumption \ref{ass:nrm-finite-types} is at most a constant $C$ which depends on the parameters $\eta > 2, \nu > 0,$ number of resources $d$ and the distribution $F$.

        \item[$(b)$] \underline{\it Under Assumption \ref{ass:nrm-infinite-types}.} From Lemma 5, 8, 9 and 10 of \cite{jiang2022degeneracy}, for the {\sf Bid Price} algorithm described in Algorithm 3 of \cite{jiang2022degeneracy}, we have that $\sup_{B_t \geq \bm{0}} \partial \mathcal{R}_t(B_t, a_t^{{\sf Bid Price}}) \leq C / \tau := \Delta_t({\sf Bid Price})$ for $t \leq T - T_0$ where $\tau = T - t + 1$ and $C, T_0$ are constants which depend only on the distribution $F$. Using the fact that in the last constant $T_0$, the regret accrued is atmost $d T_0$ and $\int_1^T C/ \tau d\tau \leq C \log T$, we have that the total regret accrued by {\sf RAMS} under Assumption \ref{ass:nrm-infinite-types} is at most $C\log T$ where the constant  depends on the parameters $\eta > 2, \nu > 0,$ number of resources $d$ and the distribution $F$.

        \item[$(c)$] \underline{\it Under Assumption \ref{ass:nrm-semi-infinite-types}.} From Theorem 1 of \cite{jiang2022degeneracy}, for the {\sf Boundary Attracted} algorithm described in Algorithm 2 of \cite{jiang2022degeneracy}, we have that $\sup_{B_t \geq \bm{0}} \partial \mathcal{R}_t(B_t, a_t^{{\sf Boundary Attracted}}) \leq C\log \tau / \tau := \Delta_t({\sf Boundary Attracted})$ for $t \leq T - T_0$ where $\tau = T - t + 1$ and $C, T_0$ are constants which depend only on the distribution $F$. Using the fact that in the last constant $T_0$, the regret accrued is atmost $d T_0$ and $\int_1^T C\log \tau / \tau d\tau \leq C \log^2 T$, we have that the total regret accrued by {\sf RAMS} under Assumption \ref{ass:nrm-semi-infinite-types} is at most $C\log^2 T$ where the constant which depends on the parameters $\eta > 2, \nu > 0,$ number of resources $d$ and the distribution $F$.
    \end{itemize}
    This concludes the proof for all three cases. \QEDB
    \endproof
    
    \proof{\underline{Proof of Corollary \ref{cor:online-matching-constant}.}}
        The proof follows analogously to the proof of Corollary \ref{cor:nrm-regret-profiles} under Assumption \ref{ass:nrm-finite-types}. \QEDB
    \endproof
}

\new{
\section{Relating the order fulfillment problem to the multisecretary problem}
\label{onapp:fulfillment-problem-mapping}

\subsection{Motivating Example}
Let's explore the following example to illuminate our point: Consider two Amazon fulfillment centers, located respectively in Salt Lake City, Utah, and Sacramento, California, as presented in Figure \ref{fig:fulfillment-example}. Both states have a total of over two thousand zip codes, which are spatially clustered and represent distinct demand locations.

The United States sees an estimated total demand volume of around eighty million Amazon packages delivered weekly. Without a precise state-wise breakdown of these deliveries, we can reasonably assume that combined deliveries in California and Utah amount to no more than five million each week. Based on these figures, we calculate a total demand volume ($T$) of $5 \times 10^6$, and a total number of demand locations or types ($D$) of $2 \times 10^3$.

Assuming uniform demand across these locations, our model suggests that at any given time $t$, the probability of receiving a demand request from type $j$ is $D^{-1} \approx T^{-\frac{1}{2}}$, which scales with the total demand volume. This differs from settings studied previously, which considered atomic distributions with a few types and implicitly assumed that the probability of receiving a demand request at a given time was independent of the total demand volume - an assumption inconsistent with the example we have described.

Alternatively, we could consider the setting where infinitely many types exist over a contiguous support. However, this approximation falls short in the presence of natural geographical features like the Sierra Nevada desert which creates gaps, as depicted in Figure \ref{fig:fulfillment-example}. Neither of these previously explored models satisfactorily fit this stylized order fulfillment problem. Instead, what we encounter is a scenario characterized by many types with gaps.demand request at time $t$ is independent of the total demand volume and hence does not align well with the aforementioned example. On the other extreme, one could consider the setting with infinitely many types over a continguous support but clearly such an approximation is wanting in the presence of gaps introduced by natural geographical features like the desert in Nevada as shown in Figure \ref{fig:fulfillment-example}. Therefore neither of the previously studied models are a good fit for this stylized order fulfillment problem. What we have are essentially \emph{many types with gaps}.

\begin{figure}[ht]
    \centering
    \includegraphics[width = 1\linewidth]{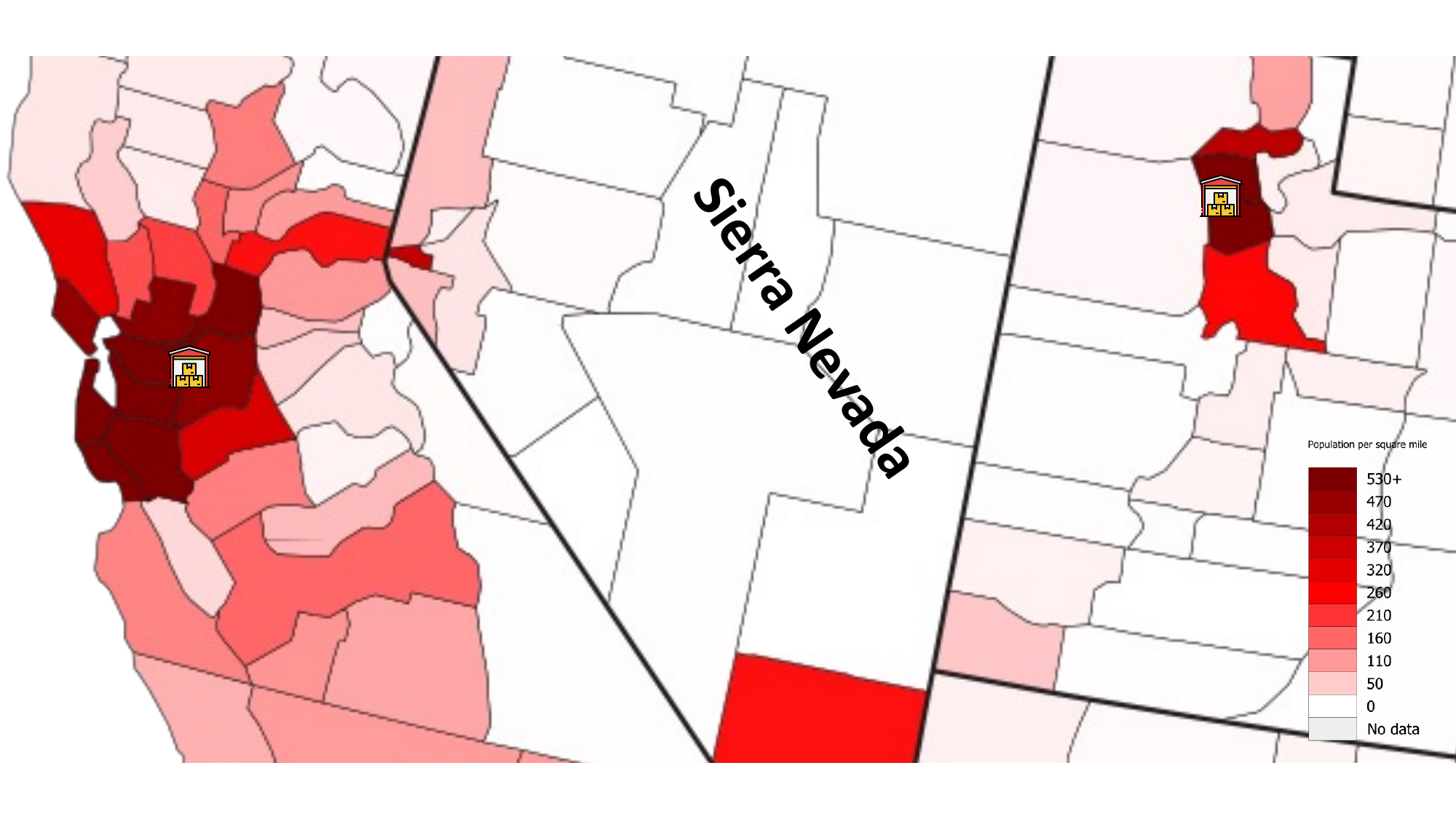}
    \caption{Illustration of spatially distributed demand with two fulfillment centers for the order fulfillment problem 
    }
    \label{fig:fulfillment-example}
\end{figure}

\subsection{Stylized model of order fulfillment}
Inspired by our example illustrated above, we consider a stylized order fulfillment problem with the demand locations being spatially distributed over the unit square $[0,1]^2$ and two fulfillment centers (FCs) denoted as ${\sf FC}A$ and ${\sf FC}B$ with a total inventory in the two warehouses being $T$. The initial inventory in ${\sf FC} A$ and ${\sf FC} B$ is denoted as $I^A_1$ and $I^B_1$ respectively. Now at each time $t$, a request $\xi_t$ arrives given by the coordinates $(x_t,y_t) \in [0,1]^2$ which is drawn from some spatial demand distribution $Q$ with measure $\mu_Q$. Given the inventory levels $I_t^A$ and $I_t^B$ at time $t$, the order fulfillment problem is to decide which fulfillment center to serve the request $\xi_t$ from. The goal is to minimize the total matching distance between the requests and the fulfillment center from which they are served. It is easy to see that this problem can be easily translated into the multisecretary problem. We will illustrate this correspond via an example as shown in Figure \ref{fig:stylized-order-fulfillment-model}.

\begin{figure}
    \centering
    \begin{tikzpicture}
            
            \fill[green!30] (0, 4.8) -- (4.8, 0) -- (0, 0) -- cycle;
            \fill[green!30] (6, 1.2) -- (1.2, 6) -- (6, 6) -- cycle;
            \draw[ultra thick] (0, 0) rectangle (6, 6);
            \node at (1.2,1.2) { \large $\mathfrak{R}_1$};
            \node at (3,3) {\large $\mathfrak{R}_2$};
            \node at (4.8, 4.8) {\large $\mathfrak{R}_3$};
            \warehouse{-0.3,-0.3}{blue};
            \warehouse{5.7,5.7}{blue};
            \node at (-0.7,0) {$\textsf{FC}A$};
            \node at (6.7,6) {$\textsf{FC}B$};
        \end{tikzpicture}
        \caption{Stylized example for order fulfillment with no demand from region $\mathfrak{R}_2$}
    \label{fig:stylized-order-fulfillment-model}
\end{figure}

In the stylized example illustrated in Figure \ref{fig:stylized-order-fulfillment-model}, we assume that the demand locations are uniformly distributed in regions $\mathfrak{R}_1$ and $\mathfrak{R}_3$ with no demand in region $\mathfrak{R}_2$. The two fulfillment centers ${\sf FC}A$ and ${\sf FC}B$ are located at $(0,0)$ and $(1,1)$ respectively. Let $d_A((x,y)) = |x| + |y|$ and $d_{B}((x,y)) = |1 - x| + |1 - y|$ denote the (Manhattan) distance from the fulfillment centers ${\sf FC}A$ and ${\sf FC} B$ respectively. The hindsight optimal problem for the order fulfillment problem is the following integer program. \begin{align}
    \label{eq:order-fulfillment-maximization}
    \min_{z_t} & \quad \sum_{t = 1}^{T} d_{A}((x_{t}, y_t)) z_t + d_B((x_t, y_t)) (1 - z_t) \\
    \nonumber
    \text{s.t.} & \quad \sum_{t = 1}^T z_t = I_1^A \\
    \nonumber
    &\quad z_t \in \{0,1\}, \quad  \forall t
\end{align}
The objective in \eqref{eq:order-fulfillment-maximization} can be equivalently written as $\sum_{t = 1}^T \left(d_A((x_t, y_t)) - d_B((x_t, y_t)) \right) z_t + d_B((x_t, y_t))$ and hence we can cast the minimization problem into the following maximization problem.
\begin{align}
    \label{eq:multisecretary-maximization}
    \max_{z_t} & \quad \sum_{t = 1}^T \left( d_{B}((x_t, y_t)) - d_{A}((x_t, y_t)) \right) z_t \\ 
    \nonumber
    \text{s.t.} & \quad \sum_{t = 1}^T z_t = I^A_1 \\
    \nonumber
    & \quad z_t \in \{0,1\}, \quad \forall t
\end{align}

The optimization problem in \eqref{eq:multisecretary-maximization} is the multisecretary problem with reward $\tilde{r}((x_t, y_t)) = d_B((x_t, y_t)) - d_A((x_t, y_t))$ we consider in this work with appropriate scaling. Observe that $\tilde{r}((x_t, y_t)) \in [-2,2]$, therefore we can scale the reward and define the types as $\theta_t = (\tilde{r}((x_t, y_t)) + 2) / 4 \in [0,1]$. We can translate the spatial demand distribution $Q$ into the distribution over the types $\theta_t$ as follows.
\begin{align*}
    \mathbb{P}\left(\theta_t \leq z \right) &= \mathbb{P}\left(\tilde{r}((x,y)) \leq 4 z - 2 \right) \\
    &= \mathbb{P}\left(d_{B}((x,y)) - d_A((x,y)) \leq 4z - 2 \right) \\
    &= \mathbb{P}\left( (1 - x  + 1 - y) - (x + y) \leq 4 z - 2 \right) \\
    &= \mathbb{P}\left( x + y \geq 2 - 2z \right)
\end{align*}
Using the fact that demand locations are uniformly distribution in regions $\mathfrak{R}_1$ and $\mathfrak{R}_3$, we have that 
\begin{align*}
    \mathbb{P}\left( \theta_t \leq z \right) = \begin{cases}
        \frac{25}{8}z^2 & \quad z \in [0,\frac{2}{5}], \\
        \frac{1}{2} & \quad z \in [\frac{2}{5}, \frac{3}{5}], \\
        1 - \frac{25}{8}(1 - z)^2 & \quad z \in [\frac{3}{5}, 1].
    \end{cases}
\end{align*}
Note that this is a $(\beta = 0, \varepsilon_0 = \frac{1}{2})$-clustered distribution with a gap interval $[\frac{2}{5}, \frac{3}{5}]$. Note that the gap in the demand location for the order fulfillment translates into a gap in the type distribution for the multisecretary problem. Moreover for simplicity we assume that the demand locations in regions $\mathfrak{R}_1$ and $\mathfrak{R}_3$ are distributed over a contiguous support but we can further discretize these regions into many small types which are clustered close to each other. This captures the more realistic setting where various zipcodes are spatially close to each other. In the context of the multisecretary problem, this is captured via the $(\beta, \varepsilon_{0}, \delta)$-clustered distribution (recall Definition \ref{def:beta-epsilon-delta-clustered}).
}

\new{
\section{Some details on the $(\beta, \varepsilon_0, \delta)$-clustered distributions}
\label{onapp:details-beta-varepsilon-delta-clustered}
As mentioned in Section \ref{sec:beta-epsilon-distribution}, there is some flexibility is how we may model a distribution or define clusters. Additionally, the parameter $\delta$ allows us to model distributions with many {\it small} types. 
In this section, we will discuss how the same distribution can have different characterizations due to different clustering and choice of parameters $\beta, \varepsilon_0$ and $\delta$. We will discuss this using two examples. For each of the examples, we will discuss two different possible clusterings and their impact on the regret guarantees.
The two examples we will consider will be atomic distributions and let $F$ be the continuous limit of those atomic distributions described below.
\begin{align*}
    F(x) =  \begin{cases}
            -8(1/4 - x)^{2} + 1/2, &\quad \quad 0 \leq x \leq 1/4 \\
            1/2, & \quad \quad 1/4 \leq x \leq 3/4 \\
            8(x - 3/4)^2 + 1/2, & \quad \quad 3/4 \leq x \leq 1.
    \end{cases}
\end{align*}
Note that it can be easily verified that the distribution $G$ above is a $(\beta = 1 , \varepsilon_0 = 1/2, \delta = 0)$-clustered distribution. Let us consider two other atomic distributions with probability mass functions denoted as $p_1$ and $p_2$ respectively and defined using the parameters $\eta_1$ and $\eta_2$ as follows,
\begin{align*}
    p_1(1/4 - k\eta_1) = p_1(3/4 + k \eta_1) = 16\eta_1^2 / (4 \eta_1 + 1) , \forall k \in \{0, 1, \dots, 1/4\eta_1\} \\
    p_2(1/4 - k\eta_2) = p_2(3/4 + k \eta_2) = 16\eta_2^2 / (4 \eta_2 + 1) , \forall k \in \{0, 1, \dots, 1/4\eta_2\}
\end{align*}
For $\eta_1 = 1/24$, we get the distribution with probability mass function $p_1$ is supported on twelve points and is an example of distribution with a few types (refer to center figure in Figure \ref{fig:discretization-of-beta-clustered}) with the minimum probability mass being $1/42$. For $\eta_2 = 1/2400$, we get the distribution with probability mass function $p_2$ is supported on twelve hundred points and can be considered an example of many small points since the number of types are large (1200) and each type has a small probability mass (at most $2 \times 10^{-3}$) (refer to the right figure in Figure \ref{fig:discretization-of-beta-clustered}). Note that as $\eta_1, \eta_2 \to 0$, we have that $p_1, p_2 \to f$. 

There are two natural ways that the distribution $p_1$ can be modelled as a $(\beta, \varepsilon_0, \delta)-$clustered distribution. First way is as a distribution with a few types (as modelled in Example \ref{ex:discret-distributions}), where we have twelve mass clusters $H = \cup_{k = 0}^{k = 5} \{k/24, (19 + k) / 24\}$ corresponding to the twelve points on the which the distribution is supported with eleven gap intervals $G = \left(\cup_{k = 0}^{k = 4} (k/24, (k+1)/24) \cup ((19 + k)/24, (20 + k) / 24) \right) \cup (5/24,19/24)$. We can easily verify that the distribution $p_1(x)$ satisfies the conditions in Definition \ref{def:beta-epsilon-delta-clustered} with $\beta = 0, \varepsilon_0 = \min_{\{x : p_1(x) > 0\}} p_1(x) = 1/42$ and $\delta = 0$. The second way to model this distribution is by having only two mass clusters $H_1 = [0,1/4]$ and $H_2 = [3/4,1]$ with one gap interval $G = (1/4,3/4)$. By considering only two clusters, we have that $\varepsilon_0 = 1/2$ since the total probability mass in both the clustered is $1/2$ each. It is easy to see that for any choice of $\beta \in [0,\infty)$, to satisfy condition (a) in Definition \ref{def:beta-epsilon-delta-clustered}, we must choose $\delta = \eta_1 > 0$. While both ways are valid in terms of modelling the distribution, the theoretical guarantees implied by the two different characterizations of the same distribution lead to two different regret scalings. Under the first way where $p_1$ is modelled as a $(\beta = 0, \varepsilon_0 = 1/42, \delta = 0)-$clustered, we get constant regret scaling, while under the second way where $p_1$ is modelled as a $(\beta = 0, \varepsilon_0 = 1/2, \delta = 1/24)-$clustered, we get that the regret will scale as $\tilde{O}(\sqrt{T})$. Note that these regret scalings not only follow from the bounds in Theorem \ref{thm:general-cwg-upper-bound} but also due to the fact that {\sf CwG} algorithm in Algorithm \ref{alg:CwG} will operate differently under the two different characterizations of the same distribution $p_1$ since the gaps are defined differently under the two different characterizations.

\begin{figure}[ht]
    \centering
        \pgfplotsset{width = 0.31\linewidth, ylabel near ticks}
        \begin{tikzpicture}
            \begin{axis}[
                axis lines = left,
                xlabel = \(x\),
                ylabel = {\(f(x)\)},
            ]
            \addplot[domain = 0:0.25, color = blue, thick]{4 - 16*x};
            \addplot[domain = 0.75:1, color = blue, thick, forget plot]{16*x - 12};
            \end{axis}
        \end{tikzpicture}
        \begin{tikzpicture}
  \begin{axis}[
      axis lines= left,
      xlabel={$x$},
      ylabel={\(p_1(x)\)},
      ylabel style={
      },
      xtick={0,0.2,0.4,0.6,0.8,1},
      xmin=0, xmax=1,
      ymin=0, ymax=0.25,
    ]
    \addplot [ycomb, red, thick, mark=*] table [x={n}, y={xn}] {discrete.dat};
  \end{axis}
\end{tikzpicture}
        \begin{tikzpicture}
  \begin{axis}[
      axis lines= left,
      xlabel={$x$},
      ylabel={\(p_2(x)\)},
      ylabel style={
      },
      xtick={0,0.2,0.4,0.6,0.8,1},
      xmin=0, xmax=1,
      ymin=0, ymax=0.002,
    ]
    \addplot [ycomb, teal, mark = x] table [x={n}, y={xn}] {manysmall.dat};
  \end{axis}
\end{tikzpicture}
    \caption{ (Left) PDF $f_{\beta}$ of distribution $F$ (Center) a few types (Right) many small types}
    \label{fig:discretization-of-beta-clustered}
\end{figure}
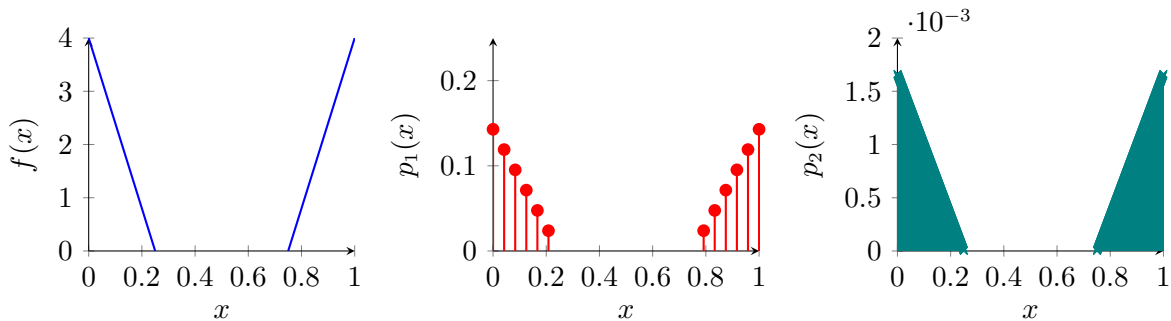

Coming to the distribution $p_2$, again there are two ways that the distribution $p_2$ can be modelled as a $(\beta, \varepsilon_0, \delta)-$clustered distribution. Since strictly speaking, $p_2$ is an atomic distribution albeit with many types, we can model is similar to how we modelled an atomic distribution with a few types. Building on that, we would have that 1200 mass clusters $H = \cup_{k = 0}^{k = 598} \{k/2400, (1801 + k) / 2400\}$ with 1199 gap intervals $G = \left(\cup_{k = 0}^{598}(k/2400,(k + 1) / 2400) \cup ((1801 + k) /2400, (1802 + k) / 2400) \right) \cup (599/2400, 1801/2400)$. We can easily verify that the distribution $p_2(x)$ satisfies the conditions in Definition \ref{def:beta-epsilon-delta-clustered} with $\beta = 0, \varepsilon_0 = \min_{\{x: p_2(x) > 0\}} p_2(x) = 1/360600$ and $\delta = 0$. The second way to model this distribution is having only two mass clusters $H_1 = [0,1/4]$ and $H_2 = [3/4,1]$ with one gap interval $G = (1/4, 3/4)$. By considering only two clusters, we have that $\varepsilon_0 = 1/2$. It is easy to verify that for $\beta = 1$ and $\delta = \eta_2$, we satisfy the condition (a) in Definition \ref{def:beta-epsilon-delta-clustered}. Note that under the first way, we have that $\varepsilon_0$ is very small and for most reasonable and practical values of the time horizon $T$, we may have that $1/\varepsilon_0 \sim T$ and hence the theoretical guarantees implied by Theorem \ref{thm:general-cwg-upper-bound} may be vacuous. On the other hand, in the second characterization as $(\beta = 1, \varepsilon_0 = 1/2, \delta = \eta_2)$-clustered distribution, we have that $\delta \sim 1/\sqrt{T}$ for reasonable values of $T$ and implied regret scaling is $\tilde{\mathcal{O}}(T^{\frac{1}{4}})$ (sublinear regret). Note that the {\sf CwG} algorithm (Algorithm \ref{alg:CwG}) operates differently under the two different characterizations of the same distribution $p_2$.
}

\end{document}